\documentclass[11pt]{amsart}
\usepackage{amsmath,amssymb,amsfonts,url}
\usepackage{color}
\usepackage{appendix}
\numberwithin{equation}{section}

\newtheorem{thm}{Theorem}[section]
\newtheorem{lem}{Lemma}[section]
\newtheorem{rem}{Remark}[section]
\newtheorem{prop}{Proposition}[section]
\newtheorem{cor}{Corollary}[section]

\usepackage{palatino,mathpazo}

\usepackage{amsthm,a0size,bm,indentfirst,wasysym}

\usepackage{amsmath}

\allowdisplaybreaks[4]

\begin{document}
	\title[Uniqueness of blowup solution]{Asymptotic Analysis and Uniqueness of blowup solutions of non-quantized singular mean field equations}
	
	%\author{Daniele Bartolucci, Wen Yang, Lei Zhang}\footnote{L.Z. is partially supported by Simon's collaboration grant 584918}
	
	\author[D. Bartolucci]{Daniele Bartolucci}
	\address{Daniele Bartolucci, Department of Mathematics, University of Rome "Tor Vergata", Via della ricerca scientifica 1, 00133 Roma,
		Italy}
	\email{bartoluc@mat.uniroma2.it}
	
	\author[W. Yang]{Wen Yang}
	\address{Wen Yang, Department of Mathematics, Faculty of Science and Technology, University of Macau, Taipa, Macau}
	\email{wenyang@um.edu.mo}
	
	\author[L. Zhang]{Lei Zhang}
	\address{Lei Zhang, Department of Mathematics, University of Florida, 1400 Stadium Rd, Gainesville, FL 32611, USA}
	\email{leizhang@ufl.edu}

	\date{\today}
	%\thanks{$^{\sharp}$D.B. is partially supported by the MIUR Excellence Department Project MatMod@TOV awarded to the Department of Mathematics, University of Rome "Tor Vergata" and by PRIN project 2022, ERC PE1\_11, "{\em Variational and Analytical aspects of Geometric PDEs}" }
	%%%%%%%%%%%%%%%%%%%%%%%%%%%%%%%%%%%%%%%%%%%%%
	\begin{abstract}
		For singular mean field equations defined on a compact Riemann surface, we prove the uniqueness of bubbling solutions as far as blowup points are either regular points or non-quantized singular sources. In particular the uniqueness result covers the most general case extending or improving all previous works of Bartolucci-Jevnikar-Lee-Yang \cite{bart-4,bart-4-2} and Wu-Zhang \cite{wu-zhang-ccm}. For example, unlike previous results, we drop the assumption of singular sources being critical points of a suitably defined Kirchoff-Routh type functional. Our argument is based on refined estimates, robust and flexible enough to be applied to a wide range of problems requiring a delicate blowup analysis. In particular we come up with several new estimates of independent interest about the concentration phenomenon for Liouville-type equations.
	\end{abstract}
	%	{\bf Mathematics Subject Classification (2007):
		%35J60, 35B45, 53C21}
	
	%  {\bf Keywords:}
	% Liouville equation, Blowup analysis.
	
	%%%%%%%%%%%%%%%%%%%%%%%%%%%%%%%%%%%%%%%%%%%%%
	
	\maketitle
	
	\section{Introduction}
	
	The main goal of this article is to study the uniqueness property of the following mean field equation with singular sources:
	\begin{equation}\label{m-equ}
		\Delta_g \nu+\rho\bigg(\frac{he^\nu}{\int_M h e^\nu{\rm d}\mu}-\frac{1}{vol_g(M)}\bigg)=\sum_{j=1}^N 4\pi \alpha_j \left(\delta_{q_j}-\frac{1}{vol_g(M)}\right) \quad {\rm in} \ \; M,
	\end{equation}
	where $(M,g)$ is a Riemann surface with metric $g$, $\Delta_g$ is the Laplace-Beltrami operator ($-\Delta_g\ge 0$), $h$ is a positive $C^5$ function on $M$, $q_1,\cdots,q_N$ are distinct points on $M$, $\rho>0,\alpha_j>-1$ are constants, $\delta_{q_j}$ is the Dirac measure at $q_j\in M$, $vol_g(M)$ is the volume of $(M,g)$. Throughout the article we assume $vol_g(M)=1$ for the sake of simplicity. Equation (\ref{m-equ}) arises in many branches of Mathematics and Physics, such as conformal geometry (\cite{B1,Zchen-lin,fang-lai,KW,troy,wei-zhang-pacific}), Electroweak and Self-Dual Chern-Simons vortices (\cite{ambjorn, spruck-yang,taran-1,taran-2,y-yang}), the statistical mechanics description of turbulent Euler flows,  plasmas and self-gravitating systems (\cite{bart-5,caglioti-2,SO,TuY,wolan}), Cosmic Strings (\cite{BT-2,pot}), the theory of hyperelliptic
	curves and modular forms (\cite{chai}) and CMC immersions in hyperbolic 3-manifolds (\cite{taran-4}). Among many other things which we cannot list here, these were some of the motivations for the huge effort made in the study of (\ref{m-equ}), including existence (\cite{BdM1,BdM2,BdMM,BLin2,BMal,luca-b,dAWZ,chan-fu-lin,chen-lin-sharp,chen-lin-deg,chen-lin,chen-lin-deg-2,DKM,DJ,EGP,lin-yan-cs,lin-yan-cs,Mal1,Mal2,mal-ruiz,nolasco-taran}), uniqueness (\cite{BGJM,bart-4,bart-4-2,bart-4-3,bj-lin,bart-lin,BLin3,BLT,chang-chen-lin,GM1,lee-lin-jfa,Lin1,lin-yan-uniq,suzuki,wu-zhang-ccm}), concentration-compactness and bubbling behavior (\cite{BM3,BT,BM,gu-zhang-1,gu-zhang-2,gluck,kuo-lin,li-cmp,ls,T3,taran-3,wei-zhang-19,wei-zhang-22, wei-zhang-jems, zhang1,zhang2}) just to quote a few. In particular, started in 
	\cite{chen-lin-wang,lin-Lwang}, a series of groundbreaking results has been recently pushed forward in the case of the flat two torus with singular sources, see  \cite{Zchen-lin-1,lin22} and references quoted therein.
	
	\medskip

	To write the main equation in an equivalent form, we invoke the standard Green's function $G_M(z,p)$, defined to be the unique (\cite{Aub}) solution of:
	\begin{equation*}
		%\label{gf}
		\begin{cases}
			-\Delta_g G_M(x,p)=\delta_p-1\quad {\rm in}\ \; M
			\\
			\\
			\int_{M}G_M(x,p){\rm d}\mu=0.
		\end{cases}
	\end{equation*}
	Then
	in local isothermal coordinates $z=z(x)$ centered at $p$, $0=z(p)$,
	$G(z,0)=G_M(x(z),p)$ can be written as
	follows
	$$G(z,0)=-\frac 1{2\pi }\log |z|+R(z,0),$$
	where $R(z,0)$ denotes the regular part in local coordinates. By using $G_M(x,p)$ we write (\ref{m-equ}) as follows,
	\begin{equation}\label{r-equ}
		\Delta_g  {\rm w}+\rho\bigg(\frac{He^{\rm w}}{\int_M H e^{\rm w}{\rm d}\mu}-1\bigg)=0 \quad {\rm in}\ \; M,
	\end{equation}
	where
	\begin{equation}\label{r-sol}
		{\rm w}(x)=\nu(x)+4\pi \sum_{j=1}^N \alpha_j G_M(x,q_j),
	\end{equation}
	and
	\begin{equation*}
		%\label{H1}
		H(x)=h(x)\prod_{j=1}^N e^{-4\pi\alpha_j G_M(x,q_j)}.
	\end{equation*}
	Note that in a local coordinates system near $q_j$, $0=z(q_j)$,
	\begin{equation*}
		%\label{H2}
		H(z)=h_j(z)|z|^{2\alpha_j},\quad |z|\ll 1,\quad 1\leq j\leq N,
	\end{equation*}
	for some $h_j(z)>0$.
	
	We say that $\{\nu_k\}$ is a sequence of bubbling solutions of (\ref{m-equ}) if the $L^{\infty}$ norm of the corresponding ${\rm w}_k$ defined by (\ref{r-sol}) tends to infinity as $k$ goes to infinity. In this situation it is well known (\cite{BM3,BT,li-cmp}) that ${\rm w}_k$ blows up at a finite number of points. To say that the set of blowup points is $\{p_1,\cdots,p_m\}$ means that there exist $m$ sequences of points $p_{k,1},\cdots,p_{k,m}$ such that, possibly along a subsequence, ${\rm w}_k(p_{k,j})\to +\infty$ and $p_{k,j}\to p_j (j=1,\cdots,m)$ as $k\to +\infty$. If none of the blowup points is a singular source,
	the uniqueness of the bubbling solutions has been proved in  \cite{bart-4}. Therefore we consider here the case where at least one blow up point is a singular source.  By assumption the strength of the singular source at $p$ is $4\pi \alpha_p$, where if $p$ is a regular point, the corresponding $\alpha_p=0$. Thus we use $\alpha_1,\cdots,\alpha_m$ to denote the strengths at $p_1,\cdots,p_m$, respectively and we let the first $\tau\geq 1$ of them be singular blowup points:
	\begin{equation}\label{largest-s}
		\alpha_i>-1,\quad \alpha_i\not \in \mathbb N, \quad 1\le i\le \tau,\quad \alpha_{\tau+1}=\cdots=\alpha_m=0,
	\end{equation}
	where $\mathbb N$ is the set of natural numbers including $n=0$. We will say that the singular source located at $p$ is \emph{non-quantized} if $\alpha_p$ is not a positive integer. Here and in the rest of this work we assume that all the singular sources, as far as they happen to be blowup points, are non-quantized and let $\alpha_M$ denote the maximum of $\{\alpha_1,...,\alpha_m\}$.\\ 
	There is a variational structure of equation (\ref{r-equ}) which is in fact the Euler-Lagrange equation of the functional:
	$$I_{\rho}({\rm w})=\frac 12 \int_M |\nabla {\rm w}|^2+\rho\int_M {\rm w}-\rho \log \int_M He^{\rm w},\quad {\rm w}\in H^1(M).$$
	Since adding a constant to any solution of (\ref{r-equ}) still yields a solution, there is no loss of generality in defining $I_{\rho}$ on the subspace of $H^1(M)$ functions with vanishing mean. A complete discussion about the variational structure of (\ref{r-equ}) can be found in \cite{mal-ruiz}.
	
	\smallskip
	
	To state the main results we adopt the following notations:
	\begin{equation*}
		%\label{G_j*}
		G_j^*(x)=8\pi (1+\alpha_j)R_M(x,p_j)+8\pi \sum_{l,l\neq j}(1+\alpha_l)G_M(x,p_l),
	\end{equation*}
	where $R_M(\cdot, p_j)$ is the regular part of $G_M(\cdot, p_j)$,
	\begin{equation}\label{L-p}
		L(\mathbf{p})=\sum_{j\in I_1}[\Delta \log h(p_j)+\rho_*-N^*-2K(p_j)]h(p_j)^{\frac{1}{1+\alpha_M}}e^{\frac{G_j^*(p_j)}{1+\alpha_M}},
	\end{equation}
    where $K(\cdot)$ denotes the Gaussian curvature of the ambient manifold, 
	\begin{equation}
		\label{a-note}
		\begin{cases} 
			\alpha_M=\max_{i}\alpha_i,\quad
			I_1=\{i\in \{1,\cdots,m\};\quad \alpha_i=\alpha_M\}, \\
			\\
			\rho_*=8\pi\sum_{j=1}^m(1+\alpha_j),\quad  N^*=4\pi\sum_{j=1}^m\alpha_j. 
		\end{cases}
	\end{equation}
	In other words, $L(\mathbf{p})$ takes the sum on those indexes whose $\alpha_j$ equals the largest among them.
	As mentioned above, whenever some blowup points are non-quantized singular sources while some others are regular, we assume without loss of generality that $p_1,\cdots,p_{\tau}$ are singular sources and $p_{\tau+1},\cdots,p_m$ are regular points. In this case
	for $(x_{\tau+1},\cdots,x_m)\in M\times\cdots \times M$ we define
	\begin{equation*}
		%\label{f*}
		\begin{aligned}
			f^*(x_{\tau+1},\cdots,x_m)
			=&\sum_{j=\tau+1}^m\big[\log h(x_j)+4\pi R_M(x_j,x_j)\big]\\
			&+4\pi \sum_{l\neq j}^{\tau +1,\cdots,m}G_M(x_l,x_j)\\
			&+\sum_{j=\tau+1}^m\sum_{i=1}^{\tau}8\pi(1+\alpha_i)G_M(x_j,p_i). 
		\end{aligned}
	\end{equation*}

	It is well known (\cite{chen-lin-sharp}) that  $(p_{\tau+1},\cdots,p_m)$ is a critical point of $f^*$. In this article we consider all possible cases of combinations of blowup points and split them in two classes:
	$$\mbox{\bf Class One:}\quad \alpha_M>0,\quad \mbox{\bf Class Two:}\quad \alpha_M\le 0.$$
	Obviously in \emph{Class One} the set of blowup points contains at least one  positive singular source 
	while in \emph{Class two}  either the blowup set is a combination of regular blowup points ($\alpha_M=0$) and negative singular sources, or it only contains negative singular sources $(\alpha_M<0)$. Our first result concerns \emph{Class One}.

	\begin{thm}\label{main-theorem-2}
		Let $\nu_k^{(1)}$ and $\nu_k^{(2)}$ be two sequences of bubbling solutions of {\rm (\ref{m-equ})}  with the same $\rho_k$: $\rho_k^{(1)}=\rho_k=\rho_k^{(2)}$ and the same blowup set: $\{p_1,\cdots,p_m\}$. Suppose $(\alpha_1,\cdots,\alpha_N)$ satisfies {\rm (\ref{largest-s})}, $\alpha_M>0$, $L(\mathbf{p})\neq 0$ and, as far as $\tau<m$,  $\det \big(D^2f^*(p_{\tau+1},\cdots,p_m)\big)\neq 0$. Then there exists $n_0>1$ such that
		$\nu_k^{(1)}=\nu_k^{(2)}$ for all $k\ge n_0$.
	\end{thm}
	Here $D^2f^*$ denotes for the Hessian tensor field on $M$ and we observe that, when evaluating $\det(D^2f^*)$ at $p_j$ for $j=\tau+1,\cdots,m$, we use $\phi_j$ to denote the conformal factor that satisfies 
	\begin{equation*}
		%\label{Def-phi-j}
		\Delta_g=e^{-\phi_j}\Delta,\,\, \phi_j(0)=|\nabla \phi_j(0)|=0,\,\,e^{-\phi_j}\Delta \phi_j=-2K,
	\end{equation*}
	for $j=\tau+1,\cdots,m$ and $h_j$ is understood as $h_je^{\phi_j}$ in a small neighborhood of $p_j$ ($j=\tau+1,\cdots,m$).
	Therefore, if $\tau=m$, which means all blowup points are singular sources, then as far as $\alpha_M>0$ we see that 
	the unique relevant assumption is
	$L(\mathbf{p})\neq 0$. 
	
	To state the result about \emph{Class two}, we first observe that the set of blowup points consists only of regular points and negative sources. We shall introduce new quantities. We use $B(q,r)$ to denote the geodesic ball of radius $r$ centered at $q\in M$, while $\Omega(q,r)$ denotes the pre-image of the Euclidean ball of radius $r$, $B_r(q)\subset \mathbb R^2$, in a suitably defined isothermal coordinate system. If $m-\tau\ge 2$ we fix a family of open, mutually disjoint sets $M_j$, the closure of their union being the whole $M$. Then we set
	$$D(\mathbf{p})=\lim_{r\to 0}\sum_{j=1}^mh(q_j)e^{G_j^*(q_j)}\bigg (
	\int_{M_j\setminus \Omega(q_j,r_j)}
	e^{\Phi_j(x,\mathbf{q})}d\mu(x)-\frac{\pi}{1+\alpha_j}r_j^{-2-2\alpha_j}\bigg ),$$
	where $M_1=M$ if $m=1$, 
    \[r_j=\left\{\begin{array}{ll}r \bigg (8h(q_j)e^{G_j^*(q_j)}\bigg)^{1/2},\quad \tau<j\le m,\\
    r,\quad 1\le j\le \tau,
    \end{array}
    \right.\] and
	\begin{align*}\Phi_j(x,\mathbf{q})=&\sum_{l=1}^m8\pi(1+\alpha_l)G(x,q_l)-G_j^*(q_j)\\
		&+\log h_j(x)-\log h_j(q_j)+4\pi\alpha_j(R_M (x,q_j)-G_M(x,q_j)),
	\end{align*}
	where we remark that in local isothermal coordinates centered at $q_j$, $0=z(q_j)$, we have $4\pi\alpha_j(R(z,0)-G(z,0))=2\alpha_j\log|z|$. It is important to note that in this case, regular blow-up points are critical points of 
    \begin{equation}\label{new-f}G_j^*(x)-G_j^*(q_j)
		+\log h_j(x)-\log h_j(0).
        \end{equation}
     For the integrability of $\Phi_j$ we make this assumption on singular blow-up points. Thus for all blow-up points we require that,
   \begin{equation}\label{neg-crit}
	\nabla \bigg ( G_j^*(x)
	+\log h_j(x)\bigg )|_{x=q_j}=0,\quad 1\le j\le m.
\end{equation}
    Concerning \eqref{neg-crit}, let us remark that in this part we are just concerned with \emph{Class two} solutions, whence $\alpha_M\le 0$ and then, if a singular source exists with strength $\alpha_j$, then it satisfies $-1<\alpha_j<0$.

	\begin{thm}\label{mainly-case-2}
		Under the same assumptions of Theorem \ref{main-theorem-2}, except that here $\alpha_M\le 0$ and (\ref{neg-crit}) holds for all blowup points, then in any one of the following cases:
		\begin{enumerate}
			\item $\alpha_M=0$,\, $L(\mathbf{p})\neq 0$,\, $\det \big(D^2f^*(p_{\tau+1},\cdots,p_m)\big)\neq 0$
			\item $\alpha_M=0$, $L(\mathbf{p})=0$,  $D(\mathbf{p})\neq 0$, $\det \big(D^2f^*(p_{\tau+1},\cdots,p_m)\big)\neq 0$
			\item $\alpha_M<0$, $D(\mathbf{p})\neq 0$
		\end{enumerate} the same uniqueness property as in Theorem \ref{main-theorem-2} holds true.
	\end{thm}
\begin{rem}\label{about-D}
The limit in the definition of $D(\mathbf{p})$ exists because $\Omega(q,r)$ is the pre-image of a ball, a non-integrable term in the neighborhood of $q$ cancels out in the evaluation as $r\to 0$.
\end{rem}

%\textcolor{red}{Thus the quantity D(q) is defined only in the sense of the principal value, and since the non integrable part changes sign, by taking a different cut we could obtain for D(q) any real number. As a consequence I don't see any chance to treat the uniqueness for this part as we did in \cite{bart-4}, where we use in a crucial way not only that the quantity D(q) is well defined, but in fact that it determines the expansion of $\rho_k$ as far as L(p)=0. For example the sign of D(q) determines the sign of $\rho_k-8\pi m$, which is clearly not true in this case. If this were not enough, exactly in the discussion of this point, what we say in the proof of Theorem \ref{mainly-case-2} is that we argue exactly as in \cite{bart-4}. In conclusion: we propose a proof, without providing the details, which we claim that can be done by arguing as in \cite{bart-4}, in a situation where we miss exactly one of the crucial assumptions used in \cite{bart-4}. It looks a bit too much, isn't it?\\ I see three possibilities to fix the problem at once, both here and in the non degeneracy file:\\ (1) add the assumptions needed to run the argument in \cite{bart-4} (this is not my favorite, the assumptions look really bad for many reasons and are very bad for the work);\\ 
%(2) remove the bad negative case $(-1,-1/2]$ and remove the application about the Onsager theory in the non degeneracy file. Treat the 
%bad case in a new work which will include the application;\\

    \begin{rem}\label{merits} Theorems \ref{main-theorem-2} and \ref{mainly-case-2} extend or improve previous results \cite{bart-4,bart-4-2,wu-zhang-ccm} significantly in at least three aspects:
	\begin{enumerate}
		\item They include all possible blowup points combinations: blowup points can be either positive or negative singular sources or both, together with regular points as well;
		\item They do not anymore require positive singular blowup points to be critical points of certain Kirchoff-Routh type functionals. In previous work \cite{bart-4-2,wu-zhang-ccm}, for technical reasons, this assumption was made providing a relevant simplification of the proof. On the one hand, the removal of this hypothesis is crucial for applications and also rather surprising from the physical point of view, as long-lived vortex structures are known to concentrate at critical points of those finite-dimensional functionals (\cite{caglioti-2}). On the other side, uniqueness becomes much harder to prove whenever we drop this assumption;   
		\item Let $\alpha_M$ be the largest index among positive singular sources. It was assumed in \cite{wu-zhang-ccm} that $0<\alpha_M<1$, an hypothesis that we drop here as well. Also, \cite{bart-4-2} was concerned with just one positive singular source.  
		
	\end{enumerate}
	With the removal of restrictions of previous results, the generality of the main theorems is a handy tool suitable for a number of applications. However, as we will shortly discuss below, the proof requires new ideas and several new estimates of independent interest about the concentration phenomenon.
    \end{rem}

	Here we also make a striking comparison about the quantities we use. First it is easy to see that $L(\mathbf{p})$ is local in nature because it only depends on information at certain blowup points, while $D(\mathbf{p})$ is involved with global integration. Both quantities have played important roles in various contexts. For example, $L(\mathbf{p})$ has been critical in Chen-Lin's degree counting program for Liouville equations \cite{chen-lin-sharp,chen-lin-deg}, $D(\mathbf{p})$ plays a crucial role in
	\cite{chang-chen-lin, chen-lin-wang} and in \cite{bart-lin, BLin3} for regular and singular Liouville equations (see in particular \cite{lin22} and references therein) in \cite{huang-zhang-cvpde, lin-zhang-jfa} for Liouville systems and in \cite{lin-yan-cs} for Chern-Simon's equation. Since Theorem \ref{main-theorem-2} and Theorem \ref{mainly-case-2} include all possible combinations of singular and regular blow up points, it is worth to explain the roles of $L(\mathbf{p})$ and $D(\mathbf{p})$ in an unified manner:
	When $\alpha_M>0$, $L(\mathbf{p})$ plays the leading role. When $\alpha_M=0$, $L(\mathbf{p})$ is still more important than $D(\mathbf{p})$. But when $L(\mathbf{p})=0$, $D(\mathbf{p})$ plays the leading role together with the non-degeneracy assumption on $f^*$. The interplay of these two quantities reveals the highly nontrivial nature of the uniqueness problem. Again the case of the flat two torus plays a special role, see \cite{lin-wang-cga}. Let us also remark that around negative singular sources, locally defined quantities play no role, in the sense that from Theorem \ref{mainly-case-2} we infer that only $D(\mathbf{p})$ matters as far as $\alpha_M<0$ and $D(\mathbf{p})\neq 0$.

	One can also obtain similar conclusions for the corresponding Dirichlet problem.
	Let $\Omega$ be an open and bounded domain in $\mathbb{R}^2$ with regular boundary $\partial\Omega\in C^2$, $\nu$ be a solution of
	\begin{equation}\label{equ-flat}
		\begin{cases}
			\Delta \nu+\rho \frac{he^\nu}{\int_{\Omega} h e^\nu{\rm d}x}=\sum_{j=1}^N 4\pi \alpha_j \delta_{q_j}  &{\rm in} \;\ \Omega,
			\\
			\\
			\nu=0  &{\rm on} \;\ \partial\Omega,
		\end{cases}
	\end{equation}
	where $h>0$ is a $C^1$ function in $\Omega$, $q_1,\cdots,q_N$ are distinct points in $\Omega$, $\rho>0$, $\alpha_j>-1$ are constants. As above we need some local/global quantities suitable to describe the combinations of blowup points. Obviously we keep the analogue definitions adopted above about bubbling solutions of (\ref{equ-flat}). \\
	Let $G_{\Omega}$ be the Green's function uniquely defined by
	\begin{equation*}
		\begin{cases}
			-\Delta G_{\Omega}(x,p)=\delta_{p}  &{\rm in} \;\ \Omega,
			\\
			\\
			G_{\Omega}(x,p)=0  &{\rm on} \;\ \partial\Omega,
		\end{cases}
	\end{equation*}
	$R_{\Omega}(x,p)=G_{\Omega}(x,p)+\frac{1}{2\pi}\log |x-p|$ be the regular part of $G_{\Omega}(x,p)$ and 
	\begin{equation*}
		G_{j,\Omega}^*(x)=8\pi (1+\alpha_j)R_{\Omega}(x,p_j)+8\pi \sum_{l\neq j}^{1,\cdots,m}(1+\alpha_l)G_{\Omega}(x,p_l),
	\end{equation*}
	As above we write,
	\begin{equation*}
		%\label{H1}
		H(x)=h(x)\prod_{j=1}^N e^{-4\pi\alpha_j G_\Omega(x,q_j)},
	\end{equation*}
	where we remark that, 
	for $q_j$ fixed, we have 
	\begin{equation*}
		%\label{H2}
		H(x)=h_j(x)|x-q_j|^{2\alpha_j},\quad x\neq q_j
	\end{equation*}
	where $h_j(x)>0$. Similar to notations for the first part, we still assume (\ref{largest-s}) for $(\alpha_1,..,\alpha_N)$ and keep the same conventions about $I_1$ and $\alpha_M$. Next let us define,
	$$L_{\Omega}(\mathbf{p})=\sum_{j\in I_1}\Delta \log h(p_j)h(p_j)^{\frac{1}{1+\alpha_M}}e^{\frac{G_{j,\Omega}^*(p_j)}{1+\alpha_M}},$$
	\begin{align*}
		f_{\Omega}^*(x_{\tau+1},\cdots,x_m)&=\sum_{j=\tau+1}^m\big[\log h(x_j)+4\pi R(x_j,x_j)\big]+4\pi \sum_{l\neq j}^{\tau +1,\cdots,m}G(x_l,x_j),\\
		&\quad+\sum_{j=\tau+1}^m\sum_{i=1}^{\tau}8\pi(1+\alpha_i)G(x_j,p_i),
	\end{align*}
	and let $D^2f_{\Omega}^*$ be the Hessian on $\Omega$. Of course, in this case $(p_{\tau+1},\cdots,p_m)$ is a critical point of $f_{\Omega}^*$ (\cite{ma-wei}). 
	Concerning {\em Class one} we have,
	
	\begin{thm}\label{main-theorem-4}
		Let $\nu_k^{(1)}$ and $\nu_k^{(2)}$ be two sequences of bubbling solutions of {\rm (\ref{equ-flat})}   with $\rho_k^{(1)}=\rho_k=\rho_k^{(1)}$, $(\alpha_1,\cdots,\alpha_N)$ satisfying {\rm (\ref{largest-s})}, $\alpha_M>0$ and $L_{\Omega}(\mathbf{p})\neq 0$. If $\tau<m$ we assume
		$\det \big(D^2f^*(p_{\tau+1},\cdots,p_m)\big)\neq 0$.
		Then there exists $n_0>1$ such that
		$\nu_k^{(1)}=\nu_k^{(2)}$ for all $k\ge n_0$.
	\end{thm}
	
	Concerning {\em Class two} we set
	\begin{align*}\Phi_{j,\Omega}(\mathbf{q}):=~&\sum_{l=1}^m 8\pi(1+\alpha_l) G_{\Omega}(x,q_l)-G_{j,\Omega}^*(q_j)+\log h_j(x)-\log h_j(q_j)\\
		&+2\alpha_j\log |x-q_j|
	\end{align*}
	and
	$$
	D_{\Omega}(\mathbf{p})
	=\lim_{r\to 0}\sum_{j=1}^m h(q_j)e^{G_{j,\Omega}^*(q_j)}\bigg (\int_{\Omega_j\setminus B_{r_j}(p_j)}e^{\Phi_{j,\Omega}(\mathbf{q})}dx-\frac{\pi}{1+\alpha_j}r_j^{-2\alpha_j-2}\bigg ),
	$$
	where $\displaystyle{r_j=r\left(\frac{h(p_j)e^{G_{j,\Omega}^*(p_j)}}{8(1+\alpha_j)^2}\right)^{\frac{1}{2\alpha+2}}}$, $\Omega_1=\Omega$ if $ 1=m$, otherwise we have $\Omega_l\cap \Omega_s=\emptyset$ for $l\neq s$ and $\cup_{j=\tau+1}^m\bar \Omega_j=\bar \Omega$.
   Similar to the previous case we require 
   \begin{equation}\label{neg-crit-2}
	\nabla \bigg (G_{j,\Omega}^*(x)+\log h_j(x)\bigg)\bigg |_{x=q_j}=0,\quad 1\le j\le m.
\end{equation}
	Then we have,
	
	\begin{thm}\label{main-theorem-1}
		Under the same assumptions of Theorem \ref{main-theorem-4}, except that here $\alpha_M\le 0$ and (\ref{neg-crit-2}) holds, in any one of the following cases:
		\begin{enumerate}
			\item $\alpha_M=0$, $L_{\Omega}(\mathbf{p})\neq 0$, $\det \big(D^2f_{\Omega}^*\big )(p_{\tau+1},\cdots,p_m)\neq 0$.
			\item $\alpha_M=0$, $L_{\Omega}(\mathbf{p})=0$, $D_{\Omega}(\mathbf{p})\neq 0$, $\det \big(D^2f_{\Omega}^*\big )(p_{\tau+1},\cdots,p_m)\neq 0$
			\item $\alpha_M<0$, $D_{\Omega}(\mathbf{p})\neq 0$
		\end{enumerate} the same uniqueness property as in Theorem \ref{main-theorem-4} holds true.
	\end{thm}
	
	Since the proof of Theorem \ref{main-theorem-2} is  
	very long and requires several estimates of independent interest, 
	we summarize below the main ideas behind it. One of the crucial steps is to obtain highly refined estimates around each blowup point. The mixture of blowup points forces one to consider different contributions in the many needed expansions each one with its own order of magnitude. The difference of these orders due to the singular sources leads to tremendous difficulties in the asymptotic analysis. This is a major obstacle which we need to overcome, the first tool to achieve this goal being new refined asymptotic estimates around each blowup point, actually far better than those at hand so far. As a matter of fact, known estimates always stops at order two, whenever the second derivatives of coefficient functions shows up. Here we need to carry out at least a fourth order expansion, taking advantage as well of some elegant algebraic identities. Among other things which will be clear in the course of the proof, in this context this has a major advantage since, unlike \cite{bart-4,bart-4-2,wu-zhang-ccm},
	in the very final step of the uniqueness proof, we do not need anymore to consider refined Pohozaev identities, but just rely on the classical one.\\ We do not want to state here such involved expansions of independent interest and refer, respectively, to Theorem \ref{int-pos} and \eqref{def-wk-vip}, \eqref{great-1} for positive singular sources, Theorem \ref{integral-neg} and \eqref{def-wkneg-vip}, \eqref{great-neg} for negative singular sources, and Theorem \ref{reg-int} and Remark \ref{new-inter-e}  for regular blow-up points. 
	The analysis is rather hard, which is why, concerning singular sources, we provide full details in the positive case \eqref{def-wk-vip}, \eqref{great-1}, while we will be more sketchy about the negative case \eqref{def-wkneg-vip}, \eqref{great-neg}. Unlike positive singular sources and regular blow up points, at least to our knowledge, the estimate in Theorem \ref{integral-neg} is the first of its kind. No shortcut seems to be at hand for the regular case, which we have to discuss in full detail.  Actually, in this case, we have to come up with a subtle refinement of an estimate (\cite{gluck,zhang1}) about the vanishing of the gradient of the normalized weight (say $\nabla \bar h_k$, see \eqref{new-rate-1}) as well as, interestingly enough, with the entirely new and rather surprising vanishing estimate of $\nabla \Delta \bar h_k$, see \eqref{new-rate-2}.\\
	
	However, even with those refined estimates at hand, the proof of the uniqueness result is not trivial at all. Indeed previous uniqueness results \cite{bart-4,bart-4-2,wu-zhang-ccm} are based on the following idea: the normalized difference of two solutions (say $\xi_k$) should converge after suitable
	rescaling to some element of the kernel of a linear operator on the plane, the proof being concluded by showing that the coefficients defining this function 
	(say roughly speaking $b_0,b_1,b_2$) vanish one by one. However, since we miss the assumptions mentioned above in \cite{bart-4-2,wu-zhang-ccm}, we face at this stage two genuine new difficulties. First of all the proof of $b_0=0$ requires a further improvement of the estimates about the oscillation of $\xi_k$ and then a sort of boostrap argument, showing a self improving property of this decay. This part of the proof is really delicate, see Lemma \ref{osci-xi-better} and the part
	running from Remark \ref{rem 4.7} to Remark \ref{rem 4.8}. If this were not enough, in this situation to prove that $b_1=0=b_2$ one needs to know the rate of vanishing of the approximating sequence of $b_0$, say $b_{0,k}$, see Lemma \ref{b0-va}. Actually one needs somehow to really restart the contradiction argument ab initio by taking into account the contribution of those terms proportional to $b_{0,k}$ with the aim of estimating again the oscillation of $\xi_k$. This is very long and non trivial and requires other delicate estimates, see for example Lemma \ref{crucial-q12} and Proposition \ref{small-osc-psi}. We split this last part of the proof into seven steps, in an attempt to make it easier to read, see subsection \ref{pf-uni-2} for further details.
	
	\bigskip
	
	Our uniqueness result has several applications. First of all, it furnishes an alternative argument for computing the topological degree of the resolvent operator for equation \eqref{m-equ}, which has been previously established in \cite{chen-lin-deg,chen-lin-deg-2} with and without singularities respectively. In more specific terms,  the authors in \cite{chen-lin-sharp,chen-lin} derived a sharp estimate for characterizing the blow-up rates of bubbling solutions and demonstrated that all bubbling solutions exist within a special set. Subsequently, they computed the degree of the associated Fredholm operator within this set and the desired topological degree of equation \eqref{m-equ} is obtained. With the assistance of this uniqueness result, we can now accurately identify bubbling solutions as the parameter $\rho$ crosses critical positions,  compute the Morse index of certain operators, and calculate the contribution of degree from the bubbling solutions, which yields at once a much direct computation of the topological degree for \eqref{m-equ}. Secondly, it is well-known (\cite{chen-lin-sharp}) that the regular blowup points are critical points of the Kirchhoff-Routh function associated with the ambient manifold. Since according to our result, for a fixed set of singular blow up points, there is only one bubbling solution attached to each such point, we could in principle exactly count the total number of bubbling solutions of this sort, once those critical points are completely classified. This is a relevant 
	fact since the exact counting of solutions is in general a challenging open problem about \eqref{m-equ}, \eqref{equ-flat}. However, as mentioned above, a series of groundbreaking results concerning this point have been recently push forward in the case of the flat two torus with singular sources, see \cite{lin22} and the references quoted therein. Last but not least, we expect the uniqueness result to be the first step to establish non-degeneracy of the bubbling sequences of \eqref{m-equ}, which would allow in turn the evaluation of the specific heat in the negative temperature regime induced by fixed vortices in the statistical mechanics mean field description of 2D flows and plasmas, in the same spirit of \cite{bart-5} for the regular case.

	A word of warning about notations. $B_\tau=B_\tau(0)$ will always denote the ball centered at the origin
	of some local isothermal coordinates $y\in B_\tau$. Whenever $B_\tau$ is any such ball centered at some point $0=y(p), p\in M$, then we will denote by $\Omega(p,\tau)\subset M$ the pre-image of $B_\tau$. 
	On the other side, $B(p,\tau)\subset M$
	will always denote a geodesic ball. Also, in many different estimates we will need some positive small number $\epsilon_0>0$, which will change in general  from estimate to estimate.
	
	\medskip
	
	The organization of this paper is as follows. In section two we collect some preliminary estimates needed for the proof of the main theorems. In section three we derive pointwise estimates for local equations. Then in section four and section five we provide the proofs of the main theorems.
	
	\medskip
	\noindent{\bf Acknowledgement} L. Zhang acknowledges support from Simon's foundation grant 584918. W. Yang acknowledges support from National
Key R\&D Program of China 2022YFA1006800, NSFC No.
12171456 and No. 12271369, FDCT No. 0070/2024/RIA1, UMDF No. TISF/2025/006/FST, No. MYRG-GRG2024-00082-FST and Startup
Research Grant No. SRG2023-00067-FST. D. Bartolucci's Research is partially supported by the MIUR Excellence Department Project MatMod@TOV awarded to the Department of Mathematics, University of Rome Tor Vergata and by PRIN project 2022, ERC PE1\_11,
	"{\em Variational and Analytical aspects of Geometric PDEs}". D. Bartolucci is member of the INDAM Research Group  "Gruppo Nazionale per l’Analisi Matematica, la Probabilità e le loro Applicazioni".
	
	\section{Preliminary Estimates}\label{preliminary}
	
	Since the proof of the main theorems requires delicate analysis, in this section we list some estimates established in \cite{BCLT,BM3,BT,BT-2,chen-lin-sharp,chen-lin,li-cmp,zhang1,zhang2}.

	Let ${\rm w}_k$ be a sequence of solutions of (\ref{r-equ}) with $\rho =\rho_k$  and assume that ${\rm w}_k$ blows up at $m$ points $\{p_1, \cdots,p_m\}$. To describe the bubbling profile of ${\rm w}_k$ near $p_j$, we set
	\begin{equation*}
		u_k={\rm w}_k-\log\bigg(\int_M He^{w_k}{\rm d}\mu \bigg),
	\end{equation*}
	so that $$\int_{M}He^{u_k}{\rm d}\mu=1,$$ and the equation for $u_k$ takes the form,
	\begin{equation*}
		\Delta_g u_k+\rho_k(He^{u_k}-1)=0\quad {\rm in} \ \; M.
	\end{equation*}
	
	From well known results about Liouville-type equations (\cite{BT-2,li-cmp}),
	\begin{equation*}
		u_k-\overline{u_k} \ \to \sum_{j=1}^m 8\pi(1+\alpha_j)G(x,p_j) \quad {\rm in} \ \; {\rm C}_{\rm loc}^2(M\backslash \{p_1,\cdots,p_m\})
	\end{equation*}
	where $\overline{u}_k$ is the average of $u_k$ on $M$:
	$\overline{u}_k=\int_{M}u_k{\rm d}\mu$. For later convenience we fix $r_0>0$ small enough and $M_j\subset M, 1\leq j\leq m$ such that
	\begin{equation*}
		M=\bigcup_{j=1}^m \overline{M}_j;\quad M_j\cap M_l=\varnothing,\  {\rm if}\ j\neq l;\quad B(p_j,3r_0)\subset M_j, \ j=1,\cdots,m.
	\end{equation*}
	According to this definition $M_1=M$, if $m=1$.
	
	The notation about local maximum is of particular relevance in this context. If $p_j$ is a regular blow up point (i.e. $\alpha_j=0$) we define $p_{k,j}$ and $\lambda_{k,j}$  as follows,
	\begin{equation*}
		\lambda_{k,j}=u_k(p_{k,j}) \mathrel{\mathop:}=\max_{B(p_j,r_0)}u_k,
	\end{equation*}
	while if  $p_j$ is a singular blow up point (i.e. $\alpha_j\neq 0$), then we define,
	\begin{equation*}
		p_{k,j}:=p_j\quad \mbox{ and }\quad \lambda_{k,j}:=u_k(p_{k,j}).
	\end{equation*}
	Next, let us define the so called "standard bubble" $U_{k,j}$ to be the solution of
	\begin{equation}\nonumber
		\Delta U_{k,j}+\rho_kh_j(p_{k,j})|x-p_{k,j}|^{2\alpha_j}e^{U_{k,j}}=0 \quad  {\rm in} \ \; \mathbb{R}^2
	\end{equation}
	which takes the form (\cite{CL1,CL2,Parjapat-Tarantello}),
	\begin{equation}\nonumber
		U_{k,j}(x)=\lambda_{k,j}-2\log\Big(1+\frac{\rho_k h_j(p_{k,j})}{8(1+\alpha_j)^2}e^{\lambda_{k,j}}|x-p_{k,j}|^{2(1+\alpha_j)}\Big).
	\end{equation}
	
	It is well-known (\cite{BCLT,BT-2,li-cmp}) that $u_k$ can be approximated by the standard bubbles $U_{k,j}$ near $p_j$ up to a uniformly bounded error term:
	\begin{equation*}
		\big|u_k(x)-U_{k,j}(x)\big| \leq C, \quad x\in B(p_j,r_0).
	\end{equation*}
	
	As a consequence, in particular we have,
	\begin{equation}\nonumber
		|\lambda_{k,i}-\lambda_{k,j}|\leq C, \quad 1\leq i,j \leq m,
	\end{equation}
	for some $C$ independent of $k$.
	
	\medskip
	In case $\tau<m$, it has been shown in \cite{chen-lin-sharp} that,
	\begin{equation}\label{first-deriv-est}
		\nabla(\log h+G_j^*)(p_{k,j})=O(\lambda_{k,j}e^{-\lambda_{k,j}}),\quad \tau+1\leq j\leq m,
	\end{equation}
	which, in view of the non-degeneracy condition
	$$\det\big(D^2f^*(p_{\tau+1},\cdot,p_m)\big)\neq0,$$
	readily implies that,
	\begin{equation}\label{p_kj-location}
		|p_{k,j}-p_j|=O(\lambda_{k,j}e^{-\lambda_{k,j}}),\quad \tau+1\leq j\leq m.
	\end{equation}
	Later, sharper estimates were obtained in \cite{chen-lin,zhang2} for $1\leq j \leq \tau$ and
	in \cite{chen-lin-sharp,gluck,zhang1} for $\tau+1\leq j\leq m$.
	
	By using $\lambda_i^k$ to denote the maximum of $u_i^k$, $i=1,2$ which share the same value of the parameter $\rho_k$, then it is a simple consequence of $L(\mathbf p)\neq 0$ (see \cite{bart-4,wu-zhang-ccm}) that 
	\begin{equation}\label{initial-small}
		|\lambda_1^k-\lambda_2^k|\leq Ce^{-\gamma \lambda_1^k},\mbox{ \rm for some }\gamma>0.
	\end{equation}

	Let us also recall that it has been established in \cite{BM3,BT} that (see \eqref{a-note})  $\rho_*=\lim\limits_{k\to +\infty}\rho_k$.
	Concerning the difference between $\rho_k$ and $\rho_*$,
	we set
	$$\rho_{k,j}=\rho_k\int\limits_{\Omega(p_{j},\tau)}He^{u_{k}}.$$For $\alpha_j\geq 0$, the following estimates hold (see  \cite{chen-lin-sharp}, \cite{chen-lin}):
	\begin{thm}
		%\label{Theorem chen-lin2}
		There exists $\epsilon_0>0$ and $d_j>0$ such that,
		\begin{align*}
			\begin{cases}\rho_{k,j}-8\pi(1+\alpha_j)=2\pi d_je^{-\frac{\lambda_{k,j}}{1+\alpha_j}}+
				O(e^{-\frac{1+\epsilon_0}{1+\alpha_j}\lambda_{k,j}}),  &1\leq j\leq \tau, \\
				\\
				\rho_{k,j}-8\pi=O\big(\lambda_{k,j}e^{-\lambda_{k,j}}\big), &\tau+1\leq j\leq m.
			\end{cases}
		\end{align*}

		% &\rho_k-\rho_*=L^*e^{-\frac{\lambda_M^k}{1+\alpha_M}}+O(e^{-\frac{\lambda_M^k}{1+\alpha_M}-\tau_1}) && \mbox{\rm where }\tau_1>0, \label{rho-k}
		
		%	with $$L^*=\dfrac{2\pi^2}{(1+\alpha_M)\sin\frac{\pi}{1+\alpha_M}}L(\mathbf{p}),\quad \lambda_M^k=\max\{\lambda_{k,1},...,\lambda_{k,m}\}.$$
	\end{thm}
	
	As mentioned above, we were not aware of any result of this sort about the case $\alpha_j\in (-1,0)$ and we refer to the new estimates in Theorem  \ref{integral-neg} below concerning this point.\\ 
	Assuming without loss of generality that $\alpha_1=\alpha_M$, then it has been proved in \cite{chen-lin} that
	\begin{equation}\label{small-rho-k}
		\rho_k-\rho_*= L(\mathbf{p})e^{-\frac{\lambda_1^k}{1+\alpha_M}}+O(e^{-\frac{1+\epsilon_0}{1+\alpha_M}\lambda_{k,j}})
	\end{equation}
	for some $\epsilon_0>0$ and $L(\mathbf{p})$ as defined in (\ref{L-p}).
	Next let us state two relevant results about suitably defined global linear problems.
	
	\subsection{A uniqueness lemma}
	\begin{lem}
		\label{lem1}
		Let $\alpha>-1$, $\alpha\not \in \mathbb N\cup \{0\}$ and $\phi$ be a $C^2$ solution of
		$$\begin{cases}
			\Delta \phi+8(1+\alpha )^2|x|^{2\alpha}e^{U_{\alpha}}\phi=0\;\;
			\hbox{ {\rm in} }\;\; \mathbb R^2,\\
			\\
			|\phi(x)|\le C(1+|x|)^{\tau},
			\quad x\in \mathbb R^2,
		\end{cases}
		$$
		where
		$$
		U_{\alpha}(x)=\log \frac{1}{(1+|x|^{2\alpha+2})^2}
		\quad \hbox{ and } \quad
		\tau\in [0,1).
		$$
		Then there exists some constant $b_0$ such that
		\begin{equation*}
			\phi(x)= b_0\frac{1-|x|^{2(1+\alpha)}}{1+|x|^{2(1+\alpha)}}.
		\end{equation*}
	\end{lem}
	
\begin{proof}[Proof of Lemma \ref{lem1}.] The proof for $\alpha>0$ can be found in \cite{chen-lin,wu-zhang-ccm}. We just focus on $-1<\alpha<0$.  Let $b_0=\phi(0)$, we set
	$$\Phi=\phi-b_0\frac{1-|x|^{2(1+\alpha)}}{1+|x|^{2(1+\alpha)}}. $$
	Let $l\ge 1$ be an
	integer. We define
	$$\phi_l(r)=\frac 1{2\pi}\int_0^{2\pi}\Phi(r\cos \theta, r\sin
	\theta)\cos (l\theta)d\theta.$$
	Then $\phi_l$ satisfies
	\begin{equation}
		\label{nov9eq2}
		\begin{cases}\phi_l''(r)+\frac
			1r\phi_l'(r)+(8(1+\alpha )^2r^{2\alpha}e^{U_{\alpha}}-\frac{l^2}{r^2})\phi_l(r)=0, \quad
			0<r<\infty,\\
			\\
			\lim_{r\to 0}\phi_l(r)=0,\quad |\phi_l(r)|\le C(1+r)^{\tau}, \quad
			r>0.
		\end{cases} 
	\end{equation}
	
	We claim that $\phi_l \equiv 0$ for all $l\ge 1$.
	Note that by assumption $\phi_l(r)=o(r)$ at infinity and $\phi_l(0)=0$. Then by direct computation
	we verify that the following functions are two linearly independent fundamental
	solutions of (\ref{nov9eq2}):
	\begin{equation}\label{two-fun}
		\begin{split}
			F_{1,l}(r)&=\frac{(\frac{l}{1+\alpha}+1)r^l+(\frac{l}{1+\alpha}-1)r^{l+2(1+\alpha)}}{1+r^{2(1+\alpha)}}, \\
			F_{2,l}(r)&=\frac{(\frac{l}{1+\alpha}+1)r^{-l+2(1+\alpha)}+(\frac{l}{1+\alpha}-1)r^{-l}}{1+r^{2(1+\alpha)}}.
		\end{split}
	\end{equation}

	Consequently $\phi_l=c_1F_{1,l}+c_2F_{2,l}$ where $c_1$ and $c_2$ are
	two constants. From (\ref{nov9eq2}) we see that $c_1=0$ because
	otherwise $|\phi_l|\sim r^{l}$ at $\infty$. Similarly by
	observing the behavior of $\phi_l$ at $0$ we have $c_2=0$ because
	otherwise $|\phi_l|\sim r^{-l}$ near $0$. So $\phi_l\equiv 0$ for all
	$l\ge 1$. The same argument also shows that the projection of
	$\phi$ along $\sin (l\theta)$ ($\forall l\ge 1$) is $0$. Finally,
	from $\Phi(0)=0$ the projection of $\Phi$ along the constant direction has a zero initial condition.
	Besides the known fundamental solution of this ODE
	(which is of course $\frac{1-|x|^{2(1+\alpha)}}{1+|x|^{2(1+\alpha)}}$),
	the other one has a logarithmic growth near the origin and infinity.
	This is done by standard ODE techniques,
	see for example \cite{zhang2}. The boundedness of $\Phi$ near zero rules out the
	second fundamental solution while the zero initial condition rules out the first,
	whence we conclude that $\Phi(x)\equiv 0$.
\end{proof}
	
	\bigskip
	
	For $\alpha=0$,  the following lemma has been proved in \cite{chen-lin-sharp}.
	\begin{lem}
		%\label{linear-lem-2}
		Let  $ \varphi $ be a $ C^2 $ solution of
		\begin{equation*}
			\begin{cases}
				\Delta \varphi+e^U\varphi=0\quad & {\rm in} \ \;\mathbb{R}^2,
				\\
				\\
				|\varphi| \leq c\big(1+|x|\big)^{\kappa} \quad & {\rm in} \ \;\mathbb{R}^2,
			\end{cases}
		\end{equation*}
		where $ U(x)=\log\frac{8}{(1+|x|^2)^2} $ and $ \kappa \in[0,1) $. Then there exist constants $b_0$, $b_1$, $b_2$ such that
		\begin{equation*}
			\varphi= b_0\varphi_0+b_1\varphi_1+b_2\varphi_2,
		\end{equation*}
		where
		\begin{equation*}
			\varphi_0(x)= \frac{1-|x|^2}{1+|x|^2},\quad \varphi_1(x)= \frac{x_1}{1+|x|^2},\quad \varphi_2(x)= \frac{x_2}{1+|x|^2}.
		\end{equation*}
		
	\end{lem}

	\medskip
	
	\section{Local asymptotic analysis near a blowup point}
	%\label{sec3}
	In this section we study the asymptotic expansion of a sequence of blowup solutions near a non-quantized singular source. By a straightforward change of variables and with an obvious abuse of notations, we assume that $u_k$ denotes in local coordinates a sequence of blowup solutions of
	\begin{equation*}
		\Delta u_k+|x|^{2\alpha}h_k(x)e^{u_k}=0 \quad \mbox{in}\quad B_{2\tau}
	\end{equation*}  
	where $B_{\tau}\subset \mathbb{R}^2$ for $\tau>0$ is the ball centered at the origin with radius $\tau>0$,
	\begin{equation}\label{assum-H}
		\frac{1}{C}\le h_k(x)\le C, \quad \|D^{m}h_k\|_{L^{\infty}(B_{2\tau})}\le C
		\quad \forall x\in B_{2\tau}, \quad \forall m\le 5,
	\end{equation}
	and
	\begin{equation}\nonumber
		\int_{B_{2\tau}}|x|^{2\alpha}h_k(x)e^{u_k}<C,
	\end{equation}
	for some uniform $C>0$. The sequence $\{u_k\}$ has its only blowup point at the origin:
	\begin{equation}\label{only-bu}
		\mbox{there exists}\quad  z_k\to 0, \mbox{ such that} \quad \lim_{k\to \infty} u_k(z_k)\to \infty
	\end{equation}
	and for any fixed $K\Subset B_{2\tau}\setminus \{0\}$, there exists $C(K)>0$ such that
	\begin{equation}\label{only-bu-2}
		u_k(x)\le C(K),\quad x\in K.
	\end{equation}
	It is also standard to assume that $u_k$ has bounded oscillation on $\partial B_{\tau}$:
	\begin{equation}\label{u-k-boc}
		|u_k(x)-u_k(y)|\le C\quad \forall x,y\in \partial B_{\tau},
	\end{equation}
	for some uniform  $C>0$.
	The following harmonic function $\psi_k$ is needed to encode the boundary oscillation of $u_k$:
	\begin{equation*}
		\begin{cases}
			\Delta \psi_k=0  \quad \hbox{in } B_{\tau},\\
			\\
			\psi_k=u_k-\frac{1}{2\pi \tau}\int_{\partial B_{\tau}}u_kdS
			\quad \hbox {on } \partial B_{\tau}.    
		\end{cases}
	\end{equation*}
	Since $u_k$ has a bounded oscillation on $\partial B_{\tau}$ (see \eqref{u-k-boc}), $u_k-\psi_k$ is constant on $\partial B_{\tau}$, the bounded oscillation of $u_k$ means   $\|D^m\psi_k\|_{B_{\tau/2}}\le C(m)$ for any $m\in \mathbb N$.
	The mean value property of harmonic functions also gives $\psi_k(0)=0$.
	For the sake of simplicity let us define,
	$$
	V_k(x)=h_k(x)e^{\psi_k(x)}, \quad \varepsilon_k=e^{-\frac{u_k(0)}{2(1+\alpha)}}
	$$
	and
	\begin{equation} \label{defvk-1}
		v_k(y)=u_k(\varepsilon_ky)+2(1+\alpha)\log \varepsilon_k-\psi_k(\varepsilon_ky),\quad y\in \Omega_k:=B_{\tau\varepsilon_k^{-1}}(0).
	\end{equation}
	Obviously we have $V_k(0)=h_k(0)$, $\Delta \log V_k(0)=\Delta \log h_k(0)$.
	Let
	\begin{equation}\label{st-bub}
		U_k(y)=
		-2\log
		\left(1+\frac{h_k(0)}{8(\alpha+1)^2}|y|^{2\alpha+2} \right),
	\end{equation}
	be the standard bubble which satisfies
	\begin{equation}\nonumber \Delta U_k+|y|^{2\alpha}h_k(0)e^{U_k(y)}=0 \quad
		\hbox{ in } \mathbb R^2.
	\end{equation}
	
	In the remaining part of this section we split the discussion into three subsections according to the sign of $\alpha$ ($\alpha>0$, $\alpha<0$ and $\alpha=0$).
	
	\subsection{Asymptotic analysis for $\alpha>0$}
	
	We say that a function is separable if it is the product of a polynomial of two variables $P(x_1,x_2)$ times a radial function with certain decay at infinity, where $P$ is the sum of monomials $x_1^nx_2^m$ where $n,m$ are non negative integers and at least one among them is odd. The main result of this subsection is
	
	\begin{thm}\label{int-pos}
		\begin{align*}
			\int_{B_{\tau}}h_k(x)|x|^{2\alpha}e^{u_k}
			=~&8\pi(1+\alpha)\left(1-\frac{e^{-u_k(0)}}{e^{-u_k(0)}+\frac{h_k(0)}{8(1+\alpha)^2}\tau^{2\alpha+2}}\right)\\
			&+d_{1,k} \Delta \log h_k(0)\varepsilon_k^2
			+T_1(V_k)\varepsilon_k^4+O(\varepsilon_k^{4+\epsilon_0}),
		\end{align*}
		for some $\epsilon_0>0$, where,
		\begin{align*}
			d_{1,k}=&-\frac{2\pi^2}{(1+\alpha ) \sin \frac{\pi}{1+\alpha}}\left(\frac{8(1+\alpha)^2}{V_k(0)}\right)^{\frac 1{1+\alpha}}\\
			&+\left(\frac{16\pi(1+\alpha)^2}{\alpha }\frac{|\nabla \log V_k(0)|^2}{V_k(0)\Delta\log V_k(0)}+\frac{16\pi(1+\alpha)^4}{\alpha V_k(0)}\right)\tau^{-2\alpha}\varepsilon_k^{2\alpha},
		\end{align*}
		and $T_1(V_k)$ is a bounded function of $D^{\beta}V_k(0)$ for $|\beta |=0,1,2,3,4$.
	\end{thm}

\noindent {\em Proof of Theorem \ref{int-pos}.}
	Let $\varepsilon_k=e^{-\frac{u_k(0)}{2+2\alpha}}$, we denote by $z_k$ the maximum point of $u_k$. It is proved in
	\cite{BCLT,BT-2} that $\varepsilon_k^{-1}z_k\to 0$. Note that we are allowed to use
	$u_k(0)$ to define $\varepsilon_k$ since it follows from the results in \cite{BCLT,BT-2} that
	$u_k(0)/u_k(z_k)=1+o(1)$. From the definition of $\psi_k$,
	$u_k-\psi_k$ has vanishing oscillation on $\partial B_1$ and it satisfies
	$$
	\Delta (u_k-\psi_k)+|x|^{2\alpha}V_k(x)e^{(u_k-\psi_k)(x)}=0 \quad
	\hbox{ in } B_1.$$
	Without loss of generality we assume
	$$V_k(0)=h_k(0)\to 8(1+\alpha)^2.$$
	A straightforward inspection of the arguments in  \cite{BCLT,BT-2} shows that, letting $x_k$ be the maximum point of $u_k-\psi_k$,
	then we still have $\varepsilon_k^{-1}x_k\to 0$. Next let us define $v_k$ as in \eqref{defvk-1},
	$$
	v_k(y)=u_k(\varepsilon_ky)-\psi_k(\varepsilon_ky)-u_k(0),
	\quad y \in \Omega_k := B_{\tau\varepsilon_k^{-1}}(0),
	$$
	which satisfies,
	$$\begin{cases}
		\Delta v_k(y)+|y|^{2\alpha}V_k(\varepsilon_ky)e^{v_k(y)}=0,\quad  y\in
		\Omega_k,\\
		\\
		v_k(0)=0,\quad y_k:=\varepsilon_k^{-1}x_k\to 0,\\
		\\
		v_k(y)\to (-2)\log(1+|y|^{2\alpha+2})\quad
		\mbox{in}\quad C^2_{\rm loc}(\mathbb R^2),\\
		\\
		v_k(y_1)=v_k(y_2),\quad \forall y_1,y_2\in \partial
		\Omega_k=\partial
		B_{\tau\varepsilon_k^{-1}}(0).
	\end{cases}
	$$
	Note that in the second equation above we use $y_k$ to denote the
	maximum point of $v_k$. Let $U_k$ be the standard bubble defined above,
	$$U_k(y)=(-2)
	\log
	\left(1+\frac{V_k(0)}{8(\alpha+1)^2}|y|^{2\alpha+2} \right),$$
	(recall that $V_k(0)=h_k(0)$) then it is well known that for all $\alpha>-1$ (see \cite{BCLT,BT-2,zhang2} )
	$$|v_k(y)-U_k(y)|\le C ,\quad |y|\le \varepsilon_k^{-1}. $$
	Therefore $U_k$ is the first term in the expansion of $v_k$ and to identify
	the second term we need the radial function $g_k$ which satisfies:
	\begin{align*}
		\begin{cases}
			g_k''(r)+\frac
			1rg_k'(r)+(r^{2\alpha}V_k(0)e^{U_k(r)}-\frac
			1{r^2})g_k(r)
			=-V_k(0)r^{2\alpha+1}e^{U_k(r)},~r\in (0,\infty)\\
			\\
			%qquad \mbox{for }\quad 0<r<\infty, 
			\lim\limits_{r\to 0+}g_k(r)=\lim\limits_{r\to \infty}g_k(r)=0,\quad |g_k(r)|\le
			C\frac r{1+r^{2\alpha+2}}.
		\end{cases}
	\end{align*}
	By direct computation one checks that in fact
	\begin{equation}
		\label{mar27e1} g_k(r)=-\frac{2(1+\alpha)}{\alpha }
		\frac{r}{1+\frac{V_k(0)}{8(1+\alpha)^2}r^{2\alpha+2}},
	\end{equation}
	solves that equation. Let
	\begin{equation}
		\label{nov10e6}
		c_{1,k}(y)=\varepsilon_kg_k(r)\nabla (\log V_k)(0)\cdot \theta, \quad \theta=(\theta_1,\theta_2),\; \theta_j=y_j/r, j=1,2.
	\end{equation}
	Then by a straightforward evaluation we see that $c_{1,k}$ satisfies,
	\begin{equation}
		\label{mar31e1} \Delta c_{1,k}+|y|^{2\alpha}V_k(0)e^{U_k(y)}c_{1,k}
		=-\varepsilon_k(\nabla (\log V_k)(0)\cdot y) V_k(0)|y|^{2\alpha}e^{U_k(y)},\quad
	\end{equation}
	in $\Omega_k$.  We claim that $c_{1,k}$ is the second term in the expansion of $v_k$. Let
	\begin{equation}\nonumber
		w_k(y)=v_k(y)-U_k(y)-c_{1,k},
	\end{equation}
	and observe that the oscillation of $w_k$ on $\partial \Omega_k$ is $O(\varepsilon_k^{2+2\alpha})$.
	
	\begin{prop}
		\label{nov7p1}
		For any $\delta>0$  small enough, there exists $C(\delta)>0$ such
		that for all $k$ large enough we have,
		$$|w_k(y)|\le
		C\varepsilon_k^{2}(1+|y|)^{\delta}\quad \mbox{{\rm in}}\quad \Omega_k.
		$$
	\end{prop}

	\begin{proof}[Proof of Proposition \ref{nov7p1}.]
	By using the equations satisfied by $U_k$ and $c_{1,k}$ we see that
	$$\Delta (U_k+c_{1,k})+|y|^{2\alpha}V_k(\varepsilon_ky)e^{U_k+c_{1,k}}=
	O(\varepsilon_k^2)(1+r)^{-2-2\alpha}\quad \mbox{in}\quad \Omega_k. $$
	Then we can write the equation for $w_k$ as follows,
	$$\begin{cases}
		\Delta
		w_k+|y|^{2\alpha}V_k(\varepsilon_k y)e^{\xi_k}w_k=O(\varepsilon_k^2)(1+r)^{-2-2\alpha},\\ \\
		w_k(0)=0,\quad |w_k(y)|\le C, \quad y\in \Omega_k,\quad
		w_k|_{\partial \Omega_k}=\tilde a_k+O(\varepsilon_k^{2+2\alpha}),
	\end{cases}
	$$
	where $\xi_k$ is obtained by the mean value theorem, $e^{\xi_k}=(e^{v_k}-e^{U_k+c_{1,k}})/w_k$, and consequently satisfies $e^{\xi_k(\eta)}\le C(1+|\eta |)^{-4-4\alpha}$.  Let
	$$\tilde \Lambda_k=\max\frac{|w_k(y)|}{\varepsilon_k^{2}(1+|y|)^{\delta}},\quad y\in
	\bar \Omega_k.$$
	Our goal is to show $\tilde \Lambda_k=O(1)$. We prove this by a contradiction and assume that $\widetilde \Lambda_k\to \infty$. Then let $\widetilde y_k$  denote a point where $\widetilde \Lambda_k$ is attained and let us define
	$$\widetilde w_k(y)=\frac{w_k(y)}{\widetilde \Lambda_k\varepsilon_k^{2}(1+|\widetilde y_k|)^{\delta}},$$
	which satisfies,
	\begin{equation}\label{stan-1}\Delta \widetilde w_k+|y|^{2\alpha}V_k(\varepsilon_k y)e^{\xi_k}\widetilde w_k
		=O(1)\frac{(1+r)^{-2-2\alpha}}{\widetilde \Lambda_k(1+|\widetilde y_k|)^{\delta}}.
	\end{equation}
	Since the oscillation of $w_k$ on $\partial \Omega_k$ is of order $O(\varepsilon_k^{2+2\alpha})$, the oscillation of $\widetilde w_k$ on $\partial \Omega_k$ is of order $o(1)(1+|\widetilde y_k|)^{-\delta}$.  Moreover  by the definition of $\widetilde w_k$ we have
	\begin{equation}\label{stan-2}|\widetilde w_k(y)|=\left|\frac{w_k(y)}{\varepsilon_k^{2}\tilde \Lambda_k (1+|y|)^{\delta}}\frac{(1+|y|)^{\delta}}{(1+|\widetilde y_k|)^{\delta}}\right|\le \frac{(1+|y|)^{\delta}}{(1+|\widetilde y_k|)^{\delta}}.
	\end{equation}
	Therefore we conclude that, along a subsequence if necessary, $\widetilde w_k$ converges in $C^{\beta}_{\rm loc}(\mathbb R^2)$ (for some $\beta\in (0,1)$) to a
	solution $w$ of
	$$\left\{\begin{array}{ll}
		\Delta w+r^{2\alpha}e^{U_{\alpha}}w=0,\quad \mathbb R^2,\\
		\\
		w(0)=0,\quad |w(y)|\le C(1+|y|)^{\delta}.
	\end{array}
	\right.
	$$
	If we had, possibly along a sub-subsequence, that $\widetilde y_k$ converge to $y_0\in \mathbb R^2$, then we would also have
	$|w(y_0)|=1$ by continuity. However this is impossible because an
	application of Lemma \ref{lem1} shows that $w\equiv 0$. Therefore necessarily  $|\widetilde y_k|\to +\infty$.
	%Remark that since $\alpha>-1$ then $|y|^{2\alpha}\in L^{1+\epsilon_0}$ for some $\epsilon_0>0$, whence $w_k$ is at least Holder continuous near the origin.
	It follows from $|\widetilde w_k(\widetilde y_k)|=1$ and the Green representation formula that
	\begin{equation}
		\label{jul17e1}
		\begin{aligned}
			\pm 1=~&\widetilde w_k(\widetilde y_k)-\widetilde w_k(0)\\
			=~&\int_{\Omega_k}(G_k(\widetilde y_k,\eta)-G_k(0,\eta))\\
			&\times\bigg
			\{|\eta
			|^{2\alpha}V_k(\varepsilon_k \eta)e^{\xi_k(\eta)}
			\widetilde w_k(\eta)
			+\frac{O(1)(1+|\eta |^{-2-2\alpha})}{\widetilde \Lambda_k(1+|\widetilde
				y_k|)^{\delta}} \bigg \}d\eta+o(1)
		\end{aligned}
	\end{equation}
	where $G_k$ is the Green function  with  Dirichlet boundary condition on $\Omega_k$, which takes the form,
	$$G_k(y,\eta)=-\frac 1{2\pi}\log |y-\eta |
	+\frac 1{2\pi}\log
	\left(
	\frac{|y|}{\varepsilon_k^{-1}}\left|\frac{\varepsilon_k^{-2}y}{|y|^2}-\eta\right|
	\right).$$
	Note that in (\ref{jul17e1}) we used the fact that $\mbox{\rm osc}_{\partial \Omega_k}\widetilde w_k=o(1)$.  We deduce from  (\ref{jul17e1})  that,
	\begin{equation}\label{ccm1}
		1\le \int_{\Omega_k}|G_k(\widetilde
		y_k,\eta)-G_k(0,\eta)|  \bigg (\frac
		{(1+|\widetilde y_k|)^{-\delta}}{(1+|\eta
			|)^{4+2\alpha-\delta}}
		+o(1)\frac{(1+|\widetilde y_k|)^{-\delta}}{(1+|\eta
			|)^{2+2\alpha}}\bigg ) d\eta,
	\end{equation}
	where we used
	$$|{\widetilde w_k(\eta)}|\le \frac{(1+|\eta|)^{\delta}}{(1+|y_k|)^{\delta}},
	\quad e^{\xi_k(\eta)}\le C(1+|\eta |)^{-4-4\alpha},\quad \tilde
	\Lambda_k\to \infty.$$
	We will come up with a contradiction to (\ref{ccm1}) by showing that:
	\begin{equation} \label{ccm2}
		\int_{\Omega_k}|G_k(\widetilde y_k,\eta)-G_k(0,\eta)|\frac{(1+|\eta|)^{\delta-4-2\alpha}+o(1)(1+|\eta
			|)^{-2-2\alpha}}
		{(1+|\widetilde y_k|)^{\delta}}d\eta=o(1).
	\end{equation}
	There are only two possibilities, either along a subsequence $|\widetilde y_k|=o(1)\varepsilon_k^{-1}$, or $|\widetilde y_k|\sim \varepsilon_k^{-1}$.
	Assume first that $|\widetilde y_k|=o(1)\varepsilon_k^{-1}$ and observe that $G_k(\widetilde y_k, \eta)$ can be written as follows
	\begin{align*}
		G_k(\widetilde y_k,\eta)=-\frac 1{2\pi}\log |\widetilde y_k-\eta |+\frac 1{2\pi}\log
		\varepsilon_k^{-1}+\frac{1}{2\pi}\log \bigg |\frac{\widetilde y_k}{|\widetilde y_k|}-\varepsilon_k^2|\widetilde y_k|\eta \bigg |.
	\end{align*}
	Thus we are reduced to show that
	\begin{equation}
		\label{tem-9}
		\begin{aligned}
			&\int_{\Omega_k} \bigg ( \bigg |\log \frac{|\tilde y_k-\eta |}{|\eta |}
			\bigg |+\bigg | \log \big |\frac{\tilde y_k}{|\tilde y_k|}-\varepsilon_k^2|\tilde y_k| \eta \big |\bigg | \bigg )\\
			&\quad\times\bigg (\frac{(1+|\eta |)^{-4-2\alpha}}{(1+|\tilde
				y_k|)^{\delta}}
			+o(1)\frac{(1+|\eta
				|)^{-2-2\alpha}}{(1+|\tilde y_k|)^{\delta}}\bigg )
			d\eta=o(1),
		\end{aligned}
	\end{equation}
	which is a consequence of elementary estimates. In fact,  for fixed $\widetilde y_k$, we decompose $\Omega_k$ as follows $\Omega_k=E_1\cup E_2\cup E_3$ where
	$$E_1:=\left\{\eta\in \Omega;~|\eta|<\frac{|\widetilde y_k|}{2}\right\},\quad E_2=\left\{\eta\in \Omega_k;~|\eta -\widetilde y_k|<\frac{|\widetilde y_k|}{2}\right\},$$ and $$E_3=\Omega_k\setminus (E_1\cup E_2).$$ 
	The integrals on the left of (\ref{tem-9}) over the regions $E_j$, $j=1,2,3$ can be easily verified to be $o(1)$ in view of the following elementary inequalities:
	$$\bigg |\log \frac{|\widetilde y_k-\eta |}{|\eta |}\bigg |\le 
	\begin{cases}
		C\log |\widetilde y_k|+C|\log|\eta||,\quad &\eta\in E_1,\\
		\\
		C\log |\widetilde y_k|+|\log |\widetilde y_k-\eta ||,\quad &\eta\in E_2,\\
		\\
		C|\widetilde y_k|/|\eta |,\quad &\eta \in E_3.
	\end{cases}
	$$
	The second logarithmic term can be estimated as follows,
	$$\bigg |\log \left|\frac{\tilde y_k}{|\tilde y_k|}-\varepsilon_k^2 |\tilde y_k| \eta \,\right  | \bigg |\le C\varepsilon_k^2 |\tilde y_k| |\eta |,\quad |\tilde y_k|=o(1)\varepsilon_k^{-1}. $$
	Next we consider the case $|\tilde y_k|\sim \varepsilon_k^{-1}$. For the Green function it is enough to observe that
	%$$|G_k(\tilde y_k,\eta)-G_k(0,\eta)|\le C(\log (1+|\eta |)+\log \varepsilon_k^{-1}),$$ 
	$${
		|G_k(\tilde y_k,\eta)-G_k(0,\eta)|\le C\left(|\log |\eta ||+|\log|y_k-\eta||+\log \frac{1}{\varepsilon_k}\right),}$$
	for some uniform $C>0$. Then
	(\ref{ccm2}) follows once more by elementary estimates. This fact concludes the proof of Proposition \ref{nov7p1} and thus the first step in the expansion of $v_k$. \end{proof}

	\begin{rem}\label{stand-rem}
		We will often need elliptic estimates of the same sort as those used in the proof of Proposition \ref{nov7p1} above. To avoid repetitions we point out here a few key assumptions, which, whenever verified at the same time, allow the adoption of the same argument. First of all, in equation (\ref{stan-1}), the right hand side has to be $O(r^{-2-\delta})$ for some positive $\delta$ and $r$ large enough. Moreover we need the inequality in (\ref{stan-2}) , $w_k(0)=0$ and $w_k$ having small enough oscillation on the boundary.
	\end{rem}
	\medskip
	Before making further improvements in the expansion, we recall that
	$\theta_1=y_1/r=\cos \theta$, $\theta_2=y_2/r=\sin \theta$.
	Now we calculate the equation of $U_k+c_{1,k}$ in more detail. We first write down the equation for $U_k+c_{1,k}$ as follows,
	\begin{align*}
		&\Delta (U_k+c_{1,k})+|y|^{2\alpha}V_k(\varepsilon_ky)e^{U_k+c_{1,k}}\\
		&=\Delta U_k+\Delta c_{1,k}+|y|^{2\alpha}V_k(0)e^{U_k+c_{1,k}+\log \frac{V_k(\varepsilon_ky)}{V_k(0)}}\\
		&=\Delta U_k+\Delta c_{1,k}+|y|^{2\alpha}V_k(0)e^{U_k}\bigg (1+c_{1,k}+
		\log \frac{V_k(\varepsilon_ky)}{V_k(0)}\\
		&\quad+\frac 12\left(c_{1,k}+\log \frac{V_k(\varepsilon_ky)}{V_k(0)}\right)^2+O(\varepsilon_k^3r^3)\bigg ).
	\end{align*}
	The Taylor expansion of the logarithmic term gives
	$$\log \frac{V_k(\varepsilon_ky)}{V_k(0)}=\varepsilon_k \nabla \log V_k(0)\cdot \theta r+\varepsilon_k^2(r^2\Theta_2+\mathcal{R}_2)+O(\varepsilon_k^3r^3),$$
	where
	\begin{align}\nonumber
		\Theta_2=~&\frac 12\partial_{11}\log V_k(0)(\theta_1^2-\frac 12)
		+\frac 12 \partial_{22}\log V_k(0)(\theta_2^2-\frac 12)
		+\partial_{12}\log V_k(0)\theta_1\theta_2 \nonumber\\
		=~&\frac{1}{4}\bigg (\cos 2\theta(\partial_{11}\log V_k(0)-\partial_{22}\log V_k(0))+2\sin 2\theta \partial_{12}\log V_k(0))\bigg )\nonumber
	\end{align}
	and
	\begin{equation}\nonumber
		\mathcal{R}_2=\frac 14 r^2\Delta \log V_k(0).
	\end{equation}
	Then we have 
	\begin{align*}
		c_{1,k}+\log \frac{V_k(\varepsilon_ky)}{V_k(0)}
		=\epsilon_k\nabla \log V_k(0)\cdot \theta (r+g_k)+\varepsilon_k^2r^2\Theta_2+\varepsilon_k^2\mathcal{R}_2+O(\varepsilon_k^3)r^3,
	\end{align*}
	and 
	\begin{align*}
		&(c_{1,k}+\log \frac{V_k(\varepsilon_ky)}{V_k(0)})^2=\varepsilon_k^2(\nabla \log V_k(0)\cdot \theta)^2(g_k+r)^2+O(\varepsilon_k^3)r^3\\
		&=\varepsilon_k^2\bigg ( \frac 12 |\nabla \log V_k(0)|^2+\frac{\cos 2\theta}{2}\left((\partial_1\log V_k(0))^2-(\partial_2\log V_k(0))^2\right)\\
		&\quad +\partial_1(\log V_k)(0)\partial_2(\log V_k)(0)\sin 2\theta\bigg )(g_k+r)^2+O(\varepsilon_k^3r^3).
	\end{align*}
	
	By using these estimates together with the equation for $U_k+c_{1,k}$ we have, 
	\begin{equation}
		\label{Uk-c1k}
		\begin{aligned}
			&\Delta (U_k+c_{1,k})+|y|^{2\alpha}V_k(\varepsilon_k y)e^{U_k+c_{1,k}} \\
			&=\varepsilon_k^2V_k(0)r^{2\alpha}e^{U_k}\bigg (\frac 14\Delta (\log V_k)(0)r^2+\frac 14|\nabla (\log V_k)(0)|^2(g_k+r)^2\\
			&\quad+ \hat \Theta_2\bigg )+O(\varepsilon_k^3)(1+r)^{-1-2\alpha}, 
		\end{aligned}
	\end{equation}
	where
	\begin{equation}
		\label{E212}
		\begin{aligned}
			\hat \Theta_2&=( E_{2,1}(r)\cos 2\theta+E_{2,2}(r)\sin 2\theta)\\
			&=\bigg (\bigg(\frac 14\left(\partial_{11}(\log V_k)(0)-\partial_{22}(\log V_k)(0)\right)r^2 \\
			&\quad+\frac 14\left((\partial_1\log V_k(0))^2-(\partial_2\log V_k(0))^2\right)(g_k+r)^2\bigg)\cos 2\theta \\
			&\quad+\left(\partial_{12}(\log V_k)(0)\frac{r^2}{2}+\frac{\partial_1V_k(0)}{V_k(0)}\frac{\partial_2V_k(0)}{V_k(0)}\frac{(g_k+r)^2}{2}\right)\sin 2\theta \bigg ). 
		\end{aligned}
	\end{equation}

	Next we define $c_{0,k}$ to be the unique solution of
	\begin{equation}
		\label{c0k}
		\begin{cases}
			\frac{d^2}{dr^2}c_{0,k}+\frac{1}{r}\frac{d}{dr}c_{0,k}+V_k(0) r^{2\alpha}e^{U_k}c_{0,k}\\
			=-\frac{\varepsilon_k^2}4V_k(0)r^{2\alpha}e^{U_k}\bigg ((g_k+r)^2|\nabla \log V_k(0)|^2+\Delta \log V_k(0)r^2\bigg ).\\
			\\
			c_{0,k}(0)=c'_{0,k}(0)=0, \quad 0<r<\tau\varepsilon_k^{-1}.
		\end{cases}
	\end{equation}
	By standard ODE theory $c_{0,k}$ is the unique solution of (\ref{c0k}) on $0<r<\tau\varepsilon_k^{-1}$ and a rough estimate of $c_{0,k}$ shows that,
	\begin{equation}\nonumber
		|c_{0,k}(r)|\le C\varepsilon_k^2\log (1+r),\quad 0<r<\tau\varepsilon_k^{-1}.
	\end{equation}
	
	We will need the following relevant property of $c_{0k}$.
	
	\begin{lem}\label{inte-c0k}
		$$\int_{\partial \Omega_k}\partial_{\nu}c_{0,k}=d_{1,k}\Delta \log V_k(0)\varepsilon_k^2+O(\varepsilon_k^{4+\epsilon_0})$$
		for a small constant $\epsilon_0>0$ and $d_{1,k}$ stated as in the claim of Theorem \ref{int-pos}.
	\end{lem}
	
	\begin{proof}[Proof of Lemma \ref{inte-c0k}.] To better evaluate $c_{0,k}$ we use the following function:
	\begin{equation}\nonumber
		\xi(r)=\frac{1-\frac{V_k(0)}{8(1+\alpha)^2}r^{2+2\alpha}}{1+\frac{V_k(0)}{8(1+\alpha)^2}r^{2+2\alpha}},
	\end{equation}
	
	which solves
	$$\Delta \xi+V_k(0)r^{2\alpha}e^{U_k}\xi=0 \quad \mbox{in}\quad \mathbb R^2.$$
	Let
	\begin{equation}\nonumber
		E_1^k(r)=\frac{\varepsilon_k^2}4V_k(0)r^{2\alpha}e^{U_k}\bigg ((g_k+r)^2|\nabla \log V_k(0)|^2+\Delta \log V_k(0)r^2\bigg )
	\end{equation}
	then integration by parts gives
	\begin{equation}\label{c0nu}\int_{\partial \Omega_k}\partial_{\nu}c_{0,k}\xi-\partial_{\nu}\xi c_{0,k}=-\int_{\Omega_k}E_1^k\xi.
	\end{equation}
	On $\partial \Omega_k=\partial B_{\tau\varepsilon_k^{-1}}(0)$, we have $\xi=-1+O(\varepsilon_k^{2+2\alpha})$ and $\partial_{\nu}\xi=O(\varepsilon_k^{3+2\alpha})$.
	Clearly we also have,
	$$\int_{\partial \Omega_k}\partial_{\nu}c_{0,k}\xi=-\int_{\partial \Omega_k}\partial_{\nu}c_{0,k}+O(\varepsilon_k^{4+\epsilon_0})$$
	and
	$$\int_{\partial \Omega_k}\partial_{\nu}\xi c_{0,k}=O(\varepsilon_k^{4+\epsilon_0}),$$
	for some $\epsilon_0>0$ since on $\partial \Omega_k$ it holds, $$\xi=-1+O(\varepsilon_k^{2+2\alpha})\quad \mbox{and }\quad  \partial_{\nu}c_{0,k}=O(\varepsilon_k^3)\log \varepsilon_k^{-1}.$$
	So (\ref{c0nu}) becomes
	$$\int_{\partial \Omega_k}\partial_{\nu}c_{0,k}=\int_{\Omega_k}E_1^k\xi+O(\varepsilon_k^{4+\epsilon_0}).$$

	At this point we use the following two identities derived in \cite{chen-lin} to estimate the right hand side.  The first one is
	\begin{equation}\label{ident-1}
		\int_0^{\infty}s^{2\alpha+3}e^{U(s)}\frac{(1-s^{2\alpha+2})(\frac{2+\alpha}{\alpha}-s^{2\alpha+2})^2}{(1+s^{2\alpha+2})^3}ds=0,
	\end{equation}
	where $U(x)=-2\log (1+|x|^{2\alpha+2})$. Based on this, for $a_k=\frac{V_k(0)}{8(1+\alpha)^2}$, since $$\frac{s^{2\alpha+3}(\frac{2+\alpha}{\alpha}-s^{2\alpha+2})^2(1-s^{2\alpha+2})}{(1+s^{2\alpha+2})^5}=-s^{-2\alpha-1}+O(s^{-4\alpha-3}),\quad s>C\varepsilon_k^{-1}, $$ then we deduce that,
	\begin{equation}
		\label{ident-3}
		\begin{aligned}
			&\frac{\varepsilon_k^2}4\int_{\Omega_k}V_k(0)r^{2\alpha}e^{U_k}(g_k+r)^2|\nabla \log V_k(0)|^2\xi\\
			&=\frac{\varepsilon_k^2\pi}{2}\frac{|\nabla V_k(0)|^2}{V_k(0)}
			\int_0^{\tau\varepsilon_k^{-1}}\frac{r^{2\alpha+3}(\frac{2+\alpha}{\alpha}-a_kr^{2\alpha+2})^2(1-a_kr^{2\alpha+2})}{(1+a_kr^{2\alpha+2})^5}dr\\
			&=b_k\varepsilon_k^2\int_0^{a_k^{\frac{1}{2\alpha+2}}\tau\varepsilon_k^{-1}}\frac{s^{2\alpha+3}(\frac{2+\alpha}{\alpha}-s^{2\alpha+2})^2(1-s^{2\alpha+2})}{(1+s^{2\alpha+2})^5}ds \\
			&=-b_k\varepsilon_k^2\int_{a_k^{\frac{1}{2\alpha+2}}\tau\varepsilon_k^{-1}}^{\infty}
			\frac{s^{2\alpha+3}(\frac{2+\alpha}{\alpha}-s^{2\alpha+2})^2(1-s^{2\alpha+2})}{(1+s^{2\alpha+2})^5}ds\\
			&=b_{1,k}\varepsilon_k^{2}+O(\varepsilon_k^{4+4\alpha}),
		\end{aligned}
	\end{equation}
	where
	$$b_k=4\pi(1+\alpha)^2|\nabla \log V_k(0)|^2\left(\frac{8(1+\alpha)^2}{V_k(0)}\right)^{\frac{1}{\alpha+1}},$$ and
	$$b_{1,k}=\frac{16\pi(1+\alpha)^4}{\alpha V_k(0)}|\nabla \log V_k(0)|^2\tau^{-2\alpha}\varepsilon_k^{2\alpha}.$$
	
	The second identity reads,
	\begin{equation}\label{ident-2}
		\int_{0}^{\infty}\frac{8(1+\alpha)^2s^{2\alpha+3}}{(1+s^{2(1+\alpha)})^2}\frac{1-s^{2(1+\alpha)}}{1+s^{2(1+\alpha)}}{\rm d}s=-\frac{4\pi}{(1+\alpha)\sin\frac{\pi}{1+\alpha}}.
	\end{equation}
	
	Since
	$$\frac{s^{2\alpha+3}(1-s^{2\alpha+2})}{(1+s^{2\alpha+2})^3}=-s^{-2\alpha-1}+O(s^{-4\alpha-3}),\quad s>C\varepsilon_k^{-1}, $$
	we use (\ref{ident-2}) to obtain
	\begin{equation}
		\label{lead-2}
		\begin{aligned}
			&\frac{\varepsilon_k^2}4\int_{\Omega_k}|y|^{2\alpha}V_k(0)e^{U_k}\Delta (\log V_k)(0)|y|^2\frac{1-a_k|y|^{2\alpha+2}}{1+a_k|y|^{2\alpha+2}}dy\\
			&=\frac{\pi}2\varepsilon_k^2\Delta (\log V_k)(0)\int_0^{\tau \varepsilon_k^{-1}}\frac{V_k(0)r^{2\alpha+3}}{(1+a_kr^{2\alpha+2})^2}\frac{1-a_kr^{2\alpha+2}}{1+a_kr^{2\alpha+2}}dr \\
			&=\frac{\pi}2\varepsilon_k^2\Delta (\log V_k)(0)V_k(0)a_k^{-\frac{2\alpha+4}{2\alpha+2}}\int_0^{a_k^{\frac{1}{2\alpha+2}}\tau\varepsilon_k^{-1}}
			\frac{s^{2\alpha+3}(1-s^{2\alpha+2})}{(1+s^{2\alpha+2})^3}ds \\ 
			&=b_{2,k}\varepsilon_k^2+b_{3,k}\varepsilon_k^2+O(\varepsilon_k^{4+4\alpha})
		\end{aligned}
	\end{equation}
	where
	$$
	b_{2,k}=-\frac{2\pi^2}{(1+\alpha)\sin \frac{\pi}{1+\alpha}}\left(\frac{8(1+\alpha)^2}{V_k(0)}\right)^{\frac 1{1+\alpha}}\Delta (\log V_k)(0)
	$$
	and
	$$
	b_{3,k}=\frac{16\pi(1+\alpha)^4}{\alpha V_k(0)}\tau^{-2\alpha}\varepsilon_k^{2\alpha}\Delta (\log V_k)(0).
	$$
	Therefore, in view of \eqref{ident-3} and \eqref{lead-2} we deduce that,
	\begin{equation}\label{add-exp}
		\int_{\partial \Omega_k}{\partial_\nu c_{0,k}}=d_{1,k} \Delta \log V_k(0)\varepsilon_k^2+O(\varepsilon_k^{4+\epsilon_0})
	\end{equation}
	for some $\epsilon_0>0$, where
	\begin{align*}
		d_{1,k}=&-\frac{2\pi^2}{(1+\alpha)\sin \frac{\pi}{1+\alpha}}\left(\frac{8(1+\alpha)^2}{V_k(0)}\right)^{\frac 1{1+\alpha}}\\
		&+\left(\frac{16\pi(1+\alpha)^2}{\alpha }\frac{|\nabla \log V_k(0)|^2}{V_k(0)\Delta\log V_k(0)}+\frac{16\pi(1+\alpha)^4}{\alpha V_k(0)}\right)\tau^{-2\alpha}\varepsilon_k^{2\alpha}.
	\end{align*}
	Then Lemma \ref{inte-c0k} is established. \end{proof}
	
	\medskip
	
	Now we derive another intermediate estimate by adding one more term in the expansion of $v_k$. This term is $c_{2,k}$ which satisfies
	
	\begin{equation}\label{c2k-e}\Delta c_{2,k}+V_k(0) r^{2\alpha}e^{U_k}c_{2,k}=
		-\varepsilon_k^2V_k(0) r^{2\alpha}e^{U_k}\hat \Theta_2.
	\end{equation}
	\noindent
	Toward this goal we first construct a function $f_{2a}(r)$ that solves
	
	$$f_{2a}''(r)+\frac 1r f_{2a}'(r)+(V_k(0) r^{2\alpha}e^{U_k}-\frac{4}{r^2})f_{2a}(r)=-V_k(0) r^{2\alpha}e^{U_k}E_{21}$$
	for $0<r<\infty$. Similarly we set $f_{2b}$ to solve
	$$ f_{2b}''(r)+\frac 1r f_{2b}'(r)+(V_k(0) r^{2\alpha}e^{U_k}-\frac{4}{r^2})f_{2b}(r)
	=-V_k(0) r^{2\alpha}e^{U_k}E_{22}
	$$
	where $E_{21}$ and $E_{22}$ are defined in (\ref{E212}).
	Let $F_{12}$ and $F_{22}$ be the two fundamental solutions as readily derived from the explicit expressions in (\ref{two-fun}). We have
	$|F_{12}(r)|\sim r^2$ at $0$ and $\infty$, $|F_{22}(r)|\sim r^{-2}$ at $0$ and $\infty$. Then we define $f_{2a}$ as follows,
	$$f_{2a}(r)=F_{12}(r)\int_r^{\infty}\frac{s^{1+2\alpha}e^{U_k}E_{21}F_{22}}{4(1-\frac{4}{(1+\alpha)^2})}ds
	-F_{22}(r)\int_0^r\frac{s^{1+2\alpha}e^{U_k}E_{21}F_{12}}{4(1-\frac{4}{(1+\alpha)^2})}ds $$
	where $E_{2a}$ is the right hand side of the equation of $f_{2a}$. Obviously $f_{2b}$ can be constructed in a similar manner.
	So we have
	\begin{equation}\label{two-fun1}
		|f_{2a}(r)|+|f_{2b}(r)|=O(1+r)^{-2\alpha},\quad 0<r<\varepsilon_k^{-1}
	\end{equation}
	
	and $c_{2,k}$ is defined by
	\begin{equation}\label{c2k}
		c_{2,k}=\varepsilon_k^{2}V_k(0) f_{2a}(r)\cos (2\theta)+\varepsilon_k^{2}V_k(0) f_{2b}(r)\sin (2\theta).
	\end{equation}
	With the definition of $c_{2,k}$ we estimate how $U_k+c_{0,k}+c_{1,k}+c_{2,k}$ approximates $v_k$.
	
	\begin{lem}\label{intermediate-2}
		$$|v_k(y)-(U_k+c_{0,k}+c_{1,k}+c_{2,k})(y)|\le C\varepsilon_k^{2+\epsilon_0}(1+|y|)^{2\epsilon_0}\quad 
		\mbox{in}\quad \Omega_k, $$
		for some $\epsilon_0>0$.
	\end{lem}
	
	\begin{proof}[Proof of Lemma \ref{intermediate-2}.] We use $$w_{2,k}=v_k-(U_k+c_{0,k}+c_{1,k}+c_{2,k}).$$
	Then based on (\ref{Uk-c1k}),(\ref{c0k}),(\ref{c2k-e}) we have
	$$(\Delta+|y|^{2\alpha}V_k(\varepsilon_ky)e^{v_k})w_{2,k}
	=O(\varepsilon_k^3)(1+|y|)^{-1-2\alpha}\quad \mbox{in}\quad \Omega_k. $$
	Obviously the right hand side can be written as
	$O(\varepsilon_k^{2+\epsilon_0}(1+|y|)^{-2-\epsilon_0/2}$ for some $\epsilon_0>0$. The oscillation of $w_{2,k}$ is not zero. But analysing the oscillations of $c_{1,k}$ and $c_{2,k}$ we see that the oscillation of $w_{2,k}$ on $\partial \Omega_k$ is of order $O(\varepsilon_k^{2+2\alpha})$. Therefore, encoding such an oscillation in a suitable harmonic function $\phi_{2,k}$, we deduce that, 
	$$|\phi_{2,k}(y)|\le C\varepsilon_k^{3+2\alpha}|y|.$$
	The construction of $c_{0,k},c_{1,k},c_{2,k}$ implies that $w_{2,k}(0)=0$. Based on the key points mentioned in Remark \ref{stand-rem}, we have
	$$|w_{2,k}-\phi_{2,k}(y)|\le C\varepsilon_k^{2+\epsilon_0}(1+|y|)^{2\epsilon_0}\quad \mbox{in}\quad \Omega_k,$$
	for some $\epsilon_0>0$ small. Since $\phi_{2,k}$ can be considered as a part of the error, Lemma \ref{intermediate-2} is established. \end{proof}
	
	\medskip
	
	Lemma \ref{intermediate-2} will be used in the proof of the main theorem later on. However we still need more terms in the expansion. We shall use $U_k+c_k$ to approximate $v_k$,
	where
	\begin{equation}\label{def-ck-vip}
		c_k=c_{0,k}+\bar c_{0,k}+\sum_{i=1}^4c_{i,k},
	\end{equation}
	and $c_{0,k}$,$c_{1,k}$,$c_{2,k}$ have been already defined. To define the remaining terms we write the expansion of $\log \frac{V_k(\varepsilon_ky)}{V_k(0)}$ as follows,
	\begin{align*}\log \frac{V_k(\varepsilon_ky)}{V_k(0)}=&\sum_{|\beta |=1}^4\varepsilon_k^{|\beta |}\frac{\partial^{\beta}\log V_k(0)}{\beta !}y^{\beta}+O(\varepsilon_k^5)|y|^5\\
		=&\sum_{|\beta |=1}^4(\varepsilon_k^{|\beta |}\Theta_{|\beta |}|y|^{|\beta |}+\varepsilon_k^{|\beta |}\mathcal{R}_{|\beta |})+O(\varepsilon_k^5)|y|^5,
	\end{align*}
	where $\mathcal{R}_{|\beta |}$ stand for the radial parts (thus $\mathcal{R}_1=\mathcal{R}_3=0$),
	while $\Theta_{|\beta |}$ are the collection of projections onto non radial modes. For later convenience let us set:
	\begin{equation}
		\label{t-c1k}
		\begin{aligned}
			\tilde c_{1,k}&=c_{1,k}+\varepsilon_k\nabla \log V_k(0)\cdot \theta r\\
			&=c_{1,k}+\varepsilon_k\nabla \log h_k(0)\cdot \theta r\\
			&=\varepsilon_k\nabla \log  h_k(0)\cdot \theta (r-\frac{2(1+\alpha)}{\alpha}\frac{r}{1+a_k r^{2\alpha+2}}),\quad a_k=\frac{V_k(0)}{8(1+\alpha)^2}, 
		\end{aligned}
	\end{equation}
	and similarly,
	$$\tilde c_{2,k}=c_{2,k}+\varepsilon_k^2\Theta_2|y|^2.$$
	Let us define,
	$$L=\Delta+|y|^{2\alpha}h_k(0)e^{U_k},$$
	then the equation for $c_{3,k}$ takes the form,
	$$Lc_{3,k}=|y|^{2\alpha}h_k(0)e^{U_k}\bigg (\Theta_3r^3\varepsilon_k^3+\tilde c_{1,k}(\tilde c_{2,k}+\mathcal{R}_2\varepsilon_k^2+c_{0,k})+\frac 16 \tilde c_{1,k}^3 \bigg ),$$
	while  for $c_{4,k}$ we have, using
	$$\tilde c_{0,k}=c_{0,k}+\mathcal{R}_2\varepsilon_k^2,$$
	\begin{align*}L c_{4,k}=|y|^{2\alpha}h_k(0)e^{U_k}\bigg (&\Theta_4 r^4 \varepsilon_k^4+(\tilde c_{0,k}\tilde c_{2,k})_A+\frac 12\big( \tilde c_{1,k}^2(\tilde c_{0,k}+\tilde c_{2,k}))_A\\
		&+\frac{(\tilde c_{2,k}^2)_A}{2}+(\tilde c_{1,k}\tilde c_{3,k})_A+\frac 1{24}(\tilde c_{1,k}^4)_A\bigg )
	\end{align*}
	where $(\cdot )_A$ denotes the angular part of the expression in the brackets.
	We define $\bar c_{0,k}$ as follows,
	\begin{equation}
		\label{bar-c0k}
		\begin{aligned}
			L \bar c_{0,k}=h_k(0)|y|^{2\alpha}e^{U_k}\bigg (&\mathcal{R}_4\varepsilon_k^4+(\tilde c_{0,k}\tilde c_{2,k})_r+\frac 12 (\tilde c_{1,k}^2(\tilde c_{0,k}+\tilde c_{2,k}))_r\\
			&+\frac{(\tilde c_{2,k}^2)_r}{2}+(\tilde c_{1,k}\tilde c_{3,k})_r+\frac 12\tilde c_{0,k}^2+\frac 1{24}(\tilde c_{1,k}^4)_r\bigg ) 
		\end{aligned}
	\end{equation}
	where $(\cdot )_r$ means the radial part of the quantity in the brackets. Similar to $\tilde c_{1,k}$ and $\tilde c_{2,k}$, we define
	$$\tilde c_{i,k}=c_{i,k}+\varepsilon_k^ir^i\Theta_i,\quad i=3,4.$$
	Then direct computation leads to
	\begin{equation}\nonumber
		\Delta (U_k+c_k)+|y|^{2\alpha}V_k(\varepsilon_ky)e^{U_k+c_k}=E_w\quad \mbox{in}\quad \Omega_k.
	\end{equation}
	where $E_w=O(\varepsilon_k^5)(1+|y|)^{1-2\alpha}$.
	A  short remark is in order about the computation. The aim of $c_{i,k}$ for $i=1,2,3,4$ is to solve angular terms  proportional to $e^{i\theta}$ whose  order is $\varepsilon_k^i$, for $i=1,2,3,4$,  respectively. The purpose of $c_{0,k}$ is to solve the radial term starting with $\varepsilon_k^2$, the aim of $\bar c_{0,k}$ is to solve the radial term starting with $\varepsilon_k^4$.  Then, by writing $|y|^{2\alpha}V_k(\varepsilon_ky)e^{U_k+c_k}$ as
	\begin{align*}
		|y|^{2\alpha}V_k(\varepsilon_ky)e^{U_k+c_k}=~&|y|^{2\alpha}h_k(0)e^{U_k+c_k+\mathfrak{h}_k}\\
		=~&|y|^{2\alpha}h_k(0)e^{U_k}\left(1+\sum_{l=1}^4\frac{1}{l!}(c_k+\mathfrak{h}_k)^l+ h.o.t\right)
	\end{align*}
	where $\mathfrak{h}_k=\log\frac{V_k(\varepsilon_ky)}{V_k(0)}$,  a lengthy evaluation allows one to identify the terms in each group.
	
	Here we further remark that $c_{l,k}$ for $l=1,2,3,4$ decays as 
	\begin{equation}\label{clk-est}
		|c_{l,k}(y)|\le C\varepsilon_k^l(1+|y|)^{l-2-2\alpha},\quad l=1,2,3,4
	\end{equation} 
	because each $c_{l,k}$ is a finite sum of terms of the form $\varepsilon_k^le^{i l\theta}f_l(r)$ where
	$f_l(r)$ solves
	$$f_l''(r)+\frac 1rf_l'(r)+(r^{2\alpha}h_k(0)e^{U_k}-\frac{l^2}{r^2})f_l(r)=O(1+r)^{l-4-2\alpha},\quad 0<r<\tau\varepsilon_k^{-1}. $$ In particular, in view of \eqref{two-fun}, standard ODE techniques show that
	$$
	f_l(r)=O(1+r)^{l-2-2\alpha},
	$$
	see for example \eqref{two-fun1} above concerning $c_{2,k}$.
	Therefore the oscillation of $c_k$ on $\partial \Omega_k$ is of order  $O(\varepsilon_k^{2+2\alpha})$. If we use $\phi_{1,k}$ to encode the oscillation of $c_k$,
	$$\Delta \phi_{1,k}=0\quad \mbox{in }\quad \Omega_k,\quad \phi_{1,k}=-c_k(y)+\frac{1}{2\pi \tau\varepsilon_k^{-1}}\int_{\partial \Omega_k}c_k,$$
	then we have,
	$$|\phi_{1,k}(y)|\le C\varepsilon_k^{3+2\alpha}|y|\quad \mbox{in}\quad \Omega_k.$$
	In particular by the mean value theorem $\phi_{1,k}(0)=0$ and the first term in the Taylor expansion reads  $a_1^ky_1+a_2^ky_2$ with $|a_1^k|+|a_2^k|\le C\varepsilon_k^{3+2\alpha}$. Let us set $U_k+c_k+\phi_{1,k}$ to approximate $v_k$, then we need to consider the error caused by $\phi_{1,k}$, which at first order takes the form $|y|^{2\alpha}V_k(\varepsilon_k y)e^{U_k+c_k}\phi_{1,k}$. It is readily seen that
	\begin{align*}
		|y|^{2\alpha}h_k(\varepsilon_ky)e^{U_k+c_k}\phi_{1,k}=&-|y|^{2\alpha}h_k(0)e^{U_k}(a_1^ky_1+a_2^ky_2)\\
		&+O(\varepsilon_k^{4+2\alpha})(1+|y|)^{-2-2\alpha}.   
	\end{align*} 
	
	The first term can be eliminated by a separable function $\phi_{2,k}$, which satisfies:
	$$L\phi_{2,k}=|y|^{2\alpha}h_k(0)e^{U_k}(a_1^ky_1+a_2^ky_2),\quad \mbox{in}\quad \Omega_k.$$
	Remark that $\phi_{2,k}$ is separable whence $\phi_{2,k}(0)=0$.
	Also, recalling that $a_i^k=O(\varepsilon_k^{3+2\alpha})$, then we see that,
	$$ \phi_{2,k}(y)=O(\varepsilon_k^{3+2\alpha})(1+r)^{-1-2\alpha}$$ and the oscillation of $\phi_{2,k}$ on $\partial \Omega_k$ is of order $O(\varepsilon_k^{4+4\alpha})$. Higher order terms due to $\phi_{1,k}$ are readily seen to yield a minor contributions of order $O(\varepsilon_k^{6+4\alpha})(1+|y|)^{-2-2\alpha}$.
	
	Next we set $\phi_{3,k}$ to eliminate the term $E_w$,
	\begin{equation*}
		\begin{cases}
			-\Delta\phi_{3,k}=E_w\quad \mbox{in}\quad \Omega_k, \\
			\\
			\phi_{3,k}(0)=0,\quad \phi_{3,k}=constant \,\, \mbox{on}\quad \partial \Omega_k,
		\end{cases}
	\end{equation*}
	which we define as follows,
	$$\phi_{3,k}(y)=\int_{\Omega_k}(G_k(y,\eta)-G_k(0,\eta))E_w(\eta)d\eta .$$
	Then by standard estimates we have,
	\begin{equation*}
		|\phi_{3,k}(y)|\le C\varepsilon_k^5(\log (1+|y|)+(1+|y|)^{3-2\alpha}).
	\end{equation*}
	Recalling the definition of $c_k$ in \eqref{def-ck-vip} and setting
	\begin{equation}
		\label{def-wk-vip}
		w_k=v_k-U_k-c_k-\sum_{i=1}^3\phi_{i,k},
	\end{equation}
	we have
	\begin{align*}
		\Delta w_k+|y|^{2\alpha}h_k(\varepsilon_ky)e^{\xi_k}w_k&=O(\varepsilon_k^5)((1+|y|)^{-1-4\alpha}+O(\varepsilon_k^{4+2\alpha})(1+|y|)^{-2-4\alpha} \\
		&=O(\varepsilon_k^{4+\delta})(1+|y|)^{-2-\delta} 
	\end{align*}
	in $\Omega_k$ for some $\delta>0$. Here we note that the second term on the r.h.s of the first equality comes from $\phi_{2,k}$.  Since $w_k(0)=0$ and $w_k$ has an oscillation of order $O(\varepsilon_k^{4+4\alpha})$ on $\partial \Omega_k$, by Remark \ref{stand-rem},
	we have
	\begin{equation}\label{great-1}
		|w_k(y)|\le C\varepsilon_k^{4+\delta}(1+|y|)^{\delta}.
	\end{equation} 
	At this point, scaling the domain of integration to $\Omega_k$, we evaluate the integral over $B_{\tau}$ as follows,
	\begin{align*}
		&\int_{B_{\tau}}|x|^{2\alpha}h_k(x)e^{u_k(x)}\\
		&=\int_{\Omega_k}|y|^{2\alpha}h_k(0)e^{U_k}(1+\tilde c_{0,k}+\frac 12 (\tilde c_{1,k}^2)_r\\
		&\quad+\bar c_{0,k}+\mathcal{R}_4\varepsilon_k^4+(\tilde c_{0,k}\tilde c_{2,k})_r+\frac 12 (\tilde c_{1,k}^2(\tilde c_{0,k}+\tilde c_{2,k}))_r\\
		&\quad+\frac 12(\tilde c_{2,k}^2)_r+(\tilde c_{1,k}\tilde c_{3,k})_r+\frac 12\tilde c_{0,k}^2+\frac 1{24}(\tilde c_{1,k}^4)_r\bigg )
		+O(\varepsilon_k^{4+\epsilon_0})\\
		&=\int_{\Omega_k}|y|^{2\alpha}h_k(0)e^{U_k}+\int_{\partial \Omega_k}\frac{\partial c_{0,k}}{\partial \nu}+\int_{\partial \Omega_k}\frac{\partial \bar c_{0,k}}{\partial \nu}+O(\varepsilon_k^{4+\epsilon_0}).
	\end{align*}
	Remark that the contribution from $\phi_{1,k}$ vanishes since it is the sum of odd functions. The same observation applies to $\phi_{2,k}$ which is a separable function. Higher order contributions are readily seen to be of order $O(\varepsilon_k^{4+\epsilon_0})$. By a direct computation we have,
	$$\int_{\Omega_k}h_k(0)|y|^{2\alpha}e^{U_k}=8\pi(1+\alpha)\left(1-\frac{e^{-u_k(0)}}{e^{-u_k(0)}
		+\frac{h_k(0)}{8(1+\alpha)^2}\tau^{2\alpha+2}}\right),$$
	which gives the first term in the statement of Theorem \ref{int-pos}. From \eqref {add-exp} and the expression of $d_{1,k}$ right after it we obtain the second term in the statement while for the last term we need a remark about $\int_{\partial \Omega_k}\partial_{\nu}\bar c_{0,k}$. Using (\ref{bar-c0k}) and $\xi$ as a solution of $L\xi=0$, we have, as in the case for $c_{0,k}$, that
	\begin{equation}\nonumber
		\int_{\partial \Omega_k}\partial_{\nu}\bar c_{0,k}=-\int_{\Omega_k}\bar E_{0,k}\xi+O(\varepsilon_k^{4+\epsilon_0})
	\end{equation}
	where $\bar E_{0,k}$ stands for the right hand side of (\ref{bar-c0k}).  The proof is concluded by the identification of $T_1(V_k)\varepsilon_k^4$ with $\int_{\partial \Omega_k}\partial_{\nu}\bar c_{0,k}$. Thus Theorem \ref{int-pos} is established. $\hfill\Box$

	\subsection{Asymptotic analysis around a negative pole}
	In this case we use $\beta\in (-1,0)$ to denote the strength of a negative pole and define,
	$$\displaystyle\tilde \varepsilon_k=e^{-\frac{\lambda^k_1}{2(1+\beta)}}.$$ We use the same notations to approximate $v_k$. The main result in this case takes the form,
	\begin{thm}\label{integral-neg}
		\begin{align*}
			\int_{B_{\tau}}h_k(x)|x|^{2\beta}e^{u_k}
			=~&8\pi(1+\beta)\left(1-\frac{e^{-u_k(0)}}{e^{-u_k(0)}+\frac{h_k(0)}{8(1+\beta)^2}\tau^{2\beta+2}}\right)\\
			&+b_k(\tau,\beta)+\ell_k(\tau,\beta)+O(e^{(-2-\epsilon_0)u_k(0)})
		\end{align*}
		for some $\epsilon_0>0$, where both $b_k(\tau,\beta)$ and $\ell_k(\tau,\beta)$ depend on $|\nabla \log V_k(0)|$ and $\Delta \log V_k(0)$, $b_k(\tau,\beta)\sim e^{-u_k(0)}$, $\ell_k(\tau,\beta)\sim e^{-2u_k(0)}$ for $\tau>0$, and 
		\begin{equation}\label{neg-tau}\lim_{\tau\to 0}\frac{b_k(\tau,\beta)}{e^{-u_k(0)}}=0,\quad
			\lim_{\tau\to 0}\frac{\ell_k(\tau,\beta)}{e^{-2u_k(0)}}=0.
		\end{equation}
	\end{thm}
	
\begin{proof}[Proof of Theorem \ref{integral-neg}.]	The proof of Theorem \ref{integral-neg} is different from that of Theorem \ref{int-pos}, as one can see  from the following proposition:
	
	\begin{prop}\label{thm-3} For fixed $-1<\beta<0$  we have,
		$$|v_k(y)-U_k(y)-c_{1,k}(y)|\le
		C(\delta)\tilde \varepsilon_k^{2+2\beta-\delta}(1+|y|)^{\delta}, \quad y\in \Omega_k$$
		for  some $\delta>0$ small enough, where $c_{1,k}$ is defined as in (\ref{nov10e6}) with $\alpha$ replaced by $\beta$ in the definition of $g_k$ in (\ref{mar27e1}).
	\end{prop}
	
	\begin{proof}[Proof of Proposition \ref{thm-3}.] Direct computation shows that $w_k$ satisfies
	$$\Delta w_k+|y|^{2\beta}e^{V_k(\tilde \varepsilon_ky)}w_k=O(\tilde \varepsilon_k^2)(1+|y|)^{-2-2\beta}$$
	in $\Omega_k$. Since $\beta\in (-1,0)$, the right hand side is bounded above by 
	$$O(\tilde \varepsilon_k^{2+2\beta-\delta})(1+|y|)^{-2-\delta}$$ for a small $\delta>0$. It is easy to see that the oscillation of $c_{1,k}$ is $O(\tilde \varepsilon_k^{2+2\beta})$ on $\partial \Omega_k$. Thus all the key ingredients mentioned in Remark \ref{stand-rem} are verified and the conclusion of Proposition \ref{thm-3} holds. 
\end{proof}
	
	It is readily seen that the bound of the error in this case is different. We use the same notations for $c_{0,k}$, $c_{i,k}$ ($i=1,2,3$) and in this case we set
	\begin{equation}\label{def-wkneg-vip}
		w_k=v_k-U_k-c_k-\sum_{i=1}^3\phi_{i,k},
	\end{equation}
	and follow step by step the same argument to deduce that,
	\begin{equation}\label{great-neg}
		|w_k(y)|\le C(\delta)\tilde \varepsilon_k^{4+4\beta-\delta}(1+|y|)^{\delta}\quad \mbox{in}\quad \Omega_k,
	\end{equation}
	for some small $\delta>0$. We do not provide the details just to avoid repetitions.  By using the expansion of $v_k$ and this estimate we see that the error in the computation of the integral $\int_{B_{\tau}}|x|^{2\beta}h_ke^{u_k}$ is of order $O(\tilde \varepsilon_k^{4+4\beta-\delta})$. Since $\beta$ is negative, some care is needed about the two identities (\ref{ident-1}), (\ref{ident-2}). The estimate which replaces  (\ref{ident-3}) reads,
	\begin{equation}
		\label{ident-n1}
		\begin{aligned}
			&\frac{\tilde \varepsilon_k^2}4\int_{\Omega_k}V_k(0)r^{2\beta}e^{U_k}(g_k+r)^2|\nabla \log V_k(0)|^2\xi\\
			&=\frac{\tilde \varepsilon_k^2\pi}{2}\frac{|\nabla V_k(0)|^2}{V_k(0)}
			\int_0^{\tau\tilde\varepsilon_k^{-1}}\frac{r^{2\beta+3}(\frac{2+\beta}{\beta}-a_kr^{2\beta+2})^2(1-a_kr^{2\beta+2})}{(1+a_kr^{2\beta+2})^5}dr\\
			&=d_k\tilde \varepsilon_k^2\int_0^{a_k^{\frac{1}{2\beta+2}}\tau\tilde\varepsilon_k^{-1}}\frac{s^{2\beta+3}(\frac{2+\beta}{\beta}
				-s^{2\beta+2})^2(1-s^{2\beta+2})}{(1+s^{2\beta+2})^5}ds 
		\end{aligned}
	\end{equation}
	where $\displaystyle{a_k=\frac{V_k(0)}{8(1+\beta)^2}}$ and
	$$d_k=4(1+\beta)^2\pi|\nabla \log V_k(0)|^2\left(\frac{8(1+\beta)^2}{V_k(0)}\right)^{\frac{1}{\beta+1}}.$$
	Here we observe that,
	$$\frac{s^{2\beta+3}(\frac{2+\beta}{\beta}-s^{2\beta+2})^2(1-s^{2\beta+2})}{(1+s^{2\beta+2})^5}=-s^{-2\beta-1}+O(s^{-4\beta-3}),\quad s\to +\infty.$$
	
	The term that corresponds to $\Delta (\log V_k(0))$ in replacement of \eqref{lead-2} takes the form, 
	\begin{equation}
		\label{lead-neg}
		\begin{aligned}
			&\frac{\tilde \varepsilon_k^2}4\int_{\Omega_k}|y|^{2\beta}V_k(0)e^{U_k}\Delta (\log V_k(0))|y|^2\frac{1-a_k|y|^{2\beta+2}}{1+a_k|y|^{2\beta+2}}dy\\
			&=\frac{\pi}2\tilde \varepsilon_k^2\Delta (\log V_k(0))\int_0^{\tau \tilde \varepsilon_k^{-1}}\frac{V_k(0)r^{2\beta+3}}{(1+a_kr^{2\beta+2})^2}\frac{1-a_kr^{2\beta+2}}{1+a_kr^{2\beta+2}}dr \\
			&=\frac{\pi}2\tilde \varepsilon_k^2\Delta (\log V_k(0))V_k(0)a_k^{-\frac{2\beta+4}{2\beta+2}}\int_0^{a_k^{\frac{1}{2\beta+2}}\tau\tilde \varepsilon_k^{-1}}\frac{s^{2\beta+3}(1-s^{2\beta+2})}{(1+s^{2\beta+2})^3}ds \\
			&=\pi\tilde \varepsilon_k^2\frac{\Delta (\log V_k)(0)}{(V_k(0))^{\frac{1}{\beta+1}}}\left(8(1+\beta)^2\right)^{\frac{\beta+2}{\beta+1}}\int_0^{a_k^{\frac{1}{2\beta+2}}\tau\tilde \varepsilon_k^{-1}}\frac{s^{2\beta+3}(1-s^{2\beta+2})}{(1+s^{2\beta+2})^3}ds.
		\end{aligned}
	\end{equation}
	At this point we set $b_k(\tau,\beta)$ to be the sum of the last term of (\ref{ident-n1}) and the last term of (\ref{lead-neg}),
	\begin{align*}
		b_k(\tau,\beta)=~&d_k\tilde \varepsilon_k^2\int_0^{a_k^{\frac{1}{2\beta+2}}\tau\tilde\varepsilon_k^{-1}}\frac{s^{2\beta+3}(\frac{2+\beta}{\beta}
			-s^{2\beta+2})^2(1-s^{2\beta+2})}{(1+s^{2\beta+2})^5}ds \\
		&+\bar d_{k}\tilde \varepsilon_k^2\int_0^{a_k^{\frac{1}{2\beta+2}}\tau\tilde \varepsilon_k^{-1}}\frac{s^{2\beta+3}(1-s^{2\beta+2})}{(1+s^{2\beta+2})^3}ds, 
	\end{align*}
	where
	$$\bar d_{k}=4(1+\beta)^2\pi\Delta (\log V_k)(0)\left(\frac{8(1+\beta)^2}{V_k(0)}\right)^{\frac{1}{1+\beta}}.$$
	It is easy to see that both terms are of the order $O(\tilde \varepsilon_k^{2+2\beta})$ and in particular that  $b_k(\tau,\beta)\sim e^{-u_k(0)}$. If we still use $\bar c_{0,k}$ to denote the radial correction term that corresponds to the radial error of order $O(\tilde \varepsilon_k^{4+4\beta})$, then we identify the $\ell_k(\tau,\beta)$ term in the statement with $\int_{\partial \Omega_k}\partial_{\nu}\bar c_{0,k}$ as in the case $\alpha>0$. Thus (\ref{neg-tau}) can be verified by the definitions of $b_k(\tau,\beta)$ and $\ell_k(\tau,\beta)$. Theorem \ref{integral-neg} is proved. 
\end{proof}

	\subsection{Asymptotic analysis around a regular blowup point}
	In this subsection we analyze the solution near a regular blow up point. After a translation we have that
   $\bar \lambda_k=\max\limits_{B_{2\tau}}u_k=u_k(0)$ denotes the height of the bubble, set $\bar \varepsilon_k=e^{-\bar \lambda_k/2}$ and consider the equation,
	$$\Delta u_k+\bar h_k(x)e^{u_k}=0\quad \mbox{in}\quad B_{2\tau}, $$
    where the ball is now centered at the maximum point.
	We assume in particular that there is only one blow up point in $B_{2\tau}$, that the standard uniform bound
	$\int_{B_{2\tau}}\bar h_ke^{u_k}\le C$ holds and finally that \eqref{u-k-boc} holds, implying
	by well known results (\cite{li-cmp}) that the blow up point is simple.\\
    Let $\phi_R^k$ be the harmonic function defined by the oscillation of $u_k$ on $\partial B_{\tau}$:

    \[\Delta \phi_R^k=0,\quad \mbox{in}\, B_{\tau},\]
    \[
    \phi_R^k(x)=u_k(x)-\frac{1}{2\pi \tau}\int_{\partial B_{\tau}}u_k(y)dS_y,\, x\in \partial B_{\tau}.\]
	
	\noindent
	The main result of this subsection is the following,
	\begin{thm}\label{reg-int} Let $q_k$ denote the local maximum of $u_k-\phi_R^k$. Then 
    \begin{equation}\label{qk:gluck} q_k=\left(-\frac{2}{\bar h_k(0)}\partial_1\phi_R^k(0)\bar\varepsilon_k^2,-\frac{2}{\bar h_k(0)}\partial_2\phi_R^k(0)\bar\varepsilon_k^2\right)+O(\bar\varepsilon_k^3).
    \end{equation}
		Moreover we have, 
        \begin{equation}\label{total-q-1}
			\int_{B_{\tau}}\bar h_k(x)e^{u_k}dx=8\pi-\frac{8\pi\bar \varepsilon_k^2}{\bar\varepsilon_k^2+a_k\tau^2}-\frac{\pi\bar \varepsilon_k^2}{2}\Delta (\log \bar h_k)(0)\bar h_k(0)b_{0,k}+O(\bar\varepsilon_k^{4-\epsilon_0}),
		\end{equation}
		for some $\epsilon_0>0$, where
		\begin{equation*}
			%\label{total-q-2}
			b_{0,k}=\int_0^{\tau\bar\varepsilon_k^{-1}}\frac{r^3(1-a_kr^2)}{(1+a_kr^2)^3}dr,\quad a_k=\bar h_k(0)/8,
		\end{equation*}
      and, for any $\delta>0$ small,
        \begin{equation}\label{new-rate-1}
        \nabla \log \bar h_k(q_k)+\left(1+\frac{4}{\tau\bar h_k(0)}\bar \varepsilon_k^2\right)\nabla \phi_R^k(q_k)=O(\bar\varepsilon_k^{3-\delta}),
        \end{equation}
        \begin{equation}\label{new-rate-2}
        \nabla \Delta \log \bar h_k(q_k)=O(\bar\varepsilon_k^{1-\delta}).
        \end{equation}
	\end{thm}
	\begin{rem}\label{rem:3.2BZY}
	 For later purposes, remark that, even if one integrates on a disk centered at $q_k$ rather than at the maximum of $u_k$, whence slightly different from $B_{\tau}$, the difference would be of order $O(\bar\varepsilon_k^4)$.  Indeed $|q_k|=O(\bar\varepsilon_k^2)$ and then the symmetric difference of these two disks has an area of order $O(\bar\varepsilon_k^2)$, while on this set we also have $e^{u_k}=O(\bar\varepsilon_k^2)$. Thus the integral on the difference of these two disks is of order $O(\bar\varepsilon_k^4)$. As far as \eqref{p_kj-location} holds true as well, a similar observation applies whenever the integral is centered at the the blow up point rather than at the maximum of $u_k$. 
	\end{rem}
	
	\begin{proof} 
Set \[\tilde u_k=u_k-\phi_R^k,\quad \mbox{in}\quad B_{\tau},\] then we cite the estimates in \cite{gluck,zhang1} in the following form. Let 
    \[\tilde v_k(y)=\tilde u_k(q_k+\bar\varepsilon_ky)+2\log \bar\varepsilon_k,\]
     then 
     \[\tilde v_k(y)=-2\log \left(1+\frac{\bar h_k(q_k)}8|y|^2\right)+c_0^k(y)+O\left(\bar\varepsilon_k^2\log \frac{1}{\bar\varepsilon_k}\right)(1+|y|)^{\delta}\]
     for $|y|\le \tau\varepsilon_k^{-1}$, $\delta>0$ small, and
     \[|c_0^k(y)|\le C\bar\varepsilon_k^2\frac{|y|^2(\log (2+|y|))^2}{(1+|y|)^2}.\]
  By standard elliptic estimates we have, 
     \begin{equation}\label{for-qk}
     \left|\nabla \left(\tilde v_k+2\log \left(1+\frac{\bar h_k(q_k)}8|y|^2\right)-c_0^k(y)\right)\right|\le C\bar\varepsilon_k^2\log \frac{1}{\bar\varepsilon_k},
     \end{equation}
     for $|y|\sim 1$. Note here that the image of $0$ after scaling becomes $p_k=-\frac{q_k}{\bar\varepsilon_k}$, which satisfies (\cite{gluck,zhang1}) $|p_k|\le C\bar\varepsilon_k \log \frac{1}{\bar\varepsilon_k}$. This estimate about $p_k$ implies that $|\nabla c_0^k(p_k)|=O(\bar\varepsilon_k^3\log \frac{1}{\bar\varepsilon_k})$.
     Evaluating (\ref{for-qk}) at $p_k$ and using the fact that $\nabla u_k(0)=0$, it is not difficult to see that, 
     \[p_k=\left(\frac{2}{\bar h_k(q_k)}\partial_1\phi_R^k(0)\bar\varepsilon_k,\frac{2}{\bar h_k(q_k)}\partial_2\phi_R^k(0)\bar\varepsilon_k\right)+O(\bar\varepsilon_k^2),\]
     showing indeed that (\ref{qk:gluck}) holds true. 
    
    Now we have
	$\tilde u_k=\mbox{constant}$ on $\partial B_{\tau}$
	and (\cite{gluck,zhang1})
	\begin{equation}\label{gradhk:gluck}
		|\nabla \log \tilde h(0)|=O(\bar \varepsilon_k^2 \log(\bar \varepsilon_k^{-1})),
	\end{equation}
where 
\begin{equation}\label{new-h-1}
\log \tilde h_k(x)=\log \bar h_k(q_k+x)+\phi_R^{k}(q_k+x),\quad q_k+x\in B_{\tau}.
\end{equation}
At this point let us define,
	\begin{equation}\label{bar-v-k}
		\bar v_k(y)=\tilde u_k(q_k+\bar \varepsilon_ky)+2\log \bar \varepsilon_k, 
	\end{equation}
	then we have,
	\begin{equation}\label{R-vk}
		\Delta \bar v_k(y)+\tilde h_k(\bar \varepsilon_ky)e^{\bar v_k}=0
	\end{equation}
    in 
    \[\Omega_k:=\{y;\quad q_k+\bar\varepsilon_ky\in B(0,\tau)\},\]
	and $\bar v_k=\mbox{constant}$ on $\partial \Omega_k$. Therefore we choose
	the first term in the approximation of $\bar v_k$ to be essentially the one suggested in \cite{zhang1} and then in
	\cite{gluck}, which takes the form,
	$$\bar U_k(y)=-2\log \left(1+\frac{\tilde h_k(0)}{8}|y|^2\right).$$
	\begin{rem}
		The $\bar U_k$ suggested in \cite{zhang1} and then in \cite{gluck} is slightly different in the fact that $\tilde h_k$ is not evaluated at $q_k$ but at the maximum point before the subtraction of the harmonic function encoding
		the oscillations. However, by using \eqref{qk:gluck}, \eqref{gradhk:gluck},
		it is easy to see that Theorem 1.1 in \cite{gluck} hold with our definition as well.
	\end{rem}
Since we need to consider round domains, we focus on $B(0,\tau_1^k\bar\varepsilon_k^{-1})$, that is, we choose as a base domain
$B(q_k,\tau_1^k)\subset B(0,\tau)$ in such a way that they are tangent, i.e.  
\begin{equation}\label{two-taus}
\tau_1^k=\tau-|q_k|.
\end{equation}
 To add more terms in the approximation of $\bar v_k$ we first write the equation of $\bar v_k$ as follows,
 \begin{equation}
 \label{bar-vk-2}
 \begin{aligned}
 \Delta \bar v_k(y)+\tilde h_k(0)e^{\bar v_k}
 =&(\tilde h_k(0)-\tilde h_k(\bar\varepsilon_ky))e^{\bar v_k}\\
 =&(\tilde h_k(0)-\tilde h_k(\bar\varepsilon_ky))e^{\bar U_k}\\
 &+O(\bar\varepsilon_k^{4-\delta})(1+|y|)^{-2-\delta/2}
 \end{aligned} \quad \mbox{in}\quad \Omega_{1,k}
 \end{equation}
 where $\delta>0$ is an arbitrarily small number,
 $\Omega_{1,k}=B(0,\tau_1^k\bar\varepsilon_k^{-1})$, and the size of the error term is estimated via  the vanishing rate of the first derivative (\ref{gradhk:gluck}) and Remark 1.5 in \cite{zhang1} which gives 
 \begin{equation}\label{rough-1}\bar v_k-\bar U_k= 
 O(\bar \varepsilon_k^2 (\log(\bar \varepsilon_k^{-1}))^2).
 \end{equation}
 Next we shall derive a first estimate of $\bar v_k$. For this purpose we write 
 \begin{align*} 
 &(\tilde h_k(0)-\tilde h_k(\bar\varepsilon_ky))e^{\bar U_k}\\
 &=\sum_{l=0}^{\infty}(h_{l,k}(r)\cos(l\theta)+\bar h_{l,k}(r)\sin (l\theta)),
 \end{align*}
 where for $l=0,1,...$,
 \[h_{l,k}(r)=\frac 1{\pi}\int_0^{2\pi}(\tilde h_k(0)-\tilde h_k(\bar\varepsilon_kre^{i\theta}))e^{\bar U_k}\cos (l\theta)d\theta, \]
 \[\bar h_{l,k}(r)=\frac 1{\pi}\int_0^{2\pi}(\tilde h_k(0)-\tilde h_k(\bar\varepsilon_k r e^{i\theta}))e^{\bar U_k}\sin (l\theta)d\theta. \]
 Since $l=1$ is of particular importance, we identify the leading terms:
 \begin{align*}
    h_{1,k}(r)=-\bar\varepsilon_k\partial_1 \tilde h_k(0)re^{\bar U_k}-\frac 18\bar\varepsilon_k^3\partial_1\Delta \tilde h_k(0)r^3e^{\bar U_k}+\bar\varepsilon_k^5H_{1,k}(r),\\
    \bar h_{1,k}(r)=-\bar\varepsilon_k\partial_2 \tilde h_k(0)re^{\bar U_k}-\frac 18\bar\varepsilon_k^3\partial_2\Delta \tilde h_k(0)r^3e^{\bar U_k}+\bar\varepsilon_k^5\bar H_{1,k}(r),
 \end{align*}
 and 
 \begin{equation}\label{Handg}|H_{1,k}(r)|+|\bar H_{1,k}(r)|\le Cr^5(1+r)^{-4}.
 \end{equation}
 We set
 $c_{l,k}=f_{l,k}(r)\cos (l\theta)+\bar f_{l,k}(r)\sin (l\theta)$ where
 $f_{l,k}$ is the solution of 
 $$\left(\frac{d^2}{dr^2}+\frac 1r \frac{d}{dr}+\tilde h_k(0)e^{\bar U_k}-\frac{l^2}{r^2}\right)f_{l,k}(r)=h_{l,k}(r),\quad 0<r<\tau_1^k \bar \varepsilon_k^{-1},$$
 which satisfies $f_{l,k}(0)=\frac{d}{dr}f_{l,k}(0)=0$.  The definition of $\bar f_{l,k}$ is similar. 
Before dealing with $l=1$ we provide the estimates for all other terms.  The $c_{0,k}$ solves
	\begin{align*}
		&\frac{d^2}{dr^2}c_{0,k}+\frac 1r\frac{d}{dr}c_{0,k}+\tilde h_k(0)e^{\bar U_k}c_{0,k}\\
		&=(-\frac 14\bar \varepsilon_k^2 \Delta \log \tilde h_k(0)\tilde h_k(0)r^2+O(\bar \varepsilon_k^4r^4))e^{\bar U_k(r)},\quad 0<r<\tau_1^k\bar \varepsilon_k^{-1}, 
	\end{align*}
	with $c_{0,k}(0)=c'_{0,k}(0)=0$, see the expansion right after Remark \ref{stand-rem} for more details.
	A rough estimate of $c_{0,k}$ takes the form,
	$$|c_{0,k}(r)|\le C\bar \varepsilon_k^2 (\log (1+r))^2.$$

For $l>1$	We define $c_{l,k}$ similar to the case with singular source. 
For $l>1$, two fundamental solutions of the homogeneous equation are (for convenience we assume $\tilde h_k(0)=8$)
\[F_{1,l}=\frac{(l+1)r^l+(l-1)r^{l+2}}{1+r^2},\quad F_{2,l}=\frac{(l+1)r^{2-l}+(l-1)r^{-l}}{1+r^2},\]
and we define $f_{l,k}$ as follows, 
\begin{equation}\label{e-l}
f_{l,k}(r)=F_{1,l}(r)\int_r^{\tau_1^k \varepsilon_k^{-1}}\frac{W_1(s)}{W(s)}ds+F_{2,l}(r)\int_0^r\frac{W_2(s)}{W(s)}ds
\end{equation}
where 
\[W(s)=\frac{-2l(l^2-1)}{s},\quad W_1(s)=h_{l,k}(s)F_{2,l}(s),\quad W_2(s)=F_{1,l}(s)h_{l,k}(s).\]
For $l\ge 4$ we just use 
\[|h_{l,k}(r)|\le C\bar \varepsilon_k^4r^4(1+r)^{-4}\] for some $C>0$ independent of $k$ and $l$. Then an elementary estimate gives, 
\begin{equation}\label{est-flk}
|f_{l,k}(r)|\le \frac{C}{l^2}\bar \varepsilon_k^4r^4(1+r)^{-2}, \quad l \ge 4.
\end{equation}

 For $l=1$ two fundamental solutions are
 \[f_1(r)=\frac{r}{1+r^2},\quad f_2(r)=-\frac{1}{2r(1+r^2)}+\frac{2r\log r}{1+r^2}+\frac{r^3}{2(1+r^2)}\]
 and we set
 \[f_{1,k}(r)=f_1(r)\int_0^r\frac{W_1(s)}{W(s)}ds+f_2(r)\int_0^r\frac{W_2(s)}{W(s)}ds\]
 where
 \[W(s)=\frac 1s,\quad W_1(s)=-h_{1,k}(s)f_2(s),  \quad W_2(s)=h_{1,k}(s)f_1(s).\]

 Then $c_{1,k}$ satisfies
 \[c_{1,k}(0)=|\nabla c_{1,k}(0)|=0,\quad |c_{1,k}(y)|\le C\bar\varepsilon_k^3\log \frac{1}{\bar\varepsilon_k}\frac{r^2}{1+r}.\]

Let 
\[c_k=\sum_{l=0}^{\infty}c_{l,k},\] then the estimate (\ref{est-flk}) guarantees that $ c_k$ is well defined. Next we set 
\[\phi_k=\sum_{l=1}^\infty\phi_{l,k}+\phi_E^k,\]
where the sum eliminates the oscillation of $-c_k$ while 
$\phi_E^k$ eliminates the oscillation of 
$\bar v_k$ on $\partial \Omega_{1,k}$ respectively.
As far as $l\ge 4$, by (\ref{est-flk}) we have,
\[|\phi_{l,k}(y)|\le \frac{C}{l^2}\bar\varepsilon_k^{4}r^2,\]
showing that $\phi_k$ is well defined. The overall estimate for $\phi_k$ is 
\[|\phi_k(y)|\le C\bar\varepsilon_k^3\log \frac{1}{\bar\varepsilon_k}|y|.\]
 Let $w_k=\bar v_k-(\bar U_k+ { c_k}+\phi_k)$, then 
 $$\Delta w_k+\tilde h_k(0)e^{\xi_k}w_k=E_k,\quad |y|\le \tau_1^k \bar \varepsilon_k^{-1},$$
 where $\xi_k$ comes from the mean value theorem, which satisfies 
 \[\xi_k(y)\to -2\log (1+|y|^2) \quad \mbox{in } \, C^2_{loc}(\mathbb R^2), \]
 \begin{align*}
 \begin{cases}
 E_k=-\tilde h_k(\bar \varepsilon_ky)e^{\bar U_k}( c_{1,k} +\phi_k)
 +O(\bar\varepsilon_k^{4-\delta_0})(1+|y|)^{-2-\delta_0/2},\\
 \\
 w_k(0)=0,\quad |\nabla w_k(0)|\le C\bar \varepsilon_k^3\log\frac{1}{\bar \varepsilon_k},\quad w_k=\mbox{ constant on }\,\, \Omega_{1,k}.
 \end{cases}
 \end{align*}

Now we claim that, for any $\delta>0$ small enough,
 \begin{equation}\label{wk-1}
 |w_k(y)|\le C\bar\varepsilon_k^3\log \frac{1}{\bar \varepsilon_k}(1+|y|)^{\delta},\quad 
 y\in \Omega_{1,k}
 \end{equation}
 implying by standard elliptic estimates that,
 \begin{equation}\label{wk-2}
 |\nabla w_k(y)|\sim \bar \varepsilon_k^3\log \frac{1}{\bar \varepsilon_k}|y|^{\delta-1},\quad 1<|y|<\frac{\tau}2\bar \varepsilon_k^{-1}.
 \end{equation}

The proof of (\ref{wk-1}) is similar to the argument used above via the Green representation formula and we sketch it here for reader's convenience. Argue by contradiction and, setting 
\[\Lambda_k=\max_{y\in \bar \Omega_{1,k}}\frac{|w_k(y)|}{\bar\varepsilon_k^3\log \frac{1}{\bar \varepsilon_k}(1+|y|)^{\delta}},\]
assume that $\Lambda_k\to \infty$.
Next, let $y_k$ be any maximum point where $\Lambda_k$ is attained and set
\[\hat w_k(y)=\frac{w_k(y)}{\Lambda_k\bar\varepsilon_k^3\log\frac{1}{\bar\varepsilon_k}(1+|y_k|)^{\delta}}.\] Thus we have that, 
 \[|\hat w_k(y)|\le \frac{(1+|y|)^{\delta}}{(1+|y_k|)^{\delta}},\]
 \[\Delta \hat w_k(y)+h_k(q_k)e^{\xi_k}\hat w_k(y)=o(1)\frac{(1+|y|)^{-3}}{(1+|y_k|)^{\delta}},\quad y\in \Omega_{1,k},\]
 and $\hat w_k(0)=0$, $\nabla \hat w_k(0)=o(1)$, and $\hat w_k=$constant on $\partial \Omega_{1,k}$.  Then it is easy to see that $|y_k|\to \infty$ and there is no way for 
 $|\hat w_k(y_k)|=1$ to hold. This contradiction proves (\ref{wk-1}).

With the aid of (\ref{wk-1}) we refine the expansion of $\bar v_k$ as follows, let 
\[w_{1,k}=\bar v_k-\bar U_k.\]
Then by (\ref{wk-1}) we have
\[\Delta w_{1,k}+\tilde h_k(0)e^{\bar U_k}w_{1,k}=(\tilde h_k(0)-\tilde h_k(\bar\varepsilon_ky))e^{\bar U_k}+O(\bar\varepsilon_k^{4-\delta})(1+|y|)^{-2-\delta/2},\]
where we used the fact that, \[O(|w_{1,k}|^2)e^{\bar U_k}=O(\bar\varepsilon_k^{4-\delta})(1+|y|)^{-2-\delta}.\]
Now we abuse the notation $c_k=\sum_{l=0}^{\infty}c_{l,k}$ by defining it as the projection on $e^{il\theta}$: For $l\ge 1$,
\[c_{l,k,1}(r)=\frac{1}{2\pi}\int_0^{2\pi}\bar v_k(re^{i\theta})\cos l\theta d\theta=\frac{1}{2\pi}\int_0^{2\pi} w_{1,k}(re^{il \theta})\cos l\theta d\theta.\]
\[c_{l,k,2}(r)=\frac{1}{2\pi}\int_0^{2\pi}\bar v_k(re^{i\theta})\sin l\theta d\theta=\frac{1}{2\pi}\int_0^{2\pi} w_{1,k}(re^{i\theta})\sin l\theta d\theta.\]

Since the projection of $e^{\bar v_k}$ on $e^{i\theta}$ is particular important, we emphasize that  

\begin{equation}\label{e-c1k1}
\begin{aligned}
&\left (\frac{d^2}{dr^2}+\frac 1{r}\frac{d}{dr}+\tilde h_k(0)e^{\bar U_k}-\frac 1{r^2}\right )c_{1,k,1}\\
&=-\bar \varepsilon_k\partial_1 \tilde h_k(0)re^{\bar U_k}-\frac 18\bar\varepsilon_k^3\partial_1\Delta \tilde h_k(0)r^3e^{\bar U_k}+\bar\varepsilon_k^5H_{1,k}\\&\quad+O(\bar\varepsilon_k^{4-\delta})(1+r)^{-3+\delta},
\end{aligned}
\end{equation}
{where $H_{1,k}$ satisfies \eqref{Handg}}. 
Of course we have a corresponding equation for $c_{1,k,2}$. Here we make two important observations. First of all we have $\bar v_k(0)=|\nabla \bar v_k(0)|=0$, secondly, we recall that $\bar v_k$ is constant on $\partial \Omega_k$. The first observation implies that,
\[c_{1,k,1}(0)=\frac{d}{dr}c_{1,k,1}(0)=0,\quad c_{1,k,2}(0)=\frac{d}{dr}c_{1,k,2}(0)=0,\]
Now we employ Lemma \ref{osc-1} in the appendix to deduce that, 
\begin{equation}\label{new-rel-3}c_{1,k,1}( \tau_1^k\bar\varepsilon_k^{-1})=-2\frac{q_1^k}{\tau}+O(\bar\varepsilon_k^3),\quad c_{1,k,2}( \tau_1^k\bar\varepsilon_k^{-1})=-2\frac{q_2^k}{\tau}+O(\bar\varepsilon_k^3).
\end{equation}

Let $E_{1,k}$ be the right hand side of (\ref{e-c1k1}), then we can compare $c_{1,k,1}$ with
\[T_k(r):=-f_1(r)\int_0^r E_{1,k}(s)f_2(s)sds+f_2(r)\int_0^rf_1(s)E_{1,k}(s)sds,\]
where we recall that two fundamental solutions are
\[f_1(r)=\frac{r}{1+r^2},\quad f_2(r)=-\frac{1}{2r(1+r^2)}+\frac{2r\log r}{1+r^2}+\frac 12 \frac{r^3}{1+r^2}.\]
Putting $w(r)=c_{1,k,1}(r)-T_k(r)$, we see that $w(r)$ satisfies,
\[\begin{cases}
w''(r)+\frac 1rw'(r)+(8e^{\bar U_k}-\frac{1}{r^2})w(r)=0,\quad 0<r<\tau\bar\epsilon_k^{-1},\\
w(0)=w'(0)=0.
\end{cases}\]
Since $w=c_1f_1+c_2f_2$, we see that $w\equiv 0$. Therefore
\[c_{1,k,1}(r)=-f_1(r)\int_0^r E_{1,k}(s)f_2(s)sds+f_2(r)\int_0^rf_1(s)E_{1,k}(s)sds.\]
By an elementary evaluation we have that, 
\begin{equation}\label{new-rel-1}
\begin{aligned}
c_{1,k,1}(r)&=-\frac 18\bar\varepsilon_k\partial_1\tilde h_k(0)r+O(r^{-1}\log r)\bar\varepsilon_k\partial_1\tilde h_k(0)\\
&-\frac 18\partial_1\Delta \tilde h_k(0)\bar\varepsilon_k^3(\frac{r}{2}\log r-\frac{5}8r)+\bar\varepsilon_k^5C_H(r), 
\end{aligned}
\end{equation}
where 
\begin{equation}\label{new-rel-2}C_H(r)=-f_1(r)\int_0^r f_2(s)H_{1,k}(s)sds+f_2(r)\int_0^r sf_1(s)H_{1,k}(s)ds
\end{equation}
and $C_H(r)$ satisfies
\[|C_H(r)|\le Cr^3,\quad \mbox{for some }C \mbox{  independent by } k.\]
Thus we deduce from \eqref{new-rel-3} that,
\begin{equation}\label{important-1}
-\frac {\tau}8\partial_1\tilde h_k(0)-\frac {\tau}{16}\partial_1\Delta \tilde h_k(0)\bar\varepsilon_k^2(\log \frac{1}{\bar\varepsilon_k}-\frac 54)+\bar\varepsilon_k^5C_H(\tau\bar\varepsilon_k^{-1})=O(\bar\varepsilon_k^2).
\end{equation}
Obviously for {$c_{1,k,2}$} defined by
\begin{equation}\label{new-rel-8} c_{1,k,2}(r)=-f_1(r)\int_0^r f_2(s)\bar H_{1,k}(s)sds+f_2(r)\int_0^r sf_1(s)\bar H_{1,k}(s)ds
\end{equation}
we have 
\begin{equation}\label{important-2}
-\frac {\tau}8\partial_2\tilde h_k(0)-\frac {\tau}{16}\partial_2\Delta \tilde h_k(0)\bar\varepsilon_k^2(\log \frac{1}{\bar\varepsilon_k}-\frac 54)+\bar\varepsilon_k^5\bar C_H(\tau\bar\varepsilon_k^{-1})=O(\bar\varepsilon_k^2).
\end{equation}
The analysis for any $c_{l,k}$ can be carried out similarly. By the uniqueness the expression of $c_{l,k}$ can be also be written by (\ref{e-l}), which gives
\[\nabla c_{l,k}(0)=0,\quad |c_{l,k}|\le \frac{C}{l^2}\bar\varepsilon_k^lr^4(1+r)^{-2},\quad l\ge 4.\]

 The harmonic function  we use  to eliminate the oscillation of $\bar v_k$ is $\phi_E^k$, by Lemma \ref{osc-1} we know that the leading term of $\phi_E^k$ is
\begin{equation}\label{phi-2-k}
\phi_2^k(\bar\varepsilon_ky):=-2\frac{q_1^k}{\tau}\bar\varepsilon_kr\cos\theta-2\frac{q_2^k}{\tau}\bar\varepsilon_kr\sin\theta+O(\bar\varepsilon_k^{4-\delta})|y|,
\end{equation}
 In other words, $\phi_E^k-\phi_2^k(\bar\varepsilon_k\cdot)$ is orthogonal to $e^{i\theta}$.

Here we describe the location of the maximum point of 
\[\hat v_k(y):=\bar v_k(y)-\phi_E^k(\bar\varepsilon_ky).\] 
In order to do this we first remark that now we have 
\begin{equation}\label{imp-est}
\bar v_k(y)=\bar U_k+c_k+O(\bar\varepsilon_k^{4-\delta})(1+|y|)^{\delta},
\end{equation} 
because, after incorporating the boundary information of $\bar v_k$ on $\Omega_{1,k}$ and using the uniqueness of the ordinary differential equations that $c_{l,k}$ satisfies, we see that there is no oscillation of $\bar v_k-\bar U_k-c_k$. Indeed, going back to (\ref{bar-vk-2}) and using $\tilde w_k$ to denote $\bar v_k-\bar U_k-c_k$, we see that $\tilde w_k$ satisfies
\[\Delta \tilde w_k+\tilde h_k(0)e^{\xi_k}\tilde w_k=O(\bar\varepsilon_k^{4-\delta})(1+|y|)^{-2-\delta/2},\]
in $\Omega_{1,k}$, $\tilde w_k(0)=|\nabla \tilde w_k(0)|=0$ and $\tilde w_k=$constant on $\partial \Omega_{1,k}$. Thus it is easy to see that 
\[\tilde w_k(y)=O(\bar\varepsilon_k^{4-\delta/2})(1+|y|)^{\delta/2}\quad \mbox{in}\quad \Omega_{1,k},\]
which is equivalent to (\ref{imp-est}).

Using the expression of $\bar U_k$ and the estimates of $c_k$ we see that $D^2(\bar U_k+c_k)(0)\le -CI$ for some $C>0$, where $I$ is the identity matrix. Thus the maximum point $\hat p_k$ of $\hat v_k$ satisfies,
\[|\hat p_k|\le C\bar\varepsilon_k^3\log \frac{1}{\bar\varepsilon_k}.\]

If we use $\hat q_k$ to denote the pre-image of $\hat p_k$ before scaling, we have $|\hat q_k|=O(\bar\varepsilon_k^4\log \frac{1}{\bar\varepsilon_k})$. 

After obtaining this estimate about $|\hat p_k|$, we will only use the form of $c_{1,k}$. Instead we remove the oscillation of $\bar v_k$ on $\partial \Omega_{1,k}$ and define
\[\hat v_k(y)=v_k(\hat p_k+y)-\phi_E^k(\hat p_k+\bar\varepsilon_ky),\quad \hat p_k+\bar\varepsilon_ky\in \Omega_{1,k} \]

\begin{equation}\label{def-hat-h}\log \hat h_k(y)=\log \tilde h_k(y)+\phi_2^k(y)+(\phi_E^k-\phi_2^k)(y).
\end{equation}
The reason we write $\phi_E^k=\phi_2^k+(\phi_E^k-\phi_2^k)$ is because when we calculate the first derivative of $\log \hat h_k(0)$, we have $\nabla (\phi_E^k-\phi_2^k)(0)=0$ and $\phi_2^k$ is the most important term of $\phi_E^k$.
Now we define $\hat v_k$ on $\Omega_{2,k}:=B(0,\tau_2^k\bar\varepsilon_k^{-1})$,
where $B(\hat q_k,\tau_2^k)\subset B(q_k,\tau_1^k)$ is such that they are tangent. 

Using Lemma \ref{osc-1} we see that the oscillation on $\partial \Omega_2^k$ is now $O(\bar \varepsilon_k^3)$. Now we write the equation of $\hat v_k$ as
\[\Delta \hat v_k+\hat h_k(\bar\varepsilon_ky)e^{\hat v_k}=0,\quad \mbox{in}\quad \Omega_{2,k}.\]
Corresponding to this setting we define $\hat U_k$, $\hat c_k$ by analogy as above and we see that there is no oscillation of $\hat v_k-\hat c_k$ on $\partial \Omega_{2,k}$. Thus, by arguing as in \eqref{wk-1},\eqref{wk-2}, we have, 
\begin{equation}\label{wk-new-e}
|\nabla^j(\hat v_k-\hat U_k-\hat c_k)|\le C\bar \varepsilon_k^{4-\delta}(1+|y|)^{\delta-j} 
\end{equation}
for $j=0,1$ and $|y|\le \bar\varepsilon_k^{-1}/2$.

%\begin{lem}\label{new-lem}
%Let $\phi$ solve 
%\[\Delta \phi+8e^{U}\phi=0,\quad \mbox{in}\quad \mathbb R^2, \]
%where $e^{U}=(1+|x|^2)^{-2}$, and $\phi$ satisfies
%\[\phi(0)=0,\quad \nabla \phi(0)=0,\quad |\phi(x)|\le C(1+|x|)^{2-\delta},\]
%for some $\delta>0$. Then $\phi\equiv 0$. 
%\end{lem}

%\begin{proof}[Proof of Lemma \ref{new-lem}.]
%We just consider the projections on $e^{i\theta}$ and prove that they are all zero. Let $k$ be any nonnegative integer, we set 
%\[\phi_k(r)=\frac 1{2\pi}\int_0^{2\pi}\phi(re^{i\theta})\cos k\theta d\theta. \]
%Then $\phi_k$ satisfies the following equation
%\[\phi_k''(r)+\frac 1r\phi_k'(r)+(8e^U-\frac{k^2}{r^2})\phi_k(r)=0.\]
%and we also have 
%$\phi_k(0)=0$,\quad $\phi_k'(0)=0$ and $|\phi_k(r)|\le C(1+r)^{2-\delta}$. 
%When $k=0$, let $g_1,g_2$ be two fundamental solutions and 
%$g_1(r)=\frac{1-r^2}{1+r^2}$, %$g_2(r)=O(\log \frac 1r)$ near $0$. Then it is easy to see that $\phi_0=0$. For $k=1$, let $g_1$ and $g_2$ be two fundamental solutions such that 
%$g_1=\frac{r}{1+r^2}$ and $g_2=O(r^{-1})$ near $0$. 
%Using the zero initial condition and the zero derivative at $0$ one has $\phi_1=0$. For $l\ge 2$, two fundamental solutions are of the order $O(r^l)$ and $o(r^{-l})$ at infinity. By the assumption on the asymptotic behavior at infinity we rule out one of them. Then $\phi_l(0)=0$ rules out the other one. Thus 
%$\phi_k\equiv 0$ for all $k=0,1,2,...$. The projections on $\sin k\theta$ are obviously $0$ by the same reason.  
%Lemma \ref{new-lem} is established. 
%\end{proof}

We are now in a position to obtain the vanishing rates of $|\nabla \hat h_k(0)|$ and $|\nabla \Delta \hat h_k(0)|$. Let  $\Omega_{\sigma,k}=B(0,\bar\varepsilon_k^{-\sigma})$. Then we consider the following Pohozaev identity on $\Omega_{\sigma,k}$,
for any fixed $\xi\in \mathbb S^1$:
\begin{equation}
\label{pi-sigma}
\begin{aligned}
\bar\varepsilon_k\int_{\Omega_{\sigma,k}}\partial_{\xi}\hat h_k(\bar\varepsilon_ky)e^{\hat v_k}=&\int_{\partial \Omega_{\sigma,k}}(\hat h_k(\bar\varepsilon_ky)e^{\hat v_k}(\xi,\nu)+\partial_{\nu}\hat v_k\partial_{\xi}\hat v_k)\\
&-\frac 12\int_{\partial \Omega_{\sigma,k}} |\nabla \hat v_k|^2(\xi\cdot \nu). 
\end{aligned}
\end{equation}
To evaluate the left hand side, after using the expansion (\ref{wk-new-e}) and the current (\cite{gluck,zhang1}) vanishing rate of $|\nabla \hat h_k(0)|=O(\bar\varepsilon_k^2\log{\bar\varepsilon_k}^{-1})$, and taking $\sigma$ close to $\frac 12$, we have 
\[LHS=8\pi\bar\varepsilon_k\partial_{\xi}\log \hat h_k(0)+\frac{16\pi}{\hat h_k^2(0)}\partial_{\xi}\Delta \hat h_k(0)\bar\varepsilon_k^3(\log \bar\varepsilon_k^{-2\sigma}+O(1))+O(\bar\varepsilon_k^{4-\delta}).\]
Note that if use $H_{\xi}^k$ to denote all the high order terms in the expansion of $\partial_{\xi} \hat h_k$:
\[\partial_{\xi}\hat h_k(x)=\partial_{\xi}\hat h_k(0)+\sum_{j=1}^3\sum_{|\alpha |=j}\frac{1}{\alpha !}\partial^{\alpha}\hat h_k(0)x^{\alpha}+H_{\xi}^k(x),\]
we can choose $\sigma$ close to $\frac 12$ so the integration with $H_{\xi}^k$ is minor compared to the error $O(\bar\varepsilon_k^{4-\delta})$.  
Using the expansion of $v_k$ to evaluate the right hand side with $\sigma$ close to $\frac 12$, we have
\[RHS=C_1\bar\varepsilon_k^3\partial_{\xi}\Delta \hat h_k(0)+C_2\bar\varepsilon_k\partial_{\xi}\hat h_k(0)+O(\bar\varepsilon_k^{4-\delta}).\]
Here we note that in the evaluation of the right hand side, the first term $\int_{\partial \Omega_{\sigma,k}}\hat h_ke^{v_k}(\xi\cdot \nu)$ gives
\[\int_{\partial \Omega_{\sigma,k}}\hat h_ke^{v_k}(\xi\cdot \nu)=C_3\bar\varepsilon_k^3\partial_{\xi}\Delta \hat h_k(0)+O(\bar\varepsilon_k^{4-\delta}),\]
 while in the last two boundary terms, because of cancellations, the only term that needs to be evaluated is 
\[\int_{\partial \Omega_{\sigma,k}}\bar U_k'(r)\partial_{\xi}c_{1,k}.\]
Thus the combination of both sides gives
\begin{equation}\label{rate-1}
8\pi\bar\varepsilon_k\partial_{\xi}\log \hat h_k(0)+\frac{16\pi}{\hat h_k^2(0)}\bar\varepsilon_k^3\partial_{\xi}\Delta\hat  h_k(0)(\log \bar\varepsilon_k^{-2\sigma}+O(1))=O(\bar\varepsilon_k^{4-\delta}).
\end{equation}
Therefore by choosing $\sigma=\sigma_1$ and $\sigma_2$ with both of them close to $\frac 12$ (without affecting the error term $O(\bar\varepsilon_k^{4-\delta})$, we have
\[|(\sigma_1-\sigma_2)\partial_{\xi}\Delta \hat h_k(0)|=O(\bar\varepsilon_k^{1-\delta})\]
for $\delta>0$ arbitrarily small.
Since $\xi$ is any vector on $\mathbb S^1$, we have
\begin{equation}\label{rate-2}
\nabla\Delta \hat h_k(0)=O(\bar\varepsilon_k^{1-\delta}),
\end{equation}
which immediately implies that 
\begin{equation}\label{rate-3}\nabla \hat h_k(0)=O(\bar\varepsilon_k^{3-\delta}),
\end{equation}
by the Pohozaev identity \eqref{rate-1}. 
Tracking the definition of $\hat h_k$ back to (\ref{def-hat-h}) and (\ref{new-h-1}), we see that (\ref{rate-3}) is equivalent to 
\[\nabla \log \bar h_k(q_k)+\nabla \phi_R^k(q_k)(1+\frac{4}{\tau \bar h_k(0)}\bar\varepsilon_k^2)=O(\bar\varepsilon_k^{3-\delta}).\]
Obviously Theorem \ref{reg-int} follows from direct evaluation. Remark that even if one integrates on a disk slightly different from $B_{\tau}$, but because $|q_k|=O(\bar\varepsilon_k^2)$, the difference on these two disks has an area of order $O(\bar\varepsilon_k^2)$ and on this difference $e^{u_k}=O(\bar\varepsilon_k^2)$. So the integral on the difference of two circles is $O(\bar\varepsilon_k^4)$. 
\end{proof}

\begin{rem}As an additional piece of information, after using the new rates in (\ref{new-rel-1}), (\ref{new-rel-2}) and (\ref{new-rel-3}), we can obtain more vanishing estimates on higher derivatives from (\ref{important-2}) if 
$h_k$ is assumed to have more regularity.
\end{rem}
\begin{rem}\label{new-inter-e}
With the improved vanishing rates of $\nabla \bar h_k(q_k)$ and $\nabla \Delta \bar h_k(q_k)$ at hand, it is useful to state here the following estimates of independent interest about the blow up at regular points: 
\begin{equation*}\label{grad-v2-new-0}
			|\bar v_k(y)-(\bar U_k+c_{0,k}+c_{2,k})(y)|\le C\bar \varepsilon_k^{2+\epsilon_0}(1+|y|)^{2\epsilon_0}\quad \mbox{{\rm in}}\quad \Omega_k,
		\end{equation*}
		and
		\begin{equation*}\label{grad-v2-new}
			|\nabla (\bar v_k-\bar U_k-c_{0,k}-c_{2,k})(y)|\le C\tau\bar \varepsilon_k^3,\quad y\in \Omega_k, \quad |y|\sim \bar\varepsilon_k^{-1}.
		\end{equation*}

\end{rem}

\bigskip
	
	\section{Proof of Theorem \ref{main-theorem-2}}
	
	In this section we prove Theorem \ref{main-theorem-2}.
	By contradiction let $\nu_1^k$ and $\nu_2^k$ be two distinct blowup sequences of (\ref{m-equ})
	corresponding to the same $\rho_k$ and the same set of blowup points. The case we consider in this context is $\alpha_M>0$.
	To simplify the exposition, we assume that there exist just three blowup points: $q$, $p_1$, $p_2$, where $q$ is a regular blowup point,$p_1$ is a positive singular source with $\alpha\equiv \alpha_M>0$, $p_2$ is a negative singular source with $\beta\in (-1,0)$.

	Let
	$${\rm w}_i^k=\nu_i^k+4\pi\alpha G(x,p_1)+4\pi \beta G(x,p_2),$$
	then we have
	\begin{equation}\label{uik-tem}
		\Delta_g {\rm w}_i^k+\rho_k\left(\frac{He^{{\rm w}_i^k}}{\int_MHe^{{\rm w}_i^k}}-1\right)=0,
	\end{equation}
	where
	$$H(x)=h(x) e^{-4\pi\alpha G(x,p_1)-4\pi\beta G(x,p_2)}.$$
	Since by adding any constant to ${\rm w}_i^k$ in (\ref{uik-tem}) we still come up with a solution, we fix this free constant to make
	${\rm w}_i^k$ satisfy
	\begin{equation*}
		\int_M He^{{\rm w}_i^k}dV_g=1,\quad i=1,2.
	\end{equation*}
	The basic assumption is that
	$$\sigma_k:=\|{\rm w}_1^k-{\rm w}_2^k\|_{L^{\infty}(M)}>0.$$
	From previous works such as \cite{bart-4}, \cite{wu-zhang-ccm} we have that $\rho_1^k=\rho_2^k=\rho_k$
	implies that $\sigma_k=O(e^{-\gamma\lambda_1^k})$ for some $\gamma>0$ where we set $\lambda_1^k={\rm w}_1^k(p_1)$ and $\lambda_2^k={\rm w}_2^k(p_1)$.
	Denoting by $\hat \lambda_i^k={\rm w}_i^k(p_2)$ and $\bar \lambda^k_i=\max\limits_{B(q,\tau)}{\rm w}_i^k$, it is well known that
	\begin{equation}\label{same-ord}
		\hat \lambda^k_i-\lambda^k_i=O(1),\quad \bar \lambda^k_i-\lambda^k_i=O(1),\quad i=1,2.
	\end{equation}
	
	We work in local isothermal coordinates around $p_1,p_2,q$, the metric locally taking the form,
	$$ds^2=e^{\phi}((dx_1)^2+(dx_2)^2),$$ where $\phi$ satisfies,
	\begin{equation}\label{def-phi}
		\Delta \phi+2Ke^{\phi}=0 \quad{\rm in }\quad B_{\tau}, \quad \phi(0)=|\nabla \phi(0)|=0.
	\end{equation}
	In particular we choose local coordinates so that $p_i$, $q$ are locally the origin of a ball $x\in B_\tau=B_\tau(0)$.
	Remark that in our setting we have,
	\begin{equation}\label{rho-3}
		\begin{cases}
			\rho^*=24\pi+8\pi\alpha+8\pi\beta,\quad N^*=4\pi(\alpha+\beta),\\
			\\
			G^*_1(x)=8\pi(1+\alpha)R(x,p_1)+8\pi(1+\beta)G(x,p_2)+8\pi G(x,q). 
		\end{cases}
	\end{equation}
	
	\begin{rem}
		A word of warning about notations. $B_\tau=B_\tau(0)$ or either $B_\tau(\bar q)$ will denote the ball centered at the origin
		of some local isothermal coordinates $x\in B_\tau$, while $B(p,\tau)\subset M$
		denotes a geodesic ball. Also, if $B_\tau$ is the a ball in local isothermal coordinates centered at some point $0=x(p), p\in M$,
		we will denote by $\Omega(p,\tau)\subset M$
		the pre-image of $B_\tau$. Of course we can always choose $\tau$ small enough to guarantee that
		$\Omega(p_i,\tau)$, $i=1,2$ and $\Omega(q,\tau)$ are simply connected and at
		positive distance one from each other. After scaling,
		$x=\varepsilon y$ for some $\varepsilon>0$, we will denote
		$B_{\tau/\varepsilon}=B_{\tau/\varepsilon}(0)$.\\
		To avoid technicalities, we will some time use the same symbols,
		say $w_k$, $\xi_k$, $\tilde h_k$,  to denote functions when expressed in different local coordinates systems. Of course, it will be
		clear time to time by the context what the given symbol means.
	\end{rem}
	\noindent
	Thus, working in these local coordinates centered at $p_1$, $p_2$, $q$ respectively, we have the local variables $x\in B_\tau$ and we define $f_k$ to be any solution of,
	\begin{equation}\label{def-fk}
		\Delta f_k=e^{\phi}\rho_k\quad {{\rm in}}\quad B_{\tau}, \quad f_k(0)=0,
	\end{equation}
	so that the local equation around 
	$p_1$, $p_2$, $q$, for 
	\begin{equation}\label{uk-def1}
		u_i^k= {\rm w}_i^k-f_k \quad { {\rm in} }\quad B_\tau
	\end{equation}
	takes the form,
	\begin{equation}
		\label{around-p1-h}
		\begin{cases}
			\Delta u_i^k+|x|^{2\alpha}\tilde h_k e^{ u_i^k}=0\quad {{\rm in}}\quad B_{\tau},\quad i=1,2,\qquad 0=x(p_1),\\
			\\
			\tilde h_k=\rho_k h(x)e^{-4\pi \alpha R(x,p_1)-4\pi \beta G(x,p_2)+\phi+f_k},
		\end{cases}
	\end{equation}
	\begin{equation}
		\label{around-p2-h}
		\begin{cases}
			\Delta u_i^k+|x|^{2\beta}\hat h_k e^{ u_i^k}=0\quad {{\rm in}}\quad B_{\tau},\quad i=1,2,\qquad 0=x(p_2)\\
			\\
			\hat h_k=\rho_k h(x)e^{-4\pi \alpha G(x,p_1)-4\pi \beta R(x,p_2)+\phi+f_k},
		\end{cases}
	\end{equation}
	\begin{equation}
		\label{around-q}
		\begin{cases}
			\Delta u_i^k+\tilde h_k e^{ u_i^k}=0\quad {{\rm in}}\quad B_{\tau},\quad i=1,2,\qquad 0=x(q),\\
			\\
			\bar h_k=\rho_k h(x)e^{-4\pi \alpha G(x,p_1)-4\pi \beta G(x,p_2)+\phi+f_k}.
		\end{cases}
	\end{equation}
	
	Obviously $\tilde h_k,\hat h_k, \bar h_k$ satisfy \eqref{assum-H}. Also from this setting we see that (\ref{small-rho-k}) holds. Indeed from 
	$\rho_k=\rho_k\int_MHe^{w_i^k}$ we evaluate the integrals around blowup points and the integration outside the bubbling disks is of the order $O(e^{-\lambda_1^k})$. 
	
	\begin{rem}\label{center-q}
		Concerning \eqref{around-q}, for technical reasons, it will be useful to define $\bar q_i^k\to 0$ to be the maximum points of $u_i^k$ in $B_\tau$ and work 
		in the local coordinate system centered at $\bar q_i^k$. In particular $q_i^k\in M$ will denote the pre-images of these points via the local isothermal map and we will work sometime, possibly taking a smaller $\tau$, with \eqref{around-q} where $0=x(q^k_i)$. 
	\end{rem}
	By using also \eqref{def-phi}, \eqref{def-fk}, we see that the non-degeneracy assumption $L(\mathbf{p})\neq 0$ takes the form
	$$\Delta_g \log h(p_1)+\rho_*-N^*-2K(p_1)\neq 0, $$
	which, for any $k$ large enough, is equivalent to (see (\ref{rho-3}))
	$$\Delta \log \tilde h_k(0)\neq 0.$$ Clearly around $p_1$ we have,
	$$u_i^k=\log \frac{e^{\lambda_1^k}}{(1+\frac{\tilde h_k(0)}{8(1+\alpha)^2}e^{\lambda_1^k}|x|^{2\alpha+2})^2}+\phi_i^k+\mbox{h.o.t.}$$
	where h.o.t stands for higher order terms and $\phi_i^k$ is the harmonic function that encodes the oscillation of $u_i^k$ on $\partial B_\tau$.
	At this point let us define,
	\begin{equation*}
		\xi_k:= ({\rm w}_1^k-{\rm w}_2^k)/\sigma_k.
	\end{equation*}
	Remark that, because of \eqref{uk-def1}, locally we have
	\begin{equation}\label{ueqw}
		u_1^k-u_2^k\equiv {\rm w}_1^k-{\rm w}_2^k
	\end{equation}
	and in particular in local coordinates $\xi_k\equiv(u_1^k-u_2^k)/\sigma_k$. Then it is well known (see \cite{bart-4,bart-4-2,wu-zhang-ccm}) that, after a suitable scaling in local coordinates as in section \ref{difference} below, the limit of $\xi_k$ takes the form,
	\begin{equation}\label{linim-1}
		b_0\frac{1-A|y|^{2\alpha+2}}{1+A|y|^{2\alpha+2}}~\mbox{\rm near } p_1,\quad  A=\lim_{k\to \infty}\frac{\tilde h_k(0)}{8(1+\alpha)^2};
	\end{equation}
	\begin{equation}\label{linim-2}
		b_0\frac{1-B|y|^{2\beta+2}}{1+B|y|^{2\beta+2}}~ \mbox{\rm near } p_2,\quad B=\lim_{k\to \infty}\frac{\hat h_k(0)}{8(1+\beta)^2};
	\end{equation}
	\begin{equation}\label{linim-3}
		b_0\frac{1-C|y|^2}{1+C|y|^2}+b_1\frac{y_1}{1+C|y|^2}+b_2\frac{y_2}{1+C|y|^2}~ \mbox{\rm near } q,  \quad C=\lim_{k\to \infty}\frac{\bar h_k(0)}8.
	\end{equation}
	
	\begin{rem}\label{xi:b_0}
		The reason why we find exactly the same $b_0$ in
		the three formulas is that $\xi_k$ converges to a constant far away from blow up points,
		and it is by now well known (see \cite{bart-4,bart-4-2,wu-zhang-ccm})
		that this constant is exactly $-b_0$.
	\end{rem}

	\subsection{ First estimates about $\sigma_k$. }
	
	As mentioned above we set,
	$$\lambda_i^k={\rm w}_i^k(p_1)\equiv u_i^k(p_1)+f_k,\;i=1,2\quad \mbox{ and  }\quad \varepsilon^2_k=e^{-\frac{\lambda_1^k}{1+\alpha}},$$
	and recall from \eqref{initial-small} that $\lambda_1^k-\lambda_2^k=O(\varepsilon_k^{\gamma})$ for some $\gamma>0$.
	We derive first a rough estimate for ${\rm w}_1^k$ and ${\rm w}_2^k$ far away from blowup points:
	\begin{lem}\label{w12-away} For $x\in M\setminus (\Omega(p_1,\tau)\cup \Omega(p_2,\tau)\cup \Omega(q,\tau))$,
		\begin{equation}
			\label{uk1}
			\begin{cases}
				{\rm w}_1^k(x)=\overline{{\rm w}_1^k}+8\pi(1+\alpha)G(x,p_1)+8\pi(1+\beta)G(x,p_2)\\
				\quad\quad\quad\quad+8\pi G(x,q_1^k)+O(\varepsilon_k^2),\\
				\\
				{\rm w}_2^k(x)=\overline{{\rm w}_2^k}+8\pi(1+\alpha)G(x,p_1)+8\pi(1+\beta)G(x,p_2)\\
				\quad\quad\quad\quad+8\pi G(x,q_2^k)  +O(\varepsilon_k^2),  
			\end{cases}
		\end{equation}
		where $q_i^k$ is defined in Remark \ref{center-q}.
	\end{lem}
	
	\begin{proof}[Proof of Lemma \ref{w12-away}.]
	From the Green representation formula (recall that $\int_M G(x,\eta)=0$) we have,
	\begin{align*}
		&{\rm w}_1^k(x)-\overline{{\rm w}_1^k}=\int_M G(x,\eta)\rho_kH e^{{\rm w}_1^k}\\
		&=\sum_{i=1}^2\bigg (G(x,p_i)\int_{\Omega({p_i},\tau/2)}\rho_kH e^{{\rm w}_1^k}+
		\int_{\Omega({p_i},\tau/2)}(G(x,\eta)-G(x,p_i))\rho_kH e^{{\rm w}_1^k}\bigg )\\
		&\quad+G(x,q_1^k)\int_{\Omega(q_1^k,\tau/2)}\rho_kH e^{{\rm w}_1^k}
		+\int_{\Omega(q_1^k,\tau/2)}(G(x,\eta)-G(x,q_1^k))\rho_kH e^{{\rm w}_1^k}\\
		&\quad+\int_{M\setminus \{\Omega({p_1},\tau/2)\cup \Omega({p_2},\tau/2)\cup
			\Omega(q_1^k,\tau/2)\}}G(x,\eta)\rho_kH e^{{\rm w}_1^k}.
	\end{align*}
	It immediately follows from the estimates in section \ref{preliminary} that the last
	term is $O(e^{-\lambda_1^k})=o(\varepsilon_k^2)$, as far as $\alpha_M>0$.
	By using Theorems \ref{int-pos},  \ref{integral-neg} and  \ref{reg-int}, we have
	$$
	\int_{\Omega({p_1},\tau/2)}\rho_kH e^{{\rm w}_1^k}=\int_{B_{\tau/2}} |x|^{2\alpha}\tilde h_k e^{u_1^k}=8\pi(1+\alpha)+O(\varepsilon_k^2),
	$$
	$$
	\int_{\Omega({p_2},\tau/2)}\rho_kH e^{{\rm w}_1^k}=\int_{B_{\tau/2}} |x|^{2\beta}\hat h_k e^{u_1^k}=8\pi(1+\beta)+O(\varepsilon_k^2),
	$$
	and
	$$
	\int_{\Omega(q_1^k,\tau/2)}\rho_kH e^{{\rm w}_1^k}=\int_{B_{\tau/2}} \bar h_k e^{u_1^k}=8\pi+O(\lambda_1^ke^{-\lambda_1^k})=8\pi +O(\varepsilon_k^2),
	$$
	where in the last integral we also used \eqref{same-ord} and \eqref{p_kj-location}, that is, in local coordinates
	where $\bar q_1^k=x(q_1^k)$, $\bar q:=0=x(q)$, we have
	$\bar q_1^k-\bar q=\bar q_1^k=O(e^{-\bar \lambda^k_1}\bar \lambda^k_1)$. We sketch the argument behind the
	estimates of the other terms.
	When evaluating
	$\int_{\Omega({p_i},\tau/2)} (G(x,\eta)-G(x,p_i))\rho_kH e^{{\rm w}_1^k}$, working in local coordinates $y(\eta)$,
	we use the Taylor expansion of $G(x,\eta_y)$ with respect to $y$ around $0=y(p_i)$ and the
	expansion of $u^k_1$ in Propositions \ref{nov7p1} and \ref{thm-3}. Then it is easy to see that the
	$O(\varepsilon_k)$ term cancels out due to the integration on the symmetric set $B_{\tau/2}$
	and the remaining terms are of order at least
	$O(\varepsilon_k^2)$. In the integration of $\int_{\Omega(q_1^k,\tau/2)}(G(x,\eta)-G(x,q_1^k))\rho_kHe^{{\rm w}_1^k}$,
	we adopt the same argument with minor changes, this time based on Theorem 1.2 in \cite{zhang1},
	which yields as well to an error of order $O(\varepsilon_k^2)$.
	This fact concludes the ${\rm w}_1$ part of the proof of \eqref{uk1}, the part for ${\rm w}_2$ is worked out in the same way.\\
	In particular, by using once more Theorem 1.2 in \cite{zhang1},  we deduce that
	$${\rm w}_1^k(x)-\overline{{\rm w}_1^k}-({\rm w}_2^k(x)-\overline{{\rm w}_2^k})=O(\varepsilon_k^2),~ x\in
	M\setminus \{\Omega({p_1},\tau)\cup \Omega({p_2},\tau)\cup
	\Omega(q,\tau)\},$$
	and, say for $x_1,x_2\in \partial \Omega(q,\tau)$,
	$${\rm w}_1^k(x_1)-{\rm w}_1^k(x_2)={\rm w}_2^k(x_1)-{\rm w}_2^k(x_2)+O(\varepsilon_k^2).$$
	Lemma \ref{w12-away} is established. \end{proof}
	
	\begin{rem}
		%\label{rem-osc}
		Lemma \ref{w12-away} shows in particular that, around each blow up point,
		the difference of the harmonic functions encoding the oscillation on the boundary
		is of order $O(\varepsilon_k^2)$.
	\end{rem}
	\noindent
	Here and in the rest of the proof we denote by,
	$$c_k(x)=\begin{cases}
		\dfrac{e^{{\rm w}_1^k(x)}-e^{{\rm w}_2^k(x)}}{{\rm w}_1^k(x)-{\rm w}_2^k(x)},\quad &\mbox{if} \quad {\rm w}_1^k(x)\neq {\rm w}_2^k(x),\\
		\\
		e^{{\rm w}_1^k(x)}, \quad &\mbox{if }\quad {\rm w}_1^k(x)={\rm w}_2^k(x),
	\end{cases}
	$$
	and use the local coordinates sequences as defined in \eqref{uk-def1},
	$$\tilde c_k(x)=\begin{cases}
		\dfrac{e^{u_1^k(x)}-e^{u_2^k(x)}}{u_1^k(x)-u_2^k(x)},\quad &\mbox{if} \quad u_1^k(x)\neq u_2^k(x),\\
		\\
		e^{u_1^k(x)}, \quad &\mbox{if }\quad u_1^k(x)=u_2^k(x),
	\end{cases}
	$$
	where (see \eqref{def-fk}),
	\begin{equation}\label{ceqtc}
		c_k(x)=e^{f_k(x)}\tilde c_k(x).
	\end{equation}
	We first prove that $b_0=0$. By contradiction assume that $b_0\neq 0$, then we claim that,
	\begin{lem}\label{small-sigma}
		$\sigma_k=O(\varepsilon_k^{2+\epsilon_0})$ for some $\epsilon_0>0$.
	\end{lem}
	\begin{proof}[Proof of Lemma \ref{small-sigma}.]
	We use in a crucial way the fact that, as far as $b_0\neq 0$, then
	\begin{equation}\label{sigma:eq}
		\frac{b_0}{2}({\lambda_{1}^k-\lambda_{2}^k})\leq \sigma_k\leq ({\lambda_{1}^k-\lambda_{2}^k})2b_0,
	\end{equation}
	for any $k$ large enough, see \eqref{radial-xi} below.
	This being said, the proof  relies on the assumption that
	$\rho_1^k=\rho_2^k=\rho_k$, which we write as follows,
	\begin{align*}
		\rho_i^k=\int_M\rho_i^kHe^{{\rm w}_i^k}=&\int_{\Omega({p_1},\tau)}\rho_i^kHe^{{\rm w}_i^k}+\int_{\Omega({p_2},\tau)}\rho_i^kHe^{{\rm w}_i^k}+
		\int_{\Omega({q},\tau)}\rho_i^kHe^{{\rm w}_i^k}\\
		&+\int_{M\setminus \{\Omega({p_1},\tau)\cup \Omega({p_2},\tau)\cup \Omega({q},\tau)\}}\rho_i^kHe^{{\rm w}_i^k}
	\end{align*}
	for $i=1,2$. First of all we have,
	\begin{align*}\int_{M\setminus \{\Omega({p_1},\tau)\cup \Omega({p_2},\tau)\cup \Omega({q},\tau)\}}(\rho_1^kHe^{{\rm w}_1^k}-\rho_2^kHe^{{\rm w}_2^k})/\sigma_k\\
		=\int_{M\setminus \{\Omega({p_1},\tau)\cup \Omega({p_2},\tau)\cup \Omega({q},\tau)\}}\rho_k Hc_k(x)\xi_k.
	\end{align*}
	Since, far away from blowup points,  $\xi_k\to -b_0$ and $c_k(x)\sim e^{-\lambda_1^k}$, we have that,
	$$\int_{M\setminus \{\Omega({p_1},\tau)\cup \Omega({p_2},\tau)\cup \Omega({q},\tau)\}}(\rho_1^kHe^{{\rm w}_1^k}-\rho_2^kHe^{{\rm w}_2^k})=
	O(e^{-\lambda_1^k})
	\sigma_k.$$
	Now we focus on the comparison of the integrals on $\Omega({p_1},\tau)$, $\Omega({p_2},\tau)$ and
	$\Omega({q},\tau)$.
	We first use Theorem \ref{int-pos} to evaluate
	$$\int_{\Omega({p_1},\tau)}\rho_kH(e^{{\rm w}_1^k}-e^{{\rm w}_2^k}).$$
	Here $\tilde h_k$ is the weight function defined in  \eqref{around-p1-h}
	relative to  $u_1^k={\rm w}_1^k-f_k$ and $u_2^k={\rm w}_2^k-f_k$ in $\Omega({p_1},\tau)$.
	After applying Theorem \ref{int-pos} to ${\rm w}_1^k$ and ${\rm w}_2^k$, respectively and looking at their difference, we have 
	\begin{align*}
		\int_{\Omega({p_1},\tau)}\rho_kH(e^{{\rm w}_1^k}-e^{{\rm w}_2^k})
		=d_{1,k}\Delta \log \tilde h_k(0)
		e^{-\lambda_1^k/(1+\alpha)}%+o(\varepsilon_{k}^2)
		(\lambda_1^k-\lambda_2^k)+O(\varepsilon_{k}^{4+\epsilon_0}), 
	\end{align*}
	where $d_{1,k}\to d_1<0$. As remarked above $L(\mathbf q)\neq 0$ is equivalent to
	$\Delta \log \tilde h_k(0)\neq 0$, implying that in fact
	$d_{1,k}\Delta \log \tilde h_k(0)\varepsilon_k^2$ is the leading term.
	By using Theorem \ref{integral-neg} for $\int_{\Omega({p_2},\tau)}\rho_kH(e^{{\rm w}_1^k}-e^{{\rm w}_2^k})$ we have,
	$$\int_{\Omega({p_2},\tau)}\rho_kH( e^{{\rm w}_1^k}-e^{{\rm w}_2^k})=o(\varepsilon_k^2)(\lambda_1^k-\lambda_2^k)+
	O(\varepsilon_k^{4+\epsilon_0}).$$
	At last, for the integration on $\Omega({q},\tau)$, we use Theorem \ref{reg-int} to deduce that,
	$$\int_{\Omega({q},\tau)}\rho_kH(e^{{\rm w}_1^k}-e^{{\rm w}_2^k})=
	o(\bar \varepsilon_k^2)(\lambda_1^k-\lambda_2^k)+O(\bar \varepsilon_k^{4-\sigma}),$$
	for some small $\sigma>0$.
	In this case we note that, since $\alpha>0$, then by possibly choosing a
	smaller $\epsilon_0>0$, we can guarantee that
	$O(\bar \varepsilon_k^{4-\sigma})=O(\varepsilon_k^{4+\epsilon_0})$. Therefore we finally deduce that,
	$$
	0=(d_{1,k}\Delta \log h_k(0)e^{-\lambda_1^k/(1+\alpha)}+o(\varepsilon_{k}^2))
	(\lambda_1^k-\lambda_2^k)+O(\varepsilon_{k}^{4+\epsilon_0}),
	$$
	which immediately implies $\sigma_k=O(\varepsilon_k^{2+\epsilon_0})$. Lemma \ref{small-sigma} is proved. 
\end{proof}
	
	\medskip
	
	The following subsection is devoted to the completion of the proof of $b_0=0$, which is rather long and non trivial.

	\subsection{Proof of $b_0=0$: improved estimates about $\xi_k$.}\label{difference}
	
	The smallness of $\sigma_k$ is crucial for us to obtain a better expansion of $\xi_k$. 
	Recall that
	\begin{equation*}\xi_k(x)=\dfrac{{\rm w}_1^k(x)-{\rm w}_2^k(x)}{\parallel {\rm w}_1^k-{\rm w}_2^k\parallel_{L^{\infty}(M)}},
	\end{equation*}
	and observe that, in view of Lemma \ref{small-sigma}, $c_k(x)$ can be written as follows,
	\begin{equation}
		\label{c}
		\begin{aligned}
			c_k(x)&=\int_0^1e^{t{\rm w}_1^k+(1-t){\rm w}_2^k}dt=e^{{\rm w}_2^k}(1+\frac 12({\rm w}_1^k-{\rm w}_2^k)+O(\sigma_k^2)),\\
			&=e^{{\rm w}_2^k}(1+O(\varepsilon_k^{2+\epsilon_0})), 
		\end{aligned}
	\end{equation}
	which in view of \eqref{ceqtc} obviously holds for $\tilde c_k$ as well.
	Therefore, in local coordinates around $p_1$, so that $0=x(p_1)$, $x\in B_\tau=B_\tau(0)$, $\xi_k$ satisfies,
	\begin{equation*}
		\Delta\xi_k(x)+|x|^{2\alpha}\tilde h_k(x) \tilde c_k(x)\xi_k(x)=0\quad {{\rm in}}\quad B_{\tau},
	\end{equation*}
	where $\tilde {h}_k$ is defined in \eqref{around-p1-h}. In particular recall that
	the assumption $L(\mathbf{q})\neq 0$ is equivalent to $\Delta \log \tilde h_k(0)\neq 0$ for $k$ large
	enough and that we set $\lambda_{i}^k={\rm w}_i^k(p_1) \quad i=1,2$. Next,
	let $\psi_{\xi}^k$ be the harmonic
	function which encodes the oscillation of  $\xi_k$ on
	$\partial B_\tau$,
	\begin{align*}
		\Delta \psi_{\xi}^k=0~ \mbox{in}~ B_\tau,\quad 
		\psi_{\xi}^k(x)=\xi_k(x)-\frac{1}{2\pi\tau}\int_{\partial B_{\tau}}\xi_k(s)ds~\mathrm{for}~ x\in \partial B_\tau,
	\end{align*}
	and $\phi_i^k$,  $i=1,2$ the harmonic functions which encode the oscillations of  $u_i^k$, $i=1,2$ on
	$\partial B_\tau$. At this stage, since $\xi_k$ converges to $-b_0$ far away from blow up points,
	we can just claim that the oscillation of $\xi_k$ on $\partial B_\tau$ is of order $o(1)$.
	We shall derive first a refined estimate about $\xi_k$ around $p_1$.\\
	First of all observe that, since locally
	around $p_1$ the rescaled function $\xi_k(\varepsilon_ky)$ (recall $0=x(p_1)$) converges on compact subsets of $\mathbb{R}^2$
	to a solution of the following linearized equation (\cite{bart-4,bart-4-2,wu-zhang-ccm}),
	$$\Delta \xi+|y|^{2\alpha}\tilde h(0)e^{U}\xi=0,\quad \tilde h(0)=\lim_{k\to \infty}\tilde h_k(0),$$
	where (see \eqref{st-bub}) $U$ is the limit of the standard bubble $U_2^k$, then (see Lemma \ref{lem1}) we have that,
	\begin{equation}\label{lim-1.0}
		\xi(y)=b_0\frac{1-c|y|^{2\alpha+2}}{1+c|y|^{2\alpha+2}},\quad c=\frac{\tilde h(0)}{8(1+\alpha)^2}.
	\end{equation}
	In particular the local convergence of $\xi_k(\varepsilon_ky)$ to $\xi(y)$ readily implies that,
	\begin{equation}\label{radial-xi}
		\frac{\lambda_{1}^k-\lambda_{2}^k}{\|{\rm w}_1^k-{\rm w}_2^k\|_{\infty}}\to b_0.
	\end{equation}
	At this point we can use (\ref{radial-xi}) and the expansions of $u_1^k$ and $u_2^k$ to identify the
	leading order of the radial part of $\xi_k(\varepsilon_ky)$. Let us recall that by definition
	$v_i^k$ is the scaling of $u_i^k$:
	\begin{equation*}
		v_i^k(y)=u_{i}^k(\varepsilon_{k}y)-\phi_{i}^k(\varepsilon_{k}y)-\lambda_{i}^k,\quad i=1,2,
	\end{equation*}
	where $\varepsilon_{k}=e^{-\frac{\lambda_{1}^k}{2(1+\alpha)}}$.
	Thus, by using $U_1^k$ and $U_2^k$ to denote the leading terms of $v_1^k$ and $v_2^k$, respectively, we have,
	$$U_1^k-U_2^k=\lambda_1^k-\lambda_2^k+2\log \frac{1+a_ke^{\lambda_2^k-\lambda_1^k}|y|^{2\alpha+2}}{1+a_k|y|^{2\alpha+2}},\quad a_k=\tilde h_k(0)/8(1+\alpha)^2.$$
	Therefore, in view of  \eqref{initial-small}, we have,
	\begin{align*}
		U_1^k(y)-U_2^k(y)
		=~&\lambda_{1}^k-\lambda_{2}^k+2\log (1+\frac{a_k|y|^{2\alpha+2}(e^{\lambda_{2}^k-\lambda_{1}^k}-1)}{1+a_k|y|^{2\alpha+2}})  \\
		=~&\lambda_{1}^k-\lambda_{2}^k+2\log (1+\frac{a_k|y|^{2\alpha+2}(\lambda_{2}^k-\lambda_{1}^k+O(|\lambda_{2}^k-\lambda_{1}^k|^2)}{1+a_k|y|^{2\alpha+2}})\\
		=~&\lambda_{1}^k-\lambda_{2}^k+2\frac{a_k|y|^{2\alpha+2}(\lambda_{2}^k-\lambda_{1}^k)}{1+a_k|y|^{2\alpha+2}}
		+\frac{O(\lambda_{1}^k-\lambda_{2}^k)^2}{(1+a_k|y|^{2\alpha+2})^2}\\
		=~&(\lambda_{1}^k-\lambda_{2}^k)\frac{1-a_k|y|^{2\alpha+2}}{1+a_k|y|^{2\alpha+2}}+\frac{O(\lambda_{1}^k-\lambda_{2}^k)^2}{(1+a_k|y|^{2\alpha+2})^2}.
	\end{align*}
	
	This fact suggests that the leading order of the radial part
	of $(v_1^k-v_2^k)/\|u_1^k-u_2^k\|_{\infty}$ converges to $\xi$,
	which is a radial function as well. This is why, as a first approximation, we choose the following radial
	function,
	\begin{equation}\label{w0kxik}
		w_{0,\xi}^k=b_0^k\frac{1-a_k|y|^{2\alpha+2}}{1+a_k|y|^{2\alpha+2}},
	\end{equation}
	which is a solution of
	\begin{equation}\label{e-w0kxi}
		\Delta w_{0,\xi}^k +|y|^{2\alpha}\tilde h_k(0)e^{U^k_{2}}w_{0,\xi}^k=0,
	\end{equation}
	where $b_0^k$ is chosen to ensure that $w_{0,\xi}^k(0)=\xi_k(0)$.
	Remark that the equation for $\xi_k-\psi_{\xi}^k$ takes the form,
	\begin{equation}\label{comp-1}
		\Delta (\xi_k-\psi_{\xi}^k)+|x|^{2\alpha}\tilde h_k\tilde c_k(\xi_k-\psi_{\xi}^k)=-|x|^{2\alpha}\tilde h_k\tilde c_k\psi_{\xi}^k.
	\end{equation}
	Our next goal is to obtain an expansion of $(\xi_k-\psi_{\xi}^k)(\varepsilon_ky)$ whose leading
	term is in fact $w_{0,\xi}^k$. This is why we write \eqref{e-w0kxi} as follows,
	\begin{equation}
		\label{e-hat-xi}
		\begin{aligned}
			&\Delta w_{0,\xi}^k+|y|^{2\alpha}\tilde h_k(\varepsilon_ky)\varepsilon_k^{2+2\alpha}
			\tilde c_k(\varepsilon_ky)w_{0,\xi}^k\\
			&=|y|^{2\alpha}\tilde h_k(0)\bigg (\frac{\tilde h_k(\varepsilon_ky)}{\tilde h_k(0)}
			\varepsilon_k^{2+2\alpha}\tilde c_k(\varepsilon_ky)-e^{U_2^k}\bigg )w_{0,\xi}^k(y). 
		\end{aligned}
	\end{equation}
	At this point we consider the equation for
	\begin{equation}\label{tildewk}
		\tilde w_k(y):=\xi_k(\varepsilon_ky)-\psi_{\xi}^k(\varepsilon_ky)-w_{0,\xi}^k(y).
	\end{equation}
	From (\ref{comp-1}) and (\ref{e-hat-xi}) we have
	\begin{equation}
		\label{crucial-1}
		\begin{aligned}
			&\Delta \tilde w_k+|y|^{2\alpha}\tilde h_k(\varepsilon_k y)\varepsilon_k^{2+2\alpha}
			\tilde c_k(\varepsilon_ky)\tilde w_k\\
			&=\tilde h_k(0) |y|^{2\alpha}
			\bigg (e^{U_2^k}-\frac{\tilde h_k(\varepsilon_ky)}{\tilde h_k(0)}\varepsilon_k^{2+2\alpha}
			\tilde c_k(\varepsilon_ky)\bigg)w_{0,\xi}^k(y)\\
			&\quad-|y|^{2\alpha}\tilde h_k(\varepsilon_k y)\varepsilon_k^{2+2\alpha}
			\tilde c_k(\varepsilon_ky)\psi_{\xi}^k(\varepsilon_ky), 
		\end{aligned}
	\end{equation}
	where in particular we also know that
	$$\tilde w_k(0)=0,\quad \tilde w_k=\mbox{constant} \quad \mbox{on}\quad \partial B_{\tau \varepsilon_k^{-1}}.$$
	Remark that the last term of (\ref{crucial-1}) is of order $o(\varepsilon_k)(1+r)^{-3-2\alpha}$ and
	consequently, by Remark \ \ref{stand-rem},  we obtain a first estimate about $\tilde w_k$:
	\begin{equation}\label{first-e}
		|\tilde w_k(y)|\le C(\delta)\varepsilon_k (1+|y|)^{\delta}\quad \mbox{in}\quad B_{\tau \varepsilon_k^{-1}}.
	\end{equation}
	The analysis around $p_2$ is similar, we define (recall that now we have local coordinates where $0=x(p_2)$),
	\begin{equation}\label{hat-wk}
		\hat w_k(y)=\xi_k(\hat \varepsilon_ky)-\hat \psi^k_{\xi}(\tilde \varepsilon_ky)-\hat w^k_{0,\xi}(y), 
	\end{equation}
	where $\hat \varepsilon_k=e^{-\frac{\hat \lambda_k}{2(1+\beta)}}$, $\hat \psi^k_{\xi}$ is the harmonic
	function that encodes the oscillation of $\xi_k$ on $\partial B_\tau$, $\hat w^k_{0,\xi}$ is the same as
	$w^k_{0,\xi}$ just with $\beta$ replacing $\alpha$ and the local limit of
	$\xi_k(\hat \varepsilon_ky)$ is now \eqref{linim-2} which replaces \eqref{lim-1.0}.
	By using the same argument as that adopted above for $\tilde w_k$, we have that,
	\begin{equation}\label{sec-e}
		|\hat w_k(y)|\le C(\delta)\hat \varepsilon_k^{1+\beta}(1+|y|)^{\delta}=C(\delta)e^{-\lambda_1^k/2}(1+|y|)^{\delta},\quad
		\mbox{in}\quad B_{\tau \varepsilon_k^{-1}}.
	\end{equation}
	
	Next we describe the expansion of $\xi_k$ near $q$. In this case \eqref{lim-1.0} is replaced by \eqref{linim-3}. Recall that we denote $\bar q_i^k=x(q_i^k)$, $i=1,2$ the expressions in local isothermal
	coordinates $0=x(q)$ of the maximizers near $q$ of $u_i^k$, $i=1,2$.
	From \eqref{p_kj-location} we have
	$|\bar q_1^k-\bar q_2^k|=O(e^{-\lambda_1^k}\lambda_1^k)$. In $\Omega(q,\tau)$ we denote by
	$\bar \varepsilon_k=e^{-\bar \lambda_2^k/2}$ the scaling factor
	where $\bar \lambda_1^k$ and $\bar \lambda_2^k$ are the maximums near $q$. With an abuse of notations
	we denote by $\bar \psi^k_\xi$
	the harmonic function that, in local isothermal
	coordinates $0=x(q)$,
	encodes the oscillation of $\xi_k$ on $B_\tau$.
	Of course, since we assume $b_0\neq 0$, by analogy with \eqref{radial-xi}
	we still have that $|\bar \lambda_1^k-\bar \lambda_2^k|\sim \sigma_k$.

	Now we work in more detail about the asymptotic expansion of $\xi_k$
	around $q$ and use local coordinates
	such that $0=x(q)$. As usual the kernel functions are the first terms in the approximation of
	$\xi_k(\bar q_2^k+\bar \varepsilon_ky)$,
	$$\bar w_{0,\xi}^k(y)=\bar b_0^k\frac{1-\frac{\bar h_2^k(0)}8|y|^2}
	{1+\frac{\bar h_2^k(0)}8|y|^2}+
	\bar b_1^k\frac{y_1}{1+\frac{\bar h_2^k(0)}8|y|^2}+
	\bar b_2^k\frac{y_2}{1+\frac{\bar h_2^k(0)}8|y|^2}.$$
	Then we have 
	\begin{equation*}
		%\label{4.equ-1st}
		\xi_k(\bar q_2^k)-\bar w_{0,\xi}^k(0)=0\quad\mbox{and}\quad \nabla(\xi_k(\bar q_2^k+\bar\varepsilon_ky)-\bar w_{0,\xi}^k(\bar\varepsilon_ky))|_{y=0}=0.
	\end{equation*}
	Next we set $\bar\psi_\xi^k$ to be the harmonic function which encodes the oscillation of $\xi_k(\bar q_2^k+\bar\varepsilon_ky)-\bar w_{0,\xi}^k(\bar\varepsilon_ky)$ on $\partial B_{\tau/\bar\varepsilon_k}$. Let us define
	\begin{equation}\label{bar-w2k}
		\bar w_k(y)=\xi_k(\bar q_2^k+\bar \varepsilon_ky)-\bar w_{0,\xi}^k(y)-\bar \psi_{\xi}^k(\bar \varepsilon_ky),
	\end{equation}
	and we have
	\begin{equation}\label{eqwbar1}
		\begin{cases}
			\Delta \bar w_k(y)+\bar \varepsilon_k^2(\tilde h_k \tilde c_k)(\bar q_2^k+\bar \varepsilon_ky)\bar w_k(y)
			=\frac{O(\bar \varepsilon_k)}{(1+|y|)^3} ~ \mbox{in}~ B_{\tau /\bar\varepsilon_k},\\
			\bar w_k(0)=0\quad \mbox{and}\quad \bar w_k={\rm constant}~\mbox{on}~\partial B_{\tau/\bar \varepsilon_k}.
		\end{cases}
	\end{equation}
	It is not difficult to check that $\bar w_k(0)=0$ and $\nabla \bar w_k(0)=-\nabla\bar\psi_\xi^k(0)=O(\bar\varepsilon_k)$ by the fact that $\bar \psi_\xi^k$ is a harmonic function. Then using the Green representation formula and standard potential estimates (see for example the proof of Proposition \ref{nov7p1}) we have the following preliminary estimate:
	\begin{equation}\label{pre-xi}
		|\bar w_k(y)|\le C(\delta)\bar \varepsilon_k(1+|y|)^{\delta} \quad {in}\quad B_{\tau\bar \varepsilon_k^{-1}}.
	\end{equation}
	Now we can improve the estimate on the oscillation of $\xi_k$  
	\begin{lem}\label{osci-xi-better}
		\begin{equation}\label{osi-vars}
			\xi_k(x_1)-\xi_k(x_2)=O(\varepsilon_k^{2}+e^{-\lambda_1^k/2})=O(\varepsilon_k^2+\varepsilon_k^{1+\alpha}),
		\end{equation}
		$\forall\, x_1,x_2\in
		M\setminus \{\Omega({p_1},\tau)\cup \Omega({p_2},\tau)\cup \Omega({q},\tau)\}$.
	\end{lem}
	\begin{rem} We will need the full strength of \eqref{osi-vars} later on, for the time being we will use its weaker form
		$\xi_k(x_1)-\xi_k(x_2)=O(\varepsilon_k^{2}+\varepsilon_k^{1+\alpha})=O(\varepsilon_k^{1+\epsilon_0})$, for some
		$\epsilon_0>0$.
	\end{rem}
	\begin{proof}[Proof of Lemma \ref{osci-xi-better}.]
	We shall divide our proof into two steps, in the first step we derive the following estimate 
	\begin{equation}
		\label{4.est-step1}
		\xi_k(x_1)-\xi_k(x_2)=O(\varepsilon_k),
	\end{equation}
	for any $x_1,x_2\in
	M\setminus \left\{\Omega({p_1},\tau)\cup \Omega({p_2},\tau)\cup \Omega({q},\tau)\right\}.$ Then in the second step we shall improve the estimation and prove the desired estimate \eqref{osi-vars}.
	\smallskip
	
	\noindent {\bf Step 1.} By using the Green representation formula for
	$\xi_k$, we have, for $x\in M\setminus \{\Omega({p_1},\tau)\cup \Omega({p_2},\tau)\cup \Omega({q},\tau)\},$
	\begin{align*}
		\xi_k(x)&=\bar \xi_k+\int_M G(x,\eta)\rho_kH c_k(\eta)\xi_k(\eta)d\eta\\
		&=\bar \xi_k+\sum_{i=1}^2\int_{\Omega({p_i},\tau/2)}G(x,\eta)\rho_kH c_k(\eta)\xi_k(\eta)d\eta\\
		&\quad+\int_{\Omega({q},\tau/2)}G(x,\eta)\rho_kH c_k(\eta)\xi_k(\eta)d\eta\\
		&\quad+\int_{M\setminus \{\Omega({p_1},\tau/2)\cup \Omega({p_2},\tau/2)\cup \Omega(q,\tau/2)\}}G(x,\eta)\rho_kH c_k(\eta)\xi_k(\eta)d\eta\\
		&=\bar \xi_k+I_1+I_2+I_3+I_4,
	\end{align*}
	with obvious meaning. Clearly $I_4$ is of order $O(e^{-\lambda_{1}^k})$ and we write,
	\begin{align*} I_1&=G(x,p_1)\int_{\Omega({p_1},\tau/2)}\rho_kHc_k\xi_k+
		\int_{\Omega({p_1},\tau/2)}(G(x,\eta)-G(x,p_1))\rho_kHc_k\xi_k\\
		&=I_{11}+I_{12}.
	\end{align*}
	Concerning $I_{12}$, by the Mean Value Theorem, in local coordinates after scaling we have $$G(x,\eta)-G(x,p_1)=\varepsilon_k(a_0\cdot y)+O(\varepsilon_k^2|y|^2),$$  
	then in view of Lemma \ref{small-sigma}, \eqref{c} and (\ref{first-e}), we see that,
	\begin{align*}
		I_{12}
		=~&\varepsilon_k\int_{B_{\frac{\tau}2 \varepsilon_k^{-1}}}\tilde h_k(0)(a_0 \cdot y)|y|^{2\alpha}e^{U_2^k}
		(\psi_{\xi}^k(\varepsilon_ky)+w_{0,\xi}^k(y)\\
		&+O(\varepsilon_k(1+|y|)^{\delta}))+O(\varepsilon_k^2)\\
		=~&\varepsilon_k\int_{B_{\frac{\tau}2 \varepsilon_k^{-1}}}\tilde h_k(0)(a_0 \cdot y)|y|^{2\alpha}e^{U_2^k}w_{0,\xi}^k(y)+O(\varepsilon_k^2)=O(\varepsilon_k^2),
	\end{align*}where we used $\psi^k_\xi(\varepsilon_ky)=O(\varepsilon_k|y|)$ and the fact that
	the first term on the right vanishes,
	$$
	\varepsilon_k\int_{B_{\frac{\tau}2 \varepsilon_k^{-1}}}\tilde h_k(0)(a_0 \cdot y)|y|^{2\alpha}e^{U^k_2}
	w^k_{0,\xi}(y)=0,
	$$
	because $w^k_{0,\xi}$ is a
	radial function, whence the integrand is a separable function. The estimate about $I_{11}$ is similar,
	\begin{align*}
		I_{11}=~&G(x,p_1)\int_{\Omega({p_1},\tau/2)}\rho_kHc_k(\eta)(\xi_k-\psi_{\xi}^k)+G(x,p_1)\int_{\Omega({p_1},\tau/2)}\rho_kHc_k(\eta)\psi_{\xi}^k\\
		=~&I_{1a}+I_{1b},
	\end{align*}
	where, writing $\psi_{\xi}^k(\varepsilon_ky)=\varepsilon_k a_k \cdot y+O(\varepsilon_k^2|y|^2)$, we have
	$$
	I_{1b}=O(1)\int_{B_{\frac{\tau}2 \varepsilon_k^{-1}}}\rho_k\tilde h_k(0)|y|^{2\alpha}e^{U^k_2} \psi_{\xi}^k(\varepsilon_ky)=O(\varepsilon_k^2),
	$$
	again because the term of order $O({\varepsilon_k})$ is the integral of a
	separable function. On the other side, $I_{1a}=O(\varepsilon_k)$ since we can write,
	\begin{equation}
		\label{i1a}
		\begin{aligned}
			\int_{\Omega({p_1},\tau/2)}\rho_kHc_k(\eta)(\xi_k-\psi_{\xi}^k)
			=&\int_{B_{\frac{\tau}2 \varepsilon_k^{-1}}}|y|^{2\alpha}\tilde h_k(\varepsilon_k y)e^{U^k_2}w_{0,\xi}^kdy+O(\varepsilon_k)\\
			=&-\int_{\partial B_{\frac{\tau}2 \varepsilon_k^{-1}}}\partial_{\nu}w_{0,\xi}^k+O(\varepsilon_k)\\
			=&~O(\varepsilon_k^{2+\epsilon_0})+O(\varepsilon_k)=O(\varepsilon_k),
		\end{aligned}
	\end{equation}
	where in the first equality we used (\ref{first-e}) and in the second equality
	we used a second order Taylor expansion of $\tilde h_k(\varepsilon_k y)$,
	the fact that the integral on $B_{\frac{\tau}2 \varepsilon_k^{-1}}$ of an integrable
	separable function vanishes and at last
	the equation for $w_{0,\xi}^k$ and the fact that $\partial_{\nu}w_{0,\xi}^k=O(r^{-3-2\alpha})$
	on the boundary. {It is crucial to point out 
		$I_{1a}$ is the only term of order $\varepsilon_k$,  and all the other terms are either of order $\varepsilon_k^2$ or $e^{-\lambda_1^k/2}=\varepsilon_k^{1+\alpha}.$} While for the integral $I_2$, by \eqref{sec-e} and a similar argument we deduce that $I_2=O(\varepsilon_k^{-\lambda_1^k/2})$. Concerning $I_3$ we have,
	\begin{align*}
		I_3&=O(\varepsilon_k^{2})+\int_{\Omega(q^k_2,\tau/2)}G(x,\eta)\rho_kHc_k(x)\xi_k\\
		&=\int_{\Omega(q^k_2,\tau/2)}(G(x,\eta)-G(x,q^k_2))\rho_kHc_k\xi_k+G(x,q^k_2)\int_{\Omega(q^k_2,\tau/2)}\rho_kHc_k\xi_k\\
		&\quad+O(\varepsilon_k^{2})\\
		&=I_{31}+I_{32}+O(\varepsilon_k^{2}).
	\end{align*}
	where we also used \eqref{same-ord} and $\bar q_2^k-\bar q=\bar q_2^k=
	O(e^{-\bar \lambda^k_1}\bar \lambda^k_1)$, see \eqref{p_kj-location}. By arguing as for $I_{12}$ above we find that $I_{32}=O(e^{-\lambda_1^k/2})$. Also, by using  (\ref{pre-xi}) instead of (\ref{first-e}),
	the same argument adopted for $I_{11}$ shows that  $I_{31}=O(e^{-\lambda_1^k/2})$. We skip these details to avoid repetitions.
	Therefore, we have shown so far that,
	$$|\xi_k(x_1)-\xi_k(x_2)|=O(\varepsilon_k)+O(e^{-\lambda_1^k/2})=O(\varepsilon_k+\varepsilon^{1+\alpha}_k),$$
	for any $x_1,x_2\in M\setminus \{\Omega({p_1},\tau)\cup \Omega({p_2},\tau)\cup \Omega({q},\tau)\},$ that is exactly the estimation \eqref{4.est-step1}. As a consequence, we see immediately that $\psi^k_{\xi}$ satisfies,
	\begin{equation}\label{impxi:1}
		|\psi^k_{\xi}(\varepsilon_k y)|\le C\varepsilon_k^2(1+|y|)\quad {\rm in} \quad B_{\tau\varepsilon_k^{-1}}.
	\end{equation}
	We will use these facts in a sort of bootstrap argument to improve the estimates about the integration of $I_{1a}$ and the
	oscillation of $\xi_k$ far away from blow up points in the following.
	\medskip
	%\begin{rem}\label{rem1a} We point out that there is just one term which is not of order\\
	%$O(\varepsilon_k^2+e^{-\lambda_1^k/2})=O(\varepsilon_k^{1+\alpha}+\varepsilon_k^2)$, which is the one which we denoted by $I_{1a}$.
	%\end{rem}
	
	\noindent {\bf Step 2.} We invoke the intermediate estimate in Lemma \ref{intermediate-2} about the expansion of $v_2^k$ near $p_1$: 
	\begin{equation}\label{p-rough-2}
		v_2^k=U_2^k+c_{0,k}+c_{1,k}+c_{2,k}+O(\varepsilon_k^{2+\delta})(1+|y|)^{-\delta}.
	\end{equation}
	Although we do not need its full strength right now, by using this approximation and $\sigma_k=O(\varepsilon_k^{2+\epsilon_0})$ (see (\ref{c})),
	for later purposes we state the following expansion in local coordinates around $0=x(p_1)$,
	\begin{equation}
		\label{fullexp:1}
		\begin{aligned}
			&\tilde h_k(\varepsilon_ky)\varepsilon_k^{2+2\alpha}\tilde c_k(\varepsilon_ky)\\
			&=\tilde h_k(0)\exp \left\{v_2^k+\log \frac{\tilde h_k(\varepsilon_ky)}{\tilde h_k(0)}\right\}(1+O(\varepsilon_k^{2+\epsilon_0}))\\
			&=\tilde h_k(0)e^{U_2^k}\left(1+c_{0,k}+\tilde c_{1,k}+\tilde c_{2,k}+\frac {\Delta (\log \tilde h_k)(0)|y|^2\varepsilon_k^2}{4}+\frac{\tilde c_{1,k}^2}{2}\right)\\
			&\quad+E_k,   
		\end{aligned}
	\end{equation}
	where $c_{0,k}$ is a radial function corresponding to the $\varepsilon_k^2$ order term in the expansion of $v_2^k$ (see (\ref{c0k})), $\tilde c_{1,k}$ (defined as in  (\ref{t-c1k})) and $\tilde c_{2,k}$ are given below:
	\begin{align*}
		&\tilde c_{1,k}=\varepsilon_k\nabla\log \tilde h_k(0)\cdot \theta\left(r-\frac{2(1+\alpha)}{\alpha}\frac{r}{1+a_kr^{2\alpha+2}}\right),\quad \tilde c_{2,k}=c_{2,k}+\Theta_2\varepsilon_k^2r^2,
	\end{align*}
	where $a_k=\frac{\tilde h_k(0)}{8(1+\alpha)^2}.$ The error term $E_k$ in \eqref{fullexp:1} satisfies,
	$$|E_k|\le C\varepsilon_k^{2+\epsilon_0}(1+|y|)^{-2+\epsilon_0-2\alpha}.$$

	\noindent At this point, around $p_1$ we define $w_{1,\xi}^k$,
	the next term in the approximation of $\xi_k$, to be a suitable solution of,
	\begin{equation}\label{xi-2nd}
		\Delta w_{1,\xi}^k+|y|^{2\alpha}\tilde h_k(0)e^{U_2^k}w_{1,\xi}^k=-\tilde h_k(0)|y|^{2\alpha}e^{U_2^k}\tilde c_{1,k} w_{0,\xi}^k.
	\end{equation}
	It turns out that, see \eqref{w1-explicit} below for further details,
	\eqref{xi-2nd}  can be explicitly solved in terms of a function which also satisfies
	$w_{1,\xi}^k(0)=0=|\nabla w_{1,\xi}^k(0)|$ and
	\begin{equation*}
		|w_{1,\xi}^k(y)|\le C\varepsilon_k (1+|y|)^{-1}\quad {\rm in} \quad B_{\tau \varepsilon_k^{-1}}.
	\end{equation*}
	In fact we will see that $\left(\mbox{recall}~\theta=\left(x_1,x_2\right)/|x|\right)$,
	$$
	w_{1,\xi}^k(y)=\frac{\tilde h_k(0)}{4\alpha(1+\alpha)}(\nabla \log \tilde h_k(0)\cdot \theta)
	\frac{\varepsilon_k r^{2\alpha+3}}{(1+r^{2+2\alpha })^2},
	$$
	is a solution of \eqref{xi-2nd} which satisfies all those properties.

	%The construction as well as the proof of (\ref{w1osc}) are obtained essentially by the same argument
	%adopted for \eqref{r-c1k}.
	%To avoid repetitions we just sketch the proof of (\ref{w1osc}). Let $g_{1,k}$ denote the radial part of $w_{1,\xi}^k$ which satisfies,
	%$$g_{1,k}''+\frac 1rg_{1,k}'(r)+\tilde h_k(0)r^{2\alpha}e^{U_2^k}g_{1,k}=E_{1,k},\quad 0<r<\tau\varepsilon_k^{-1},$$
	%where
	%$$|E_{1,k}(r)|\le C\varepsilon_kr(1+r)^{-4-2\alpha}.$$
	%Then we define,
	%$$g_{1,k}(r)=f_{1,1}(r)\int_r^{\tau\varepsilon_k^{-1}}\frac{W_1(s)}{W(s)}ds+f_{1,2}(r)\int_0^r\frac{W_2(s)}{W(s)}ds,$$
	%where
	%$$|W(s)|\le Cs^{-1},\quad |W_1(s)|\le C\varepsilon_k(1+s)^{-4-2\alpha},\quad |W_2(s)|\le C\varepsilon_ks^2(1+s)^{-4-2\alpha},$$
	%and $f_{1,1}$ and $f_{1,2}$ are two fundamental solutions as in (\ref{two-fun}) with
	%$f_{1,1}\sim r$ at $0$ and $f_{1,2}\sim r^{-1}$ at $0$.
	%Then by a straightforward estimate about $g_{1,k}$ we see that (\ref{w1osc}) holds.\\
	\noindent
	Next, let us we write the equation for $w_{1,\xi}^k$ in the following form:
	\begin{equation}
		\label{e-w1k}
		\begin{aligned}
			&\Delta w_{1,\xi}^k+|y|^{2\alpha}\tilde h_k(\varepsilon_k y)
			\varepsilon_k^{2+2\alpha}\tilde c_k(\varepsilon_k y)w_{1,\xi}^k \\
			&=\tilde h_k(0)|y|^{2\alpha}\left[\left(\frac{\tilde h_k(\varepsilon_ky)}{\tilde h_k(0)}
			\varepsilon_k^{2+2\alpha}\tilde c_k(\varepsilon_k y)-e^{U_2^k}\right)w_{1,\xi}^k-e^{U_2^k}\tilde c_{1,k}w_{0,\xi}^k\right],
			%\\&\quad-\tilde h_k(0)|y|^{2\alpha}e^{U_2^k}\tilde c_{1,k}w_{0,\xi}^k,
		\end{aligned}
	\end{equation}
	and then define,
	$$z_k(y)=\tilde w_k(y)-w_{1,\xi}^k(y)=\xi_k(\varepsilon_ky)-\psi_{\xi}^k(\varepsilon_ky)-w_{0,\xi}^k(y)-w_{1,\xi}^k(y).$$
	Therefore, in view of \eqref{crucial-1}, \eqref{impxi:1} and \eqref{e-w1k}, the equation for $z_k$ reads:
	\begin{align*}
		&\Delta z_k+|y|^{2\alpha}\tilde h_k(\varepsilon_k y)\varepsilon_k^{2+2\alpha}\tilde c_k(\varepsilon_k y)z_k\\
		&=\tilde h_k(0)|y|^{2\alpha}(e^{U_2^k}- \frac{\tilde h_k(\varepsilon_ky)}{\tilde h_k(0)}
		\varepsilon_k^{2+2\alpha}\tilde c_k(\varepsilon_k y))(w_{0,\xi}^k+w_{1,\xi}^k)\\
		&\quad+\tilde h_k(0) |y|^{2\alpha}e^{U_2^k}\tilde c_{1,k} w_{0,\xi}^k+O(\varepsilon_k^{2+\epsilon_0})(1+|y|)^{-2-2\alpha+\epsilon_0}\\
		&=O\left(\varepsilon_k^2\right)(1+|y|)^{-2-2\alpha}. 
	\end{align*}
	By using the usual argument based on Green's representation formula (see Remark \ref{stand-rem}), we have
	\begin{equation}\label{est-z}
		|z_k(y)|\le C(\delta)\varepsilon_k^2(1+|y|)^{\delta},
	\end{equation}
	implying in particular that the oscillation of $z_k$ on
	$\partial B_{\tau\varepsilon_k^{-1}}$ is of order $O(\varepsilon_k^{2-\delta})$.
	
	\medskip
	\noindent
	In view of \eqref{est-z}, going back through the estimates above for $\xi_k$,
	concerning in particular $I_{1a}$ we find that,
	\begin{align*}
		I_{1a}&=G(x,p_1)\int_{B(p_1,\frac{\tau}2)}\rho_k Hc_k(\eta)(\xi_k-\psi_{\xi}^k)\\
		&=G(x,p_1)\int_{B_{\frac{\tau}{2}\varepsilon_k^{-1}}}\rho_k\tilde h_k(\varepsilon_ky)e^{U^k_2}
		(w_{0,\xi}^k+w_{1,\xi}^k)dy+O(\varepsilon_k^2).
	\end{align*}
	At this point recall that, after a second order Taylor expansion of $\tilde h_k(\varepsilon_ky)$, the
	zero-th order term relative to $w_{0,\xi}^k$ has already been shown above
	to be of order $O(\varepsilon_k^{2+\epsilon_0})$ (see \eqref{i1a}), the one relative to $w_{1,\xi}^k$ vanishes since the
	integrand is separable, the term proportional to $\varepsilon_k w_{0,\xi}^k$ vanishes
	again because the integrand is separable, while the one proportional to
	$\varepsilon_k w_{1,\xi}^k$ is already of order $O(\varepsilon_k^2)$.
	{On the other hand, we have already mentioned in Step 1 that all the terms except for $I_{1a}$ are of order either $\varepsilon_k^2$ or $\varepsilon_k^{1+\alpha}$}. Therefore, we derive that
	\begin{equation*}
		|\xi_k(x_1)-\xi_k(x_2)|=O(\varepsilon_k^2)+O(e^{-\lambda_1^k/2}),
	\end{equation*}
	$\forall x_1,x_2\in M\setminus \{\Omega({p_1},\tau)\cup \Omega({p_2},\tau)\cup \Omega({q},\tau)\}$,
	which is (\ref{osi-vars}). Hence, we finish the proof.
\end{proof}
	
	\begin{rem}\label{rem 4.7} Even with the improved estimate (\ref{osi-vars}),
		the leading order in the estimate of $I_{32}$ is still $O(e^{-\lambda_1^k/2})$.
		This is due to the fact that near $q^k_i$, see \eqref{linim-3},
		$\xi_k$ has non-radial parts after scaling which contribute in the expansion
		of $I_{32}$ at order $O(\bar \varepsilon_k)=O(e^{-\lambda_1^k/2})$.
	\end{rem}

	Before moving to the next stage we remark that (\ref{osi-vars}) further
	improves the harmonic functions that encode the oscillations of $\xi_k$ near blowup points.
	For example, for $\psi_{\xi}^k$ we now have
	\begin{equation}\label{im-ph-xi}
		|\psi_{\xi}^k(\varepsilon_ky)|\le C\varepsilon_k^{2+\epsilon_0}(1+|y|).
	\end{equation}
	We will see that this estimate implies that all the terms in the expansions involving
	$\psi_{\xi}^k$ are in fact negligible.
	
	Next we identify the radial and non-radial parts of $\tilde w_k$, see \eqref{tildewk}.
	In the second step of the estimates about $\tilde w_k$ we use $E_w$
	to denote the right hand side of (\ref{crucial-1}).
	
	By using \eqref{fullexp:1}, the equation of $z_k$ can be written as follows,
	\begin{align*}
		&\Delta z_k+|y|^{2\alpha}\tilde h_k(\varepsilon_k y)\varepsilon_k^{2+2\alpha}\tilde c_k(\varepsilon_k y)z_k\\
		&=-\tilde h_k(0)|y|^{2\alpha}e^{U_2^k}\bigg (w_{0,\xi}^k\big (c_{0,k}+\tilde c_{2,k}+\frac 14\Delta (\log \tilde h_k)(0)|y|^2\varepsilon_k^2+\frac  12(\tilde c_{1,k})^2\big )\nonumber\\
		&\quad+w_{1,\xi}^k\tilde c_{1,k}+O(\varepsilon_k^{2+\epsilon_0})(1+|y|)^{-2-2\alpha+\epsilon_0} \bigg ). \nonumber
	\end{align*}
	By using separable functions to remove separable terms, we can write the equation of the radial part of $z_k$, which we denote $\tilde z_k$,
	\begin{align*}
		&\Delta \tilde z_k+|y|^{2\alpha}\tilde h_k(0)e^{U_2^k}\tilde z_k\\
		&=-\tilde h_k(0)|y|^{2\alpha}e^{U_2^k}\left[w_{0,\xi}^k\left(c_{0,k}+\frac{\Delta (\log \tilde h_k)(0)|y|^2\varepsilon_k^2}{4}+\frac  {(\tilde c_{1,k})^2_r}{2}\right)+(w_{1,\xi}^k\tilde c_{1,k})_r\right],
		%\\&\quad-\tilde h_k(0)|y|^{2\alpha}e^{U_2^k}(w_{1,\xi}^k\tilde c_{1,k})_r,
	\end{align*}
	where we let $(A)_{\theta},(A)_r$ denote the angular part and radial part of $A$. 
	Corresponding to these terms we construct $z_0^k$ to solve
	\begin{equation}\label{eq-z0}
		\begin{cases}\dfrac{d^2}{dr^2}z_0^k+\frac 1r\dfrac{d}{dr}z_0^k+r^{2\alpha}\tilde h_k(0)e^{U_2^k}z_0^k=F_0^k,
			\quad 0<r<\tau \varepsilon_k^{-1},\\
			\\
			z_0^k(0)=\frac{d}{dr}z_0^k(0)=0,
		\end{cases}
	\end{equation} 
	where
	\begin{align*}
		F_0^k(r)
		=&-\tilde h_k(0)|y|^{2\alpha}e^{U_2^k}\bigg (w_{0,\xi}^k\big (c_{0,k}+\frac 14\Delta (\log \tilde h_k)(0)|y|^2\varepsilon_k^2\\
		&+\frac  12\varepsilon_k^2|\nabla \log \tilde h_k(0)|^2(g_k+r)^2\big )
		+(w_{1,\xi}^k\tilde c_{1,k})_r\bigg )\nonumber
	\end{align*}
	and we used the fact that the radial part of $\frac 12(\tilde c_{1,k})^2$ is just
	$$\frac 12\varepsilon_k^2|\nabla \log \tilde h_k(0)|^2(g_k+r)^2.$$
	Again by standard potential estimates it is not difficult to see that
	$$|z_0^k|\le C\varepsilon_k^2(1+r)^{-2\alpha}\log (1+r).$$
	After defining $w_{2,\xi}^k$ to remove the separable terms of the order $O(\varepsilon_k^2)$
	in $\tilde w_k$,
	\begin{equation}
		\label{first:estimate}
		\begin{aligned}
			&\Delta w_{2,\xi}^k+|y|^{2\alpha}h_k(0)e^{U_2^k}w_{2,\xi}^k\\
			&=-|y|^{2\alpha} h_k(0)e^{U_2^k}\bigg (w_{0,\xi}^k\big (\tilde c_{2,k}+\frac 12(\tilde c_{1,k}^2)_{\theta}\big )+(\tilde c_{1,k}w_{1,\xi}^k)_{\theta}\bigg ),
		\end{aligned}
	\end{equation}
	and using (\ref{two-fun1}) we see that 
	\begin{equation}\label{1st-w2k}
		|w_{2,\xi}^k|\le C\varepsilon_k^2(1+|y|)^{-2\alpha} \quad {\rm in} \quad B_{\tau \varepsilon_k^{-1}}.
	\end{equation}
	Also by standard potential estimates as usual we have
	$$
	\left|\tilde w_k-\sum_{i=1}^2w_{i,\xi}^k-z_0^k(y)\right|\le C\varepsilon_k^{2+\epsilon_0}(1+|y|)^{-2\alpha+\epsilon_0}.
	$$
	Next, by using (\ref{c}) for $\tilde c_k(\varepsilon_ky)$ and Lemma \ref{small-sigma} for $\sigma_k$ and
	recalling the argument about $w_{0,\xi}^k$ in the equation \eqref{i1a},
	we evaluate the integral around $p_1$ as follows,
	\begin{align*}
		&\int_{B_{\tau \varepsilon_k^{-1}}}|y|^{2\alpha}\tilde h_k(\varepsilon_ky)\varepsilon_k^{2+2\alpha}
		\tilde c_k(\varepsilon_ky)\xi_k(\varepsilon_k y)dy\\
		&=\tilde h_k(0)\int_{B_{\tau \varepsilon_k^{-1}}}|y|^{2\alpha}\exp\{U_2^k+c_{0,k}+\tilde c_{1,k}+\tilde c_{2,k}+
		\frac 14\Delta (\log \tilde h_k(0))\varepsilon_k^2|y|^2  \\
		&\quad+\frac 12 (\tilde c_{1,k})^2\}
		(w_{0,\xi}^k+w_{1,\xi}^k+z_0^k)+O(\varepsilon_k^{2+\epsilon_0}) \\
		&=\int_{B_{\tau \varepsilon_k^{-1}}}\tilde h_k(0)|y|^{2\alpha}e^{U_2^k}\bigg (\big (1+ c_{0,k}+\frac 12 (\tilde c_{1,k})^2+\frac 14\varepsilon_k^2\Delta (\log \tilde h_k)(0)|y|^2\big )w_{0,\xi}^k \\
		&\quad+(\tilde c_{1,k} w_{1,\xi}^k)_r+z_0^k \bigg )
		+O(\varepsilon_k^{2+\epsilon_0}).  
	\end{align*}
	\noindent
	Remark that, by using the equation for $z_0^k$ in \eqref{eq-z0}, we have,
	\begin{equation}\label{key-c-2}
		\int_{B_{\tau \varepsilon_k^{-1}}}|y|^{2\alpha}
		\tilde h_k(\varepsilon_ky)\varepsilon_k^{2+2\alpha}\tilde c_k(\varepsilon_ky)\xi_k(\varepsilon_ky)dy
		=-\int_{\partial B_{\tau \varepsilon_k^{-1}}}\frac{\partial z_0^k}{\partial \nu}+O(\varepsilon_k^{2+\epsilon_0}).
	\end{equation}
	Let $c_{1,k}$ be defined as in \eqref{nov10e6} and let us set,
	$$\hat c_1^k(x)=c_{1,k}\left(\frac{x}{\varepsilon_k}\right),$$
	that is,
	$$\hat c_1^k(x)=-\frac{2(1+\alpha )}{\alpha}(\nabla \log \tilde h_k(0)\cdot \theta)\frac{e^{-\lambda_1^k}|x|}{\left(e^{-\lambda_1^k}+\frac{\tilde{h}_k(0)}{8(1+\alpha)^2}|x|^{2\alpha+2}\right)^2}.$$
	Putting $\lambda=\lambda_1^k$ and differentiating with respect to $\lambda$ we have that,
	$$\frac{d}{d\lambda}\hat c_1^k(x)
	=\frac{\tilde{h}_k(0)}{4\alpha(1+\alpha)}(\nabla \log \tilde h_k(0)\cdot \theta) \frac{e^{-\lambda}|x|^{2\alpha +3}}{\left(e^{-\lambda_1^k}+\frac{\tilde{h}_k(0)}{8(1+\alpha)^2}|x|^{2\alpha+2}\right)^2}.$$
	Setting
	$$w_{1,\xi}^k(y)=b_0^k\left(\frac{d}{d\lambda}\hat c_1^k\right)(\varepsilon_ky)\quad \mbox{with}
	\quad \varepsilon_k=e^{-\frac{\lambda_1^k}{2(1+\alpha )}}, $$ we have
	\begin{equation}\label{w1-explicit}
		w_{1,\xi}^k(y)=b_0^k\frac{\tilde{h}_k(0)}{4\alpha(1+\alpha)}(\nabla \log \tilde h_k(0)\cdot \theta)
		\frac{\varepsilon_k r^{2\alpha+3}}{\left(1+\frac{\tilde{h}_k(0)}{8(1+\alpha)^2}r^{2+2\alpha }\right)^2},
	\end{equation}
	and it is readily seen (see also \eqref{mar31e1})  that $w_{1,\xi}^k$ satisfies \eqref{xi-2nd}.\\
	
	Similarly, let $A_0^k(x)=b_0^kc_{0,k}(|x|/\varepsilon_k)$,
	then we have,
	\begin{align*}
		%\frac{\tilde{h}_k(0)}{8(1+\alpha)^2}
		%&\Delta A_0^k+\tilde h_k(0)\frac{8(1+\alpha)^2|x|^{2\alpha}e^{-\lambda_1^k}}{(e^{-\lambda_1^k}+|x|^{2\alpha+2})^2}A_0^k=-b_0^k\frac {\tilde h_k(0)}{4} 8(1+\alpha)^2\frac{e^{-\lambda_1^k}|x|^{2\alpha+2}}{(e^{-\lambda_1^k}+|x|^{2\alpha+2})^2} \\
		%&\bigg ((1-\frac{2(1+\alpha)}{\alpha}\frac{e^{-\lambda_1^k}}{e^{-\lambda_1^k}+|x|^{2\alpha}})|\nabla \log \tilde h_k(0)|^2+\Delta \log \tilde h_k(0)\bigg ).\nonumber
		&\Delta A_0^k+\tilde h_k(0)\frac{|x|^{2\alpha}e^{-\lambda_1^k}}{\left(e^{-\lambda_1^k}+\frac{\tilde{h}_k(0)}{8(1+\alpha)^2}|x|^{2\alpha+2}\right)^2}A_0^k\\
		&=-b_0^k\frac {\tilde h_k(0)}{4}\frac{e^{-\lambda_1^k}|x|^{2\alpha+2}}{\left(e^{-\lambda_1^k}+\frac{\tilde{h}_k(0)}{8(1+\alpha)^2}|x|^{2\alpha+2}\right)^2} \\
		&\quad\times\left(\left(1-\frac{2(1+\alpha)}{\alpha}\frac{e^{-\lambda_1^k}}{e^{-\lambda_1^k}+\frac{\tilde{h}_k(0)}{8(1+\alpha)^2}|x|^{2\alpha+2}}\right)^2|\nabla \log \tilde h_k(0)|^2+\Delta \log \tilde h_k(0)\right).
	\end{align*}
	At this point we define $A_{\lambda}^k(x)=\frac{d}{d\lambda}A_0^k(x)$ and then a lengthy evaluation
	shows that,
	\begin{equation}\label{important-z0}
		z_0^k(y)=A_{\lambda}^k(\varepsilon_k y).
	\end{equation}
	By using (\ref{important-z0}) we deduce that,
	\begin{align*}
		\int_{B_{\tau \varepsilon_k^{-1}}}\Delta z_0^k=~&\varepsilon_k^2\int_{B_{\tau \varepsilon_k^{-1}}}
		\Delta A_{\lambda}^k(\varepsilon_ky)dy\\
		=~&\int_{B_{\tau}}\Delta A^k_{\lambda}(x)dx=\frac{d}{d\lambda}\int_{B_{\tau}}\Delta A_0^k(x)dx.
	\end{align*}
	Since
	$$\Delta A_0^k(x)=b_0^k\Delta c_{0,k}\left(\frac{|x|}{\varepsilon_k}\right)\varepsilon_k^{-2},$$
	we have
	\begin{align*}
		\frac{1}{b_0^k}\int_{B_{\tau}}\Delta A_0^k(x)dx=&\int_{B_{\tau \varepsilon_k^{-1}}}\Delta c_{0,k}(y)dy\\
		=&
		\int_{\partial B_{\tau \varepsilon_k^{-1}}}\frac{\partial c_{0,k}}{\partial \nu}=
		d_{1,k} \Delta \log \tilde h_k(0)\varepsilon_k^2+O(\varepsilon_k^{2+\epsilon_0}), 
	\end{align*}
	and consequently, by a straightforward inspection of \eqref{add-exp},
	\begin{align*}
		\int_{\partial B_{\tau \varepsilon_k^{-1}}}\frac{\partial z_0^k}{\partial \nu}=
		\int_{B_{\tau \varepsilon_k^{-1}}}\Delta z_0^k=-\frac{1}{1+\alpha}b_0^kd_{1,k}
		\Delta (\log \tilde h_k)(0)\varepsilon_k^2+O(b_0^k\varepsilon_k^{2+\epsilon_0}),
	\end{align*}
	which, in view of \eqref{key-c-2} eventually implies that,
	\begin{equation}\label{p1-leading}
		\int_{\Omega(p_1,\tau)}\rho_kH c_k(x)\xi_k=-\frac{1}{1+\alpha}b_0^kd_{1,k}
		\Delta (\log \tilde h_k)(0)\varepsilon_k^2+O(b_0^k\varepsilon_k^{2+\epsilon_0}).
	\end{equation}
	\medskip

	On the other side, the contribution of the integral around $p_2$ is very small,
	\begin{equation}\label{p2-minor}
		\int_{\Omega(p_2,\tau)}\rho_kH c_k(x)\xi_k=O(\varepsilon_k^{2+\epsilon_0}),
	\end{equation}
	as can be seen by an argument similar to that adopted around $p_1$. Indeed,
	in the expansion of $\xi_k$ around $p_2$, the equation for $w_{1,\xi}^k$ also comes from
	the differentiation of $\lambda$ with respect to the $c_{1,k}$ in the expansion of $u_1^k$.
	The equation for $z_0^k$, which is the second order radial term in the expansion of $\xi_k$
	around $p_2$, also comes from that of $c_{0,k}$ in the expansion of $u_1^k$ around $p_2$.
	However the scaling is now with respect to $\hat \varepsilon_k=\varepsilon_k^{\frac{1+\alpha}{1+\beta}}$, whence
	$O(\hat \varepsilon_k^2)=O(\varepsilon_k^{2+\epsilon_0})$, implying immediately that \eqref{p2-minor} holds.
	The proof would be a cut and paste of the same argument adopted above for $p_1$ with minor changes,
	where one needs just to replace $\alpha$ by $\beta$, we skip these details to avoid repetitions.
	
	\medskip
	
	Concerning the integral on $\Omega(q,\tau)$, we first improve the estimate of
	$\bar w_k$ (see \eqref{pre-xi}). At first, we can show that on $\partial B_{\tau/\bar\varepsilon_k}$ 
	\begin{align*}
		\bar\psi_\xi^k(\varepsilon_ky_1)-\bar\psi_\xi^k(\varepsilon_ky_2)=~&\xi_k(\bar q_2^k+\bar\varepsilon_ky_1)-\bar w_{0,\xi}^k(y_1)-\xi_k(\bar q_2^k+\bar\varepsilon_ky_2)+\bar w_{0,\xi}^k(y_2)\\
		=~&O(\bar\varepsilon_k+\varepsilon_k^2),
	\end{align*}
	where we used Lemma \ref{osci-xi-better} and the explicit formula of $\bar w_{0,\xi}^k(y)$. As a consequence, by a standard argument we could derive an estimate like
	(\ref{im-ph-xi}) holds for $\bar \psi^k_\xi(\varepsilon_ky)$ as well, i.e.,
	\begin{equation}
		\label{4.equ-barpsi}
		|\bar\psi_\xi^k(\varepsilon_ky)|\leq C(\bar\varepsilon_k^2+\bar\varepsilon_k\varepsilon_k^2)(1+|y|),\quad y\in B_{\tau/\bar\varepsilon_k}.
	\end{equation}
	While at $\bar q_2^k$,  by standard elliptic estimate for harmonic function we have $|\nabla\bar\psi^k_\xi(0)|=O(\bar\varepsilon_k^2+\bar\varepsilon_k\varepsilon_k^2)$. Therefore, in local
	coordinates such that $0=x(q)$ and after scaling $x=\bar q_2^k+\bar \varepsilon_k y$,
	where $\bar \varepsilon_k=e^{-\bar\lambda_1^k/2}$,
	we can write
	the equation for $\bar w_k$ as follows,
	$$\Delta \bar w_k+\bar h_k(\bar q_2^k+\bar \varepsilon_k y) \bar \varepsilon_k^{2}
	\bar c_k(\bar q_2^k+\bar \varepsilon_k y)\bar w_k=
	O(\varepsilon_k^{2+\epsilon_0})(1+|y|)^{-3}
	\quad {\rm in} \quad B_{\tau \bar \varepsilon_k^{-1}}.$$
	Here we remark that, compared with \eqref{eqwbar1}, the improvement in the estimates of the
	right hand side is obtained because of \eqref{c}, the vanishing rate of the gradient of the coefficient function
	for regular blowup points (see \eqref{first-deriv-est}),
	the improved estimate \eqref{4.equ-barpsi} for $\bar\psi_{\xi}^k$ and the fact that
	$O(\bar \varepsilon_k^2)=O(\varepsilon_k^{2+\epsilon_0})$. On the other hand, we have 
	$$\bar w_k(0)=0,\quad \nabla\bar w_k(0)=O(\varepsilon_k^{2+\epsilon_0})\quad {\rm and}\quad \bar w_k \mbox{ is a constant on } \partial B_{\tau \bar \varepsilon_k^{-1}}.$$ 
	As a consequence, by the usual potential estimates, we conclude that,
	\begin{equation*}
		|\bar w_k(y)|\le C(\delta)\varepsilon_k^{2+\epsilon_0}(1+|y|)^{\delta},\quad |y|\le \tau \bar \varepsilon_k^{-1}.
	\end{equation*}
	Next, by using the expansion of $\tilde c_k$ in \eqref{c} we have,
	$$
	\int_{\Omega(q,\tau)}\rho_kH (x)c_k(x)\xi_k
	=\int_{B_{\tau \bar \varepsilon_k^{-1}}}\bar h_k(\bar q_2^k+\bar \varepsilon_k y)e^{U_2^k}
	\xi_k(\bar q_2^k+\bar \varepsilon_k y)
	+O(\varepsilon_k^{2+\epsilon_0}).
	$$
	Remark that $e^{-\lambda_1^k/2}=O(\varepsilon_k^{1+\epsilon_0})$, then by the expansion of $\xi_k$ in (\ref{pre-xi}) we see that
	all the terms including $\bar \psi^k_\xi(\bar \varepsilon_ky)$ are of order
	$O(\varepsilon_k^{2+\epsilon_0})$.
	Also, neglecting terms which vanish due to the separability of the integrand we have that,
	\begin{align*}
		&\int_{B_{\tau \bar \varepsilon_k^{-1}}}\bar h_k(\bar q_2^k +\bar \varepsilon_k y)e^{U_2^k}
		\xi_k(\bar q_2^k+\bar \varepsilon_k y)\\
		&=\int_{B_{\tau \bar \varepsilon_k^{-1}}}\bar h_k(\bar q_2^k)e^{U_2^k}(\bar w^k_{0,\xi}+\bar w_k)+\bar \varepsilon_k\int_{B_{\tau \bar \varepsilon_k^{-1}}}(\nabla \log\bar h_k(\bar q_2^k) \cdot y)
		\frac{e^{U_2^k}(\bar b_k \cdot y)}{1+\frac{\bar h_k(\bar q_2^k)}{8}|y|^2}\\
		&\quad+\bar \varepsilon_k\int_{B_{\tau \bar \varepsilon_k^{-1}}}(\nabla \log\bar h_k(\bar q_2^k) \cdot y)e^{U_2^k}\bar w_k+O(\varepsilon_k^{2+\epsilon_0})\\
		&=\int_{B_{\tau \bar \varepsilon_k^{-1}}}\bar h_k(\bar q_2^k)e^{U_2^k}\bar w_k+
		\bar \varepsilon_k\int_{B_{\tau \bar \varepsilon_k^{-1}}}(\nabla \log\bar h_k(\bar q_2^k) \cdot y) 
		\frac{e^{U_2^k}(\bar b_k \cdot y)}{1+\frac{\bar h_k(\bar q^k_2)}{8}|y|^2}+O(\varepsilon_k^{2+\epsilon_0})\\
		&=O(\varepsilon_k^{2+\epsilon_0}),
	\end{align*}
	where $\bar b_k=(\bar b^k_1, \bar b^k_2)$ and we used the same argument as in \eqref{i1a} to show that
	$\int_{B_{\tau \bar \varepsilon_k^{-1}}}\bar h_k(\bar q^k_2)e^{U_2^k}\bar w^k_{0,\xi}$ is of order
	$O(\varepsilon_k^{2+\epsilon_0})$.
	Therefore we eventually deduce that,
	\begin{equation}\label{small-q}
		\int_{\Omega(q,\tau)}\rho_kH c_k(x)\xi_k=O(\varepsilon_k^{2+\epsilon_0}).
	\end{equation}
	In view of \eqref{p1-leading}, \eqref{p2-minor}, \eqref{small-q}, we come up with
	a contradiction to $$\int_M\rho_k H c_k(x){\xi_k}=0$$ as follows,
	\begin{align*}
		0&=\int_M\rho_k H c_k\xi_k\\
		&=\int_{\Omega(p_1,\tau)}H c_k\xi_k+\int_{\Omega(p_2,\tau)}H c_k\xi_k+\int_{\Omega(q,\tau)}H c_k\xi_k\\
		&\quad+\int_{M\setminus \{\Omega(p_1,\tau)\cup \Omega(p_2,\tau)\cup \Omega(q,\tau)\}}H c_k\xi_k\\
		&=Cb_0\Delta \log h_k(p_1)\varepsilon_k^2+o(\varepsilon_k^2),
	\end{align*}
	for some constant $C\neq 0$, since in particular the integrals on $\Omega(p_2,\tau)$, $\Omega(q,\tau)$
	and $M\setminus \{\Omega(p_1,\tau)\cup \Omega(p_2,\tau)\cup \Omega(q,\tau)\}$ are all of order
	$o(\varepsilon_k^2)$.
	Since for $k$ large $L(\mathbf{p})\neq 0$ is the same as $\Delta \log h_k(p_1)\neq 0$,
	we obtain a contradiction as far as $b_0\neq 0$.
	\begin{rem}\label{rem 4.8} To summarize the idea of the proof of $b_0=0$. The assumption $L(\mathbf{p})\neq 0$ plays an important role. First we need it to prove that $\sigma_k=O(\varepsilon_k^{2+\epsilon_0})$. This estimate makes it possible to have a good approximation of $\xi_k$ around each blowup point. Then we use $L(\mathbf{p})\neq 0$ again together with $\int_M\rho_kHc_k\xi_k=0$ to get a contradiction.
	\end{rem}
	\medskip

	\subsection{Proof of $b_1=b_2=0$}\label{pf-uni-2}
	First we remark that in all the previous works on uniqueness of bubbling solutions, proving $b_1=b_2=0$ after the establishment of $b_0=0$ usually requires little new ideas as well as effort. But in this case we need to establish major new estimates and the proof is divided in multiple non-trivial steps. At the beginning of this part we 
	recall that
	$u^k_i$ $i=1,2$ satisfy \eqref{uk-def1}, \eqref{around-q},
	in local coordinates around $q$, $0=x(q)$, 
	\begin{equation*}
		\Delta u_i^k+\bar h_k(x)e^{u_i^k}=0\quad \mbox{in}\quad B_\tau,\quad i=1,2,
	\end{equation*}
	and that $\bar q_i^k$ denote the local maximum points of $u_i^k$ $i=1,2$. Of course we have 
	$|\bar q_1^k-\bar q_2^k|=O(\lambda_1^ke^{-\lambda_1^k})$.
	As remarked right after \eqref{ueqw} we also have that,
	$$
	\xi_k=\frac{u_1^k-u_2^k}{\sigma_k}=\frac{{\rm w}_1^k- {\rm w}_2^k}{\sigma_k},
	$$
	where $\xi_k$ satisfies,
	\begin{equation*}
		\Delta \xi_k+\bar h_k \tilde c_k\xi_k=0\quad {in}\quad B_\tau.
	\end{equation*}
	The expression of $\tilde c_k$ is
	\begin{equation*}
		\tilde c_k(x)=e^{u_2^k}\int_0^1e^{t(u_1^k-u_2^k)}dt=e^{u_2^k}\left(1+\frac 12\sigma_k\xi_k+\frac 16\sigma_k^2\xi_k^2+O(\sigma_k^3)\right).
	\end{equation*}
	One major consequence of $b_0=0$ is that, since after scaling around a blow up point
	the local limit of $\xi_k$ takes the form \eqref{linim-1}, \eqref{linim-2}, \eqref{linim-3} respectively,
	then it is not difficult to see that,
	\begin{equation*}
		\frac{\lambda_1^k-\lambda_2^k}{\sigma_k}\to 0,\quad \frac{\tilde \lambda_1^k-\tilde \lambda_2^k}{\sigma_k}\to 0,
		\quad \frac{\bar \lambda_1^k-\bar \lambda_2^k}{\sigma_k}\to 0.
	\end{equation*}
	Also, see Remark \ref{xi:b_0}, at this point we know that, far
	away from blow up points, $\xi_k$ is of order $o(1)$. This fact obviously implies that,
	in local coordinates around any blow up point, defining $\phi_i^k$ $i=1,2$ to be the harmonic functions
	encoding the oscillations of $u_i^k$, then we have,
	\begin{equation*}
		(\phi_1^k-\phi_2^k)/\sigma_k=o(1).
	\end{equation*}
	On the other side, we cannot anymore rely on \eqref{sigma:eq} which was deduced
	by $b_0\neq 0$, whence neither we can rely on Lemma \ref{small-sigma} or on the
	improved estimates \eqref{osi-vars}. This is why we have to derive
	first a new estimate about $\sigma_k$, this time based on $b_0=0$. On the other side,
	since our current goal is to prove $b_1=b_2=0$, we can assume without loss of generality that
	at least one of them is not zero. 
	
	\medskip
	
	\noindent{\bf Step one: An initial estimate of $\sigma_k$}
	\begin{lem}\label{new-sigma}
		Assume that $(b_1,b_2)\neq (0,0)$, then 
		$$\sigma_k=O(\lambda_1^ke^{-\frac{\lambda_1^k}2}).$$
	\end{lem}
	
	\begin{proof}[Proof of Lemma \ref{new-sigma}.]
Let
	$$
	\bar {u}^k_i=u_i^k-\phi_i^k,
	$$
	where, in local coordinates around $q$, $\phi_i^k$ $i=1,2$ are the harmonic functions
	encoding the oscillations, of $u_i^k$ on $B_\tau(\bar q_i^k)$, whence satisfying, $\phi_i^k(\bar q_i^k)=0$.
	Thus, working in local coordinates around $q$
	such that $0=x(q)$, since $u_i^k$, $i=1,2$, satisfy \eqref{around-q} we can write
	\begin{equation}\label{v2k-q}
		\Delta \bar u_2^k+\bar h_{k,0}e^{\bar u_2^k}=0\quad {in}\quad B_\tau
	\end{equation}
	with
	\begin{equation}\label{bar-hk-q}
		\bar h_{k,0}(x)=\bar h_k(x)e^{\phi_2^k(x)}.
	\end{equation}
	Therefore the equation for $\bar u_1^k$ reads,
	\begin{equation}\label{v1k-q}
		\Delta \bar u_1^k+\bar h_{k,0}e^{\phi_1^k-\phi_2^k}e^{\bar u_1^k}=0\quad {in}\quad B_\tau
	\end{equation}
	The reason that we are not introducing a new coefficient function for $\bar u_1^k$ is because we use $\bar u_2^k$ as the base to carry out the proof and we will focus on the coefficient function $\bar h_{k,0}$ of $\bar u_2^k$. 
	With these definitions we have
	\begin{equation}\label{xi-v}
		\xi_k=\frac{\bar u_1^k-\bar u_2^k}{\sigma_k}+\bar\psi_\xi^k,\quad
		\bar\psi_\xi^k=\frac{\phi_1^k-\phi_2^k}{\sigma_k}
	\end{equation}
	and the expression of $c_k$ near $q$ becomes,
	\begin{equation*}
		\rho_kH \tilde c_k=\bar h_{k,0} e^{\bar u_2^k}\left(1+\frac12\sigma_k\xi_k+O(\sigma_k^2)\right).
	\end{equation*}
	Let $\bar \phi_1^k$ be the harmonic function on $B_\tau(\bar q_2^k)$ that encodes the oscillation of $u_1^k$ on $\partial B_\tau(\bar q_2^k)$:
	\begin{equation*}
		\begin{cases}
			\Delta \bar \phi_1^k=0,\quad &{\rm in}\quad B_\tau(\bar q_2^k),\\
			\\
			\phi_1^k(x)= u_1^k(x)-\frac{1}{2\pi\tau}\int_{\partial B_\tau(q_2^k)} u_1^k,  &{\rm on}\quad \partial B_\tau(\bar q_2^k),
		\end{cases}
	\end{equation*}
	then we observe that $(\bar \phi_1^k-\phi_2^k)/\sigma_k=o(1)$ in $B_\tau(\bar q_2^k)$ because of the following boundary estimate,
	\begin{align*}\frac{(\bar \phi_1^k-\phi_2^k)(x)}{\sigma_k}&=\frac{u_1^k(x)-u_2^k(x)}{\sigma_k}-\frac{1}{2\pi \tau}
		\int_{\partial B_\tau(\bar q_2^k)}\frac{u_1^k- u_2^k}{\sigma_k}\\
		&=\xi_k(x)-\frac{1}{2\pi \tau}\int_{\partial B_\tau(\bar q_2^k)}\xi_k=o(1), \quad x\in \partial B_\tau(\bar q_2^k).
	\end{align*}
	
	Next we claim that
	\begin{equation}\label{har-small}
		\bar \phi_1^k(x)-\phi_1^k(x)=O(\lambda_1^k e^{-\lambda_k}),\quad x\in B_{\frac{\tau}2}(\bar q_2^k).
	\end{equation}
	Indeed, by the definition of $\bar \phi_1^k$ we have,
	\begin{equation}\label{bar-phi-1-s}\bar \phi_1^k(x)=-\int_{\partial B_\tau(\bar q_2^k)}\partial_{\nu}G_2(x,\eta)\bar \phi_1^k(\eta)dS_{\eta}, \quad x\in B_{\tau/2}(\bar q_2^k),
	\end{equation}
	where $G_2$ is the Green function on $B_\tau(\bar q_2^k)$. Similarly for $\phi_1^k$ we have
	$$\phi_1^k(x)=-\int_{\partial B_\tau(\bar q_1^k)}\partial_{\nu}G_1(x,\eta)\phi_1^k(\eta)dS_{\eta},\quad x\in B_{\tau/2}(\bar q_2^k)$$
	where $G_1$ is the Green's function on $B_\tau(\bar q_1^k)$.
	After a translation that moves $B_\tau(\bar q_1^k)$ to $B_\tau(\bar q_2^k)$ we see that,
	\begin{equation}\label{bar-phi-small}
		\phi_1^k(x)=-\int_{\partial B_\tau(\bar q_2^k)}\partial_{\nu}G_2(x,\eta) \phi_1^k(\eta+\bar q_1^k-\bar q_2^k)dS_{\eta},\quad x\in B_{\tau/2}(\bar q_2^k).
	\end{equation}
	From (\ref{bar-phi-1-s}) and (\ref{bar-phi-small}) we see that $$|\bar \phi_1^k(x)-\phi_1^k(x)|\le C\int_{\partial B_{\tau}(\bar q_2^k)}| \phi_1^k(\eta+\bar q_1^k-\bar q_2^k)-\bar \phi_1^k(\eta))|dS_{\eta}\le C\lambda_1^ke^{-\lambda_1^k}$$
	for $x\in B_{\tau/2}(\bar q_2^k)$ because $|\nabla u_1^k|=O(1)$ outside bubbling area and 
	$$|\bar q_2^k - \bar q_1^k|=O(\lambda_1^ke^{-\lambda_1^k}).$$ 
	Hence, \eqref{har-small} is established. As a consequence, as far as $\sigma_k/(\lambda_1^ke^{-\lambda_1^k})\to \infty$, we have
	$$\frac{\phi_1^k-\phi_2^k}{\sigma_k}=o(1)\quad {\rm in}\quad B_{\tau/2}(\bar q_2^k).$$
	
	Here we also note that the maximum of $\bar u_i^k$ is $\bar q_i^k+O(e^{-\lambda_1^k})$ (see \cite{zhang1}). 
	It is readily seen that the corresponding error is negligible for our purposes, which is why, 
	with an abuse of notations, we will just use $\bar q_i^k$ to denote the local maximum of $\bar u_i^k$ as well. 
	Recall that the weight function in \eqref{v2k-q} is $\bar h_{k,0}$ whence the weight in \eqref{v1k-q} is
	$\bar h_{k,0}e^{\phi_1^k-\phi_2^k}$. Since 
	$$\phi_i^k(\bar q_i^k)=0
	\quad\mathrm{and}\quad
	|\bar q_1^k-\bar q_2^k|=O(\lambda_1^ke^{-\lambda_1^k}),$$ 
	we have
	$$\bar h_{k,0}(\bar q_1^k)e^{(\phi_1^k-\phi_2^k)(\bar q_1^k)}=\bar h_{k,0}(\bar q_2^k)+O(\lambda_1^ke^{-\lambda_1^k}).
	$$
	Thus, by the main theorem in \cite{gluck} we can write $\bar u_i^k$ as follows,
	$$\bar u_i^k=\bar\lambda_i^k-2\log (1+\frac{e^{\bar \lambda_i^k}\bar h_{k,0}(\bar q_i^k)}8|x-\bar q_i^k|^2)+
	O((\lambda_1^k)^2e^{-\lambda_1^k}),\quad i=1,2.$$
	To simplify the exposition let us set $h_i=\bar h_{k,0}(\bar q_i^k)$, then we have,
	$$
	\bar u_1^k-\bar u_2^k=\bar \lambda_1^k-\bar \lambda_2^k+
	2\log \left(1+\frac{e^{\bar \lambda_2^k}\frac{h_2}8|x-\bar q_2^k|^2-e^{\bar \lambda_1^k}\frac{h_1}8|x-\bar q_1^k|^2}
	{1+e^{\bar \lambda_1^k}\frac{h_1}8|x-\bar q_1^k|^2}\right)+O((\lambda_1^k)^2e^{-\lambda_1^k}). $$
	Denoting as usual $\bar \varepsilon_k=e^{-\bar \lambda_1^k/2}$ we have
	\begin{align*}
		&\bar u_1^k(\bar q_1^k+\bar\varepsilon_ky)-\bar u_2^k(\bar q_1^k+\bar\varepsilon_ky)\\
		&=\bar \lambda_1^k-\bar \lambda_2^k+2\log \left(1+\frac{\frac{e^{\bar \lambda_2^k-\bar \lambda_1^k}h_2}
			8|\frac{\bar q_2^k-\bar q_1^k}{\bar \varepsilon_k}+y|^2-\frac{h_1}{8}|y|^2}{1+\frac{h_1}8|y|^2}\right)
		+O((\lambda_1^k)^2e^{-\lambda_1^k}).
	\end{align*}
	In view of $\log (1+A)=A+O(A^2)$, we define,
	\begin{align*}
		A&=\frac{e^{\bar \lambda_2^k-\bar \lambda_1^k}\frac{h_2}8|\frac{\bar q_1^k-\bar q_2^k}{\bar \varepsilon_k}+y|^2-
			\frac{h_1}8|y|^2}{1+\frac{h_1}8|y|^2}\\
		&=\frac{(1+\bar \lambda_2^k-\bar \lambda_1^k+O(|\bar \lambda_2^k-\bar \lambda_1^k|^2))\frac{h_2}8(|y|^2+2a\cdot y+|a|^2)-
			\frac{h_1}8|y|^2}{1+\frac{h_1}8|y|^2}
	\end{align*}
	where $a=\frac{\bar q_2^k-\bar q_1^k}{\bar \varepsilon_k}$. At this point, as far as
	$\sigma_k/(\bar \lambda_1^ke^{-\bar \lambda_1^k/2})\to +\infty$, combining the leading terms we have,
	\begin{align*}
		\frac{\bar u_1^k(\bar q_1^k+\bar \varepsilon_ky)-\bar u_2^k(\bar q_1^k+\bar \epsilon_k y)}{\sigma_k}=\frac{\bar\lambda_1^k-\bar\lambda_2^k}{\sigma_k}\frac{1-\frac{h_1}8|y|^2}{1+\frac{h_1}8|y|^2}
		+\frac{h_1}{2}\frac{\frac{a}{\sigma_k}\cdot y}{1+\frac{h_1}8|y|^2}+o(1).
	\end{align*}
	Obviously the radial leading term vanishes in the limit due to $$(\bar\lambda_1^k-\bar\lambda_2^k)/\sigma_k\to 0~(\mbox{that is }b_0=0),$$
	while 
	$b_1,b_2$ should arise from the term
	$$\frac{h_1}2\frac{\frac{\bar q_1^k- \bar q_2^k}{\bar \varepsilon_k\sigma_k}\cdot y}{(1+\frac{h_1}8|y|^2)}.$$
	However, by using again $\sigma_k/(\bar \lambda_1^ke^{-\bar \lambda_1^k/2})\to \infty$, we deduce that
	$b_1=b_2=0$ because of,
	$$\frac{|\bar q_1^k-\bar q_2^k|}{\bar \varepsilon_k}=O(\bar \lambda_1^ke^{-\bar \lambda_1^k/2}).$$
	A contradiction to $(b_1,b_2)\neq (0,0)$. Thus the Lemma \ref{new-sigma} is proved. 
\end{proof}
	
	\noindent{\bf Step two: Intermediate estimates for $\bar v_2^k$}
	
	\medskip
	
	We work in local coordinates centered at $q$, $0=x(q)$, with $\bar \varepsilon_k=e^{-\bar \lambda_2^k/2}$ and define
	\begin{equation}\label{bar-v-2-q}
		\bar v_2^k(y)=\bar u_2^k(\bar q_2^k+\bar \varepsilon_ky)+2\log \bar \varepsilon_k,\quad  y\in  \Omega_k,
	\end{equation}
	and 
	\begin{equation}\label{xi-k-q}
		\bar \xi_k(y)=\xi_k(\bar q_2^k+\bar \varepsilon_ky),\quad y \in \Omega_k:=B_{\tau \bar \varepsilon_k^{-1} }.
	\end{equation}
	Note that $\bar v_2^k$ has no oscillation on $\partial \Omega_k$. If we invoke Lemma \ref{intermediate-2} for $\bar v_2^k$ we have 
	\begin{equation}\label{v2-exp-new2}
		\bar v_2^k(y)=U_2^k+c_{0,k}+c_{2,k}+O(\bar \varepsilon_k^{2+\delta})(1+|y|)^{\delta}.
	\end{equation}
	In view of Lemma \ref{new-sigma} we can write
	\begin{equation*}
		\bar \varepsilon_k^2 \rho_kH \tilde c_k(\bar q_2^k+\bar \varepsilon_ky)\xi_k=\bar h_{k,0}(\bar q_2^k+\bar \varepsilon_ky)e^{U_2^k}
		\bar \xi_k\left(1+\frac 12\sigma_k \bar \xi_k+O(\sigma_k^2)\right).
	\end{equation*}
	
	Here we point out that the function which we called $u_k$ in Lemma \ref{intermediate-2} is
	what we call $\bar u^k_i$, $i=1,2$ here, which satisfy \eqref{v2k-q}, \eqref{bar-hk-q}. Also remark that
	$\bar \lambda_i^k=\bar u^k_i(\bar q^k_i)$, $i=1,2$ and $\bar \varepsilon_k=e^{-\bar \lambda_2^k/2}$.
	Now we set
	$$
	w_2^k(y)=\bar v_2^k(y)-c_{0,k}(y)-c_{2,k}(y),
	$$
	where $c_{i,k}$ $i=0, 2$, are defined as in the proof of Theorem \ref{reg-int}. Here we explain why there is no $c_{1,k}$ term. The reason is $\bar h_{k,0}$ satisfies $\nabla \bar h_{k,0}(\bar q_2^k)=O(\bar \lambda_1^ke^{-\bar \lambda_1^k})$. The estimate of $c_{1,k}$ would be 
	$$|c_{1,k}(y)|\le C\bar \lambda_1^ke^{-\frac 32\bar \lambda_1^k}\log (1+|y|)(1+|y|)^{-1},\quad \mbox{in}\quad \Omega_k, $$
	and 
	$$|\nabla c_{1,k}(y)|\le C\bar\varepsilon_k^{3+\epsilon_0},\quad y\in \Omega_k,\quad |y|\sim \bar\varepsilon_k^{-1}.$$
	These two estimates show that at this stage the term $c_{1,k}$ is negligible.Indeed 
	$w_2^k$ satisfies,
	\begin{align*}
		\begin{cases}
			\Delta w_2^k+\bar h_{k,0}(\bar q_2^k)e^{\bar v_2^k}w_2^k=O(\bar \varepsilon_k^3)(1+|y|)^{-1},\\
			\\
			w_2^k(0)=|\nabla w_2^k(0)|=0,
		\end{cases}   
	\end{align*}
	and the oscillation of $w_2^k$ on $\partial \Omega_k$ is of order $O(\bar \varepsilon_k^2)$. Now we use Lemma \ref{intermediate-2} to obtain 
	$$|w_2^k(y)|\le C\bar\varepsilon_k^{2+\epsilon_0}(1+|y|)^{2\epsilon_0},$$
	which is equivalent to 
	\begin{equation}\label{for-v2-q}
		|\bar v_2^k-U_2^k-c_{0,k}-c_{2,k}|\le C\bar \varepsilon_k^{2+\epsilon_0}(1+|y|)^{2\epsilon_0},\quad y\in \Omega_k. 
	\end{equation}
	The gradient estimate on the boundary is
	\begin{equation}\label{grad-v2-new-2}\nabla \bar v_2^k=\nabla U_2^k+\nabla c_{0,k}+\nabla c_{2,k}+O(\tau)\bar\varepsilon_k^3,\quad y\in \Omega_k, \quad |y|\sim \bar\varepsilon_k^{-1}.
	\end{equation}
	
	We shall also need the following estimate based on (\ref{v2-exp-new2}),
	\begin{equation}
		\label{v2-exp-new}
		\begin{aligned}
			\bar v_2^k(y)+\log \frac{\bar h_{k,0}(\bar q_2^k+\bar \varepsilon_ky)}{\bar h_{k,0}(\bar q_2^k)}
			=~&U_2^k+c_{0,k}+\tilde c_{2,k}+\frac{\Delta (\log \bar h_{k,0})(\bar q_2^k)\bar \varepsilon_k^2r^2}{4}\\
			&+O(\bar \varepsilon_k^3)r^3+O(\bar \varepsilon_k^{2+\delta})(1+r)^{\delta}, 
		\end{aligned}
	\end{equation}
	We will use (\ref{v2-exp-new}) to derive a better estimate for $|\bar q_1^k-\bar q_2^k|$. 
	
	\medskip
	
	\noindent{\bf Step three: Initial estimates of the oscillation of $\xi_k$ far away from blow up points}

	\begin{lem}\label{small-osci-xi}
		For any $x_1,x_2\in M\setminus \left\{\Omega({p_1},\tau)\cup \Omega({p_2},\tau)\cup \Omega({q},\tau)\right\}$, there exists $C>0$ independent of $k$ and $\tau$ that
		$$|\xi_k(x_1)-\xi_k(x_2)|\le C\varepsilon_k $$
	\end{lem}
	
	\begin{proof}
	Since $b_0=0$ then $\xi_k\to 0$ in $C^2_{\rm loc}(M\setminus \{p_1,p_2,q\})$. Based on this information we consider again
	the expansion of $\xi_k$ around $p_1$, $p_2$ and $q$. In local coordinates around $p_1$, $0=x(p_1)$, we define 
	$$\tilde \xi_k(y)=\xi_k(p_1+\varepsilon_ky)$$ which
	satisfies
	\begin{equation}\label{around-p_1}
		\Delta \tilde \xi_k+|y|^{2\alpha}\tilde h_k(\varepsilon_ky)e^{v_2^k}\tilde \xi_k=
		\frac{O(\sigma_k)}{(1+|y|)^{4+2\alpha}} \quad {\rm in}\quad \Omega_k=
		B_{\tau/\varepsilon_k},
	\end{equation}
	with $\tilde h_k$ defined in \eqref{around-p1-h}.
	Our aim is to derive a first approximation of $\tilde \xi_k$. Let $\tilde \psi_{\xi}^k(\varepsilon_ky)$
	be the harmonic function that satisfies
	$\tilde \xi_k-\tilde \psi^k_{\xi}(\varepsilon_k\cdot)=$constant on $\partial \Omega_k$ and $\tilde \psi^k_{\xi}(0)=0$.
	Then it is easy to see that
	$$\tilde \xi_{1,k}(y):=\tilde \xi_k(y)-\tilde \psi^k_{\xi}(\varepsilon_k y)$$
	satisfies
	$$(\Delta+|y|^{\alpha}\tilde h_k(0)e^{U_2^k})\tilde \xi_{1,k}=(o(\varepsilon_k)+O(\sigma_k) )(1+|y|)^{-2\alpha-3}.$$
	To obtain a first order approximation of $\tilde \xi_{1,k}$ we use
	$$\tilde w_{0,\xi}^k=\tilde b_0^k\frac{1-a_k|y|^{2+2\alpha}}{1+a_k|y|^{2+2\alpha}},\quad a_k=\frac{\tilde h_k(0)}{8(1+\alpha)^2},$$
	where we choose $\tilde b_0^k$ such that $\tilde \xi_{1,k}(0)-w_{0,\xi}^k(0)=0$.  Note that we have $\tilde b_0^k\to 0$ but we don't have yet any information about the rate of vanishing of $\tilde b_0^k$.
	By the usual standard elliptic estimates we have,
	\begin{equation}\label{pre-1}
		|\tilde \xi_{1,k}(y)-\tilde w_{0,\xi}^k(y)|\leq C(\delta)(o(\varepsilon_k)+\sigma_k)(1+|y|)^{\delta}\quad {in}\quad \Omega_k.
	\end{equation}
	Similar estimates can also be obtained around $p_2$ and $q_2^k$. Around $p_2$ the scaling factor is
	$\hat \varepsilon_k=e^{-\frac{\lambda_1^k}{2(1+\beta)}}$ and, in local coordinates around $p_2$, $0=x(p_2)$,
	we define 
	$$\hat \xi_{1,k}=\xi_k(p_2+\hat \varepsilon_ky)-\hat \psi_k(\hat \varepsilon_ky),$$
	with the corresponding obvious meaning
	of $\hat \psi_k$, and
	$$\hat w_{0,\xi}^k(y)=\hat b_0^k\frac{1-\hat a_k|y|^{2+2\beta}}{1+\hat a_k|y|^{2+2\beta}},\quad \hat a_k=\frac{\hat h_k(0)}{8(1+\beta)^2}$$
	to be the first term in approximation of $\hat \xi_{1,k}$, where $\hat b_0^k$ is chosen so that $\hat \xi_{1,k}(0)-\hat w_{0,\xi}^k(0)=0$ and $\hat h_k$ defined in \eqref{around-p2-h}.  Thus,
	by arguing as above we have that,
	\begin{equation}\label{1st-p2}
		|\hat \xi_{1,k}-\hat w_{0,\xi}^k(y)|\le C(\delta)\sigma_k(1+|y|)^{\delta},\quad |y|\leq \tau \hat \varepsilon_k^{-1}.
	\end{equation}
	At last, around $\bar q_2^k$, the equation for $\bar \xi_k(y)$ ( see \eqref{xi-k-q} ) reads, 
	\begin{equation*}
		\Delta \bar \xi_k(y)+\bar h_{k,0}(\bar q_2^k)\exp\left(\bar v_2^k+\log \frac{\bar h_{k,0}(\bar q_2^k+\bar \varepsilon_ky)}{\bar h_{k,0}(\bar q_2^k)}\right) \bar \xi_k(y)=\frac{O(\sigma_k)}{(1+|y|)^{4}},
	\end{equation*}
	for $|y|\leq \tau \bar \varepsilon_k^{-1}$. Let us recall again (see \eqref{gradhk:gluck}) 
	that,
	$$\nabla \log \bar h_{k,0}(\bar q_2^k)=O(\lambda_1^ke^{-\lambda_1^k}). $$
	We denote by $\bar \xi_{1,k}$ the function $\xi_{1,k}$ modified to have vanishing oscillation and let
	$$\bar w_{0,\xi}^k(y)=\bar b_0^k\frac{1-\bar a_k|y|^2}{1+\bar a_k|y|^2}+b_1^k\frac{y_1}{1+\bar a_k|y|^2}+b_2^k\frac{y_2}{1+\bar a_k|y|^2}$$
	be the first term in the corresponding approximation, where $\bar a_k=\frac{\bar h_k^*(0)}{8}$ and $\bar b_0^k\to 0$, $b_1^k\to b_1$, $b_2^k\to b_2$
	are such that
	$$
	\bar \xi_{1,k}(0)=\bar w_{0,k}(0),\quad |\nabla (\bar \xi_{1,k}-\bar w_{0,k})(0)|=0.
	$$
	Here we note that the non-radial part of $\bar w_{0,\xi}^k$ has an oscillation of order $O(\bar \varepsilon_k)$ on $\partial \Omega_k$.
	However this oscillation can be removed by another harmonic function whose modulus is not larger than $O(\bar \varepsilon_k^2)|y|$,
	which is a minor contribution at this point. Thus we have
	\begin{equation}\label{1-around-q}
		|\bar \xi_{1,k}-\bar w_{0,\xi}^k|\le C(\delta)\sigma_k (1+|y|)^{\delta},\quad |y|\leq \tau \bar \varepsilon_k^{-1}.
	\end{equation}
	By using (\ref{pre-1}),(\ref{1st-p2}) and (\ref{1-around-q}) we derive the desired estimate of the oscillation of $\xi_k$ far away
	from blow up points. For any $x_1,x_2\in M\setminus (\Omega({p_1},\tau) \cup \Omega({p_2},\tau)\cup \Omega({q},\tau))$ we use
	the Green representation formula for $\xi_k$,
	\begin{equation}\label{osci-again}
		\xi_k(x_1)-\xi_k(x_2)=\int_M(G(x_1,\eta)-G(x_2,\eta))\rho_kHc_k\xi_k(\eta)d\eta.
	\end{equation}
	It is enough to evaluate the r.h.s of (\ref{osci-again}) on the three local disks surrounding the blow up points, since
	the integral over the remaining region is of order $O(e^{-\lambda_1^k})$. Concerning the integral over $\Omega(p_1,\tau)$ we have
	\begin{align*}
		&\int_{\Omega(p_1,\tau)}(G(x_1,\eta)-G(x_2,\eta))\rho_kHc_k\xi_kd\eta\\
		&=G_{12}(x_{12},p_1)\int_{B_{\tau}}|x|^{2\alpha}\tilde h_ke^{u_2^k}\xi_kd\eta\\
		&\quad+\int_{B_{\tau}}(G_{12}(x_{12},\eta)-G_{12}(x_{12},p_1))
		|x|^{2\alpha}\tilde h_ke^{u_2^k}\xi_k d\eta+O(\sigma_k),
	\end{align*}
	where $G_{12}(x_{12},\eta)=G(x_1,\eta)-G(x_2,\eta)$. Recall that (see \eqref{i1a})
	by using the equation satisfied by $\tilde w_{0,\xi}^k$ and the fact that $\partial_{\nu}\tilde w_{0,\xi}^k=O(r^{-3-2\alpha})$,
	then one can see that,
	$$
	\int_{B_{\tau}}|x|^{2\alpha}\tilde h_ke^{u_2^k}\tilde w_{0,\xi}^k=O(\varepsilon_k^{2+\epsilon_0}),
	$$
	for some $\epsilon_0>0$.
	Since the error estimates in (\ref{pre-1}) is of order $O(\varepsilon_k)$, then it is readily
	seen as well that the integral over $\Omega(p_1,\tau)$ is $O(\varepsilon_k)$.
	The integrals around $p_2$ and $q$ are of the order $O(\lambda_1^ke^{-\lambda_1^k/2})$.
	Lemma \ref{small-osci-xi} is established. \end{proof}

	\noindent{\bf Step four: a crucial new estimate of $|\bar q_1^k-\bar q_2^k|$}
	
	\begin{lem}\label{crucial-q12}
		$$|\bar q_1^k-\bar q_2^k|\le C\bar \varepsilon_k^2$$
		for $C>0$ independent of $k$. 
	\end{lem}
	\begin{proof}
	We consider the following Pohozaev identity for the equation of $\bar u_2^k$: for any $\xi\in \mathbb S^1$,
	$$\int_{B_\tau(\bar q_2^k)}\partial_{\xi}\bar h_{k,0}e^{\bar u_2^k}=\int_{\partial B_\tau(\bar q_2^k)}
	(\bar h_{k,0}e^{\bar u_2^k}(\xi\cdot \nu)+\partial_{\nu}\bar u_2^k\partial_{\xi}\bar u_2^k-\frac 12|\nabla \bar u_2^k|^2(\xi\cdot \nu))dS$$
	and setting $\Omega_k=B_{\tau \bar \varepsilon_k^{-1}}$, then the Pohozaev identity takes the form (recall that $\bar v_2^k$ is defined in (\ref{bar-v-2-q}))
	\begin{equation}
		\label{poho-q}
		\begin{aligned}
			&\int_{\Omega_k}\partial_{\xi}\bar h_{k,0}(\bar q_2^k+\bar \varepsilon_ky)e^{\bar v_2^k} \\
			&=\frac 1{\bar \varepsilon_k}\int\limits_{\partial \Omega_k}
			\bigg (\bar h_{k,0}e^{\bar v_2^k}(\xi\cdot \nu)+
			\partial_{\nu}\bar v_2^k\partial_{\xi}\bar v_2^k-\frac 12|\nabla \bar v_2^k|^2(\xi\cdot \nu)\bigg )dS. 
		\end{aligned}
	\end{equation}
	We use (\ref{v2-exp-new}) for the expansion of $\bar v_2^k$ and (\ref{grad-v2-new-2}) for $\nabla \bar v_2^k$ on $\partial \Omega_k$. In addition we use the following identities,
	$$\bar h_{k,0}(\bar q_2^k)\int_{\Omega_k}e^{U_2^k}dx=8\pi\left(1-\frac{1}{\bar h_{k,0}(\bar q_2^k)}\tau^{-2}\bar \varepsilon_k^2+O(\bar \varepsilon_k^4)\right),$$
	and
	\begin{align*}
		&\frac 14\frac{\partial_{\xi}(\Delta \log\bar h_{k,0})(\bar q_2^k)}{\bar h_{k,0}(\bar q_2^k)}\bar \varepsilon_k^2\int_{\Omega_k}\bar h_{k,0}(\bar q_2^k)|y|^2e^{U_2^k}dy\\
		&=\frac 14\bar \varepsilon_k^2\frac{\partial_{\xi}(\Delta\log \bar h_{k,0})(\bar q_2^k)}{\bar h_{k,0}(\bar q_2^k)}2\pi\int_0^{\tau \bar \varepsilon_k^{-1}}
		\frac{\bar h_{k,0}(\bar q_2^k)r^3}{(1+\frac{\bar h_{k,0}(\bar q_2^k)}8r^2)^2}dr\\
		&=16\pi\bar \varepsilon_k^2\frac{\partial_{\xi}(\Delta \log\bar h_{k,0})(\bar q_2^k)}{(\bar h_{k,0}(\bar q_2^k))^2}\left(\log \left(\frac{\bar h_{k,0}(\bar q_2^k)}8\tau^2\right)+
		2\log \frac{1}{\bar \varepsilon_k}\right)+O(\bar \varepsilon_k^4). 
	\end{align*}
	Based on (\ref{v2-exp-new2}), \eqref{grad-v2-new-2}, (\ref{v2-exp-new}) the left hand side of (\ref{poho-q}) takes the form 
	\begin{align*}
		&\partial_{\xi}\log \bar h_{k,0}(\bar q_2^k)\left(8\pi-\frac{8}{\bar h_{k,0}(\bar q_2^k)}\tau^{-2}\bar \varepsilon_k^2\right)\\
		&+16\pi\bar \varepsilon_k^2\frac{\partial_{\xi}\bar h_{k,0}(\bar q_2^k)}{\bar h_{k,0}^2(\bar q_2^k)}\Delta \log \bar h_{k,0}(\bar q_2^k)\bigg (\log \left(\frac{\bar h_{k,0}(\bar q_2^k)}8\tau^2\right)+2\log \frac{1}{\bar \varepsilon_k}\bigg )
		+O(\bar \varepsilon_k^2),
	\end{align*}
	where $O(\bar\varepsilon_k^2)$ comes from the $O(r^3)$ term in the expansion of $\partial_{\xi}\bar h_{k,0}$. 
	
	Concerning the right hand side of (\ref{poho-q}), we estimate the first term 
	by using the expansion in (\ref{for-v2-q}) and the cancellations of those integrals involving separable terms, to obtain that it is of order $O(\bar\varepsilon_k^2)$. Here we also use the fact that  $e^{v_2^k}=O(\bar\varepsilon_k^4)$ on $\partial \Omega_k$. For the two other terms, which involve $\nabla v_2^k$, we use (\ref{grad-v2-new-2}) to deduce that their contribution is of order $O(\bar\varepsilon_k^2)$. Therefore, from the evaluation of (\ref{poho-q}) around $\bar q_2^k$, we have that, 
	\begin{equation}
		\label{new-sigma-1}
		\begin{aligned}
			&16\pi\bar \varepsilon_k^2\frac{\partial_{\xi}\bar h_{k,0}(\bar q_2^k)}{\bar h_{k,0}^2(\bar q_2^k)}\Delta \log \bar h_{k,0}(\bar q_2^k)\bigg (\log \left(\frac{\bar h_{k,0}(\bar q_2^k)}8\tau^2\right)+2\log \frac{1}{\bar \varepsilon_k}\bigg )\\
			&+8\pi\partial_{\xi}\log \bar h_{k,0}(\bar q_2^k)=O(\bar \varepsilon_k^2), 
		\end{aligned}
	\end{equation}
	where we have used the fact that $|\nabla \log \bar h_{k,0}(\bar q_2^k)|=O(\bar\varepsilon_k^2\log\frac{1}{\bar\varepsilon_k})$. The same argument
	yields a similar estimate for
	$$
	\bar v_1^k(y)=\bar u_1^k(\bar q_1^k+\bar \varepsilon_{1,k}y)+2\log \bar \varepsilon_{1,k},\quad  |y|<\tau\bar \varepsilon_{1,k}^{-1},\quad
	\bar \varepsilon_{1,k}=e^{-\frac{\bar \lambda^k_{1}}{2}}.
	$$
	Indeed, denoting by $\bar h_{k,1}=
	\bar h_{k,0}e^{\phi^k_1-\phi_2^k}$, which is the weight function relative to $\bar u_1^k$, we have
	\begin{equation}
		\label{new-sigma-2}
		\begin{aligned}
			&16\pi\bar \varepsilon_{1,k}^2\frac{\partial_{\xi}\bar h_{k,1}(\bar q_1^k)}{\bar h_{k,1}^2(\bar q_1^k)}\Delta \log \bar h_{k,1}(\bar q_1^k)\bigg (\log \left(\frac{\bar h_{k,1}(\bar q_1^k)}8\tau^2\right)+2\log \frac{1}{\bar \varepsilon_{1,k}}\bigg )\\
			&+8\pi\partial_{\xi}\log \bar h_{k,1}(\bar q_1^k)=O(\bar \varepsilon_{1,k}^2). 
		\end{aligned}
	\end{equation}
	In view of $\rho_1^k=\rho_2^k=\rho_k$ we already know (see \eqref{initial-small}) that $\bar \lambda_1^k-\bar\lambda_2^k=O(e^{-\epsilon_0\lambda_1^k})$ ($i=1,2$) for some $\epsilon_0>0$, whence we also have 
	$$\bar\varepsilon_{1,k}-\bar \varepsilon_k=O(\bar \varepsilon_k^{1+\epsilon_0}).$$
	The difference between $\bar h_{k,1}$ and $\bar h_{k,0}$ is (see \eqref{xi-v})
	$$\bar h_{k,1}(x)-\bar h_{k,0}(x)=\bar h_{k,0}(x)(e^{\sigma_k\bar \psi_{\xi}^k}-1).$$
	By Lemma \ref{small-osci-xi} we have,
	$$|\bar \psi^k_{\xi}(\bar \varepsilon_ky)|\leq C\bar \varepsilon_k \varepsilon_ky, $$
	whence, for a small positive constant $\epsilon_0>0$,
	$$|\sigma_k\bar \psi^k_{\xi}(\bar \varepsilon_ky)|\leq C\bar \varepsilon_k^{2+\epsilon_0}|y|,\quad |y|<\frac{\tau}{\bar \varepsilon_k}.$$
	Thus, recalling that $\bar \psi^k_{\xi}$ is harmonic, we see that all the terms containing 
	$\bar \psi^k_{\xi}$  are negligible and we can expand the first term in (\ref{new-sigma-2}) as follows, 
	\begin{align*}
		\partial_{\xi}\log \bar h_{k,1}(\bar q_1^k)=~&
		\partial_{\xi}\log \bar h_{k,0}(\bar q_2^k)+
		\sum_{j=1}^2\partial_{\xi j}\log \bar h_{k,0}(\bar q_2^k)(\bar q_1^k-\bar q_2^k)_j\\
		&+O(\bar\varepsilon_k^2)+
		O(|\bar q_1^k-\bar q_2^k|^2).
	\end{align*}
	The same argument applies of course to the other terms on the left of (\ref{new-sigma-2}). Therefore, for $\xi=1,2$, the difference between (\ref{new-sigma-1}) and (\ref{new-sigma-2}) takes the form,
	$$\left(\begin{matrix}
		\partial_{11}\log \bar h_{k,0}(\bar q_2^k)+E_{h,11} & \partial_{12}\log \bar h_{k,0}(\bar q_2^k)+E_{h,12} \\
		\partial_{21}\log \bar h_{k,0}(\bar q_2^k)+E_{h,21} & \partial_{22}\log \bar h_{k,0}(\bar q_2^k)+E_{h,22}
	\end{matrix}\right)
	\left(\begin{matrix}
		(\bar q_1^k-\bar q_2^k)_1\\
		(\bar q_1^k-\bar q_2^k)_2
	\end{matrix}
	\right)=O(\bar\varepsilon_k^2).$$
	where $E_{h,ij}=O\left(\bar\varepsilon_k^2\log \frac{1}{\bar \varepsilon_k}\right)$. 	At this point the non-degeneracy assumption about $D^2\log(\bar h_{k,0})$ readily implies that $|q_1^k-q_2^k|=O(\bar \varepsilon_k^2)$ as well.
	Lemma \ref{crucial-q12} is established. \end{proof}
	
	\begin{rem} The non-degeneracy assumption of ${\rm det}(D^2f^*)\neq 0$ is equivalent to ${\rm det}(D^2\log(\bar h_{k,0}))(\bar q_2^k)\neq 0$.
	\end{rem}
	
	\medskip

	As a consequence of this improved estimate about $|\bar q_1^k-\bar q_2^k|$ we have,
	
	\begin{cor}\label{cor-sigma}  There exists $C>0$ independent of $k$ such that $$\sigma_k\le C\bar \varepsilon_k.$$
	\end{cor}
	\begin{proof}
We take advantage of the new estimate of $|\bar q_1^k-\bar q_2^k|$ and go through the argument used in Lemma \ref{new-sigma} to estimate $\sigma_k$, by using the expansion (\ref{v2-exp-new}) and Lemma \ref{crucial-q12}. Then we see that $b_1,b_2$ are determined by
	$$\frac{h_1}2\frac{\frac{\bar q_1^k- \bar q_2^k}{\bar \varepsilon_k\sigma_k}\cdot y}{(1+\frac{h_1}8|y|^2)},$$
	which implies, in view of the new estimate about $|\bar q_1^k-\bar q_2^k|$, that we have
	$$\sigma_k=O(\bar \varepsilon_k).$$ 
	We omit the details of this proof just to avoid repetitions.
\end{proof}
	
	\medskip
	
	\noindent{\bf Step five: Improved estimate on the oscillation of $\xi_k$ away from blowup points and vanishing rates of $b_0^k$, $\hat b_0^k$ and $\bar b_0^k$.}
	
	\medskip
	
	An major consequence of Corollary \ref{cor-sigma} is that the oscillation of $\xi_k$ far away from blow up points is of order $O(\bar \varepsilon_k)$.
	
	\begin{prop}\label{small-osc-psi}
		For any $x_1,x_2\in M\setminus \{\Omega(p_1,\tau)\cup \Omega(p_2,\tau)\cup\Omega(q,\tau)\}$,
		$$|\xi_k(x_1)-\xi_k(x_2)|\le C\bar \varepsilon_k$$
		for some $C>0$ independent by $k$.
	\end{prop}
	
	\begin{proof}  Recall that Lemma \ref{small-osci-xi} speaks that the oscillation of $\xi_k$ far away from blow up points is of order $O(\varepsilon_k)$.
	Because of this fact, the estimates about the harmonic functions used to eliminate the oscillations
	of $\xi_k$ is improved as well. For example concerning the  function $\psi_{\xi}^k$ around $p_1$ we have,
	$$\psi_{\xi}^k(\varepsilon_ky)=\sum_ja_j\varepsilon_k^2y_j+O(\varepsilon_k^3)|y|^2,\quad |y|<\tau \varepsilon_k^{-1},$$
	for suitable $a_j,j=1,2$.
	We study the expansion of $\xi_k$ around $p_1$ again to obtain an estimate about the decay of $b_0^k$. In fact at this point we know that $b_0^k\to 0$.
	In local coordinates around $p_1$, $0=x(p_1)$, the first term of the approximation of $\xi_k$ around $p_1$ is still $w_{0,\xi}^k$ defined as in (\ref{w0kxik}) that satisfies (\ref{e-w0kxi}).
	We need to write the equation for $w_{0,\xi}^k$ in the same form as that satisfied by
	$\tilde \xi_k(y)$ (see \eqref{around-p_1}). Thus by using the expansion of $v_2^k$ around $p_1$ (\ref{p-rough-2}), we write
	\begin{align*}
		&|y|^{2\alpha} \tilde h_k(0)\mbox{exp} \bigg
		(v_2^k+\log\frac{\tilde h_k(\varepsilon_ky)}{\tilde h_k(0)}\bigg )\\
		&=|y|^{2\alpha}\tilde h_k(0)\mbox{exp}\{U_2^k+c_{0,k}+\tilde c_{1,k}+\tilde c_{2,k} \\
		&\quad+\frac 14\Delta (\log \tilde h_k)(0)|y|^2\varepsilon_k^2+O(\varepsilon_k^{2+\epsilon_0})(1+|y|)^{2+\epsilon_0})\} \\
		&=|y|^{2\alpha} \tilde h_k(0)e^{U_2^k}(1+c_{0,k}+\tilde c_{1,k}+\tilde c_{2,k}
		+\frac 14\Delta (\log \tilde h_k)(0)|y|^2\varepsilon_k^2 \\
		&\quad+\frac 12(\tilde c_{1,k})^2)+O(\varepsilon_k^{2+\epsilon_0})(1+|y|)^{-2-\epsilon_0}, 
	\end{align*}
	for some small $\epsilon_0>0$ depending by $\alpha>0$. Then we can write the equation of $w_{0,\xi}^k$ as follows,
	\begin{equation}
		\label{e-hat-xi-new}
		\begin{aligned}
			&\Delta w_{0,\xi}^k+\tilde h_k(0)|y|^{2\alpha}\exp\left (v_2^k+\log \frac{\tilde h_k(\varepsilon_ky)}{\tilde h_k(0)}\right )w_{0,\xi}^k\\
			&=|y|^{2\alpha} \tilde h_k(0)e^{U_2^k}w_{0,\xi}^k
			\left(c_{0,k}+\tilde c_{1,k}+\tilde c_{2,k}+\frac{\Delta (\log  \tilde h_k)(0)|y|^2\varepsilon_k^2}{4}+\frac{\tilde c_{1,k}^2}{2}\right) \\
			&\quad+O(\varepsilon_k^{2+\epsilon_0})b_0^k(1+|y|)^{-2-\epsilon_0}. 
		\end{aligned}
	\end{equation}
	The next term in the expansion is $w_{1,\xi}^k$ defined as in (\ref{xi-2nd}). Now we need to include $b_0^k$ in the estimate of $w_{1,\xi}^k$,
	$$|w_{1,\xi}^k(y)|\le Cb_0^k \varepsilon_k (1+|y|)^{-1}, $$
	where we used $w_{0,\xi}^k(0)=b_0^k$.
	At this point we write the equation for $w_{1,\xi}^k$ in the following form:\begin{align*}
		&\Delta w_{1,\xi}^k+\tilde h_k(0)|y|^{2\alpha}\exp \bigg (v_2^k+\log \frac{\tilde h_k(\varepsilon_ky)}{\tilde h_k(0)}\bigg ) w_{1,\xi}^k \\
		&=- \tilde h_k(0)|y|^{2\alpha}e^{U_2^k}\tilde c_{1,k}w_{0,\xi}^k+|y|^{2\alpha}\tilde h_k(0)e^{U_2^k}\tilde c_{1,k}w_{1,\xi}^k
		+\frac{O(b_0^k\varepsilon_k^{3})}{(1+|y|)^{3+2\alpha}}.\nonumber
	\end{align*}
	We first consider just a rough estimate of the last two terms,
	\begin{equation}
		\label{crude-w1}
		\begin{aligned}
			&\Delta w_{1,\xi}^k+\tilde h_k(0)|y|^{2\alpha}\exp\left({v_2^k+\log \frac{\tilde h_k(\varepsilon_ky)}{\tilde h_k(0)}}\right)w_{1,\xi}^k \\
			&=- \tilde h_k(0)|y|^{2\alpha}e^{U_2^k}\tilde c_{1,k}w_{0,\xi}^k+\frac{O(b_0^k\varepsilon_k^2)}{(1+|y|)^{3+2\alpha}}. 
		\end{aligned}
	\end{equation}
	Let
	$$\underline{w}_1^k=\tilde \xi_k(y)-w_{0,\xi}^k-w_{1,\xi}^k-\psi_{\xi}^k(\varepsilon_ky),$$
	then from (\ref{around-p_1}), (\ref{e-hat-xi-new}) and (\ref{crude-w1}) we see that the equation for $\underline{w}_1^k$ takes the form,
	\begin{equation}
		\label{f-W1k}
		\begin{aligned}
			&\left(\Delta+|y|^{2\alpha}\tilde h_k(0)\exp\left(v_2^k+\log \frac{\tilde h_k(\varepsilon_ky)}{\tilde h_k(0)}\right)\right)\underline{w}_1^k \\
			&=\frac{O(\bar \varepsilon_k)}{(1+|y|)^{4+2\alpha}}+\frac{O(b_0^k\varepsilon_k^2)}{(1+|y|)^{3+2\alpha}}+\frac{O(\varepsilon_k^2)}{(1+|y|)^{3+2\alpha}}. 
		\end{aligned}
	\end{equation}
	Concerning the r.h.s of \eqref{f-W1k}, the first term comes from the expansion of $c_k$, the second term comes from $w_{1,\xi}^k$,
	the third term comes from $\psi_{\xi}^k$. Here we also note that $w_{0,\xi}^k$ is constant on $\partial \Omega_k$,
	while $w_{1,\xi}^k$ has an oscillation of order $O(b_0^k\varepsilon_k^2)$. To eliminate this oscillation we use another harmonic function of order
	$O(b_0^k\varepsilon_k^3|y|)$. For simplicity and with an abuse of notations we include this function in $\psi_{\xi}^k$.
	As a consequence standard potential estimates show that,
	\begin{equation*}
		|\underline{w}_1^k(y)|\leq C(\delta)(\varepsilon_k^2+\bar \varepsilon_k)(1+|y|)^{\delta},
		\quad y\in \Omega_k.
	\end{equation*}
	Therefore the improved estimate of $\psi_{\xi}^k$ yields an improved estimate about $\underline{w}_1^k$ which then yields to an improvement of the estimate about $\tilde\xi_k$. In other words we run a refined bootstrap argument.
	\begin{rem}
		%\label{rem:al1}
		We only need to consider $\alpha>1$, because if $0<\alpha<1$ then all these estimates are not needed.
		We will be back to this point later on.
	\end{rem}
	With this new estimate about $\underline{w}_1^k$ we can just go through the evaluation of the oscillation of $\xi_k$ far away from the blow up
	points as in Lemma \ref{small-osci-xi} to deduce that,
	\begin{equation}\label{new-osci}
		|\xi_k(x_1)-\xi_k(x_2)|\le C(\varepsilon_k^2+\bar \varepsilon_k).
	\end{equation}
	
	The new estimate (\ref{new-osci}) implies that the harmonic function $\psi_{\xi}^k$ that encodes the oscillation of
	$\xi_k$ on $\partial \Omega(p_1,\tau)$ satisfies
	\begin{equation}\label{new-psi-osc}
		|\psi_{\xi}^k(\varepsilon_ky)|\le C(\varepsilon_k^3+\varepsilon_k\bar \varepsilon_k)|y|,\quad |y|<\tau \varepsilon_k^{-1}.
	\end{equation}
	As a consequence of (\ref{new-psi-osc}) we can improve the estimate about $\underline{w}_1^k$ as follows, 
	\begin{equation}\label{int-w1}
		|\underline{w}_1^k(y)|\le C(\delta)(b_0^k\varepsilon_k^2+\varepsilon_k^3+\bar \varepsilon_k)(1+|y|)^{\delta}.
	\end{equation}
	In order to further improve the estimate about $\psi_{\xi}^k$ we point out that, in view of (\ref{new-osci}), then around $p_2$ and $q$ we have the following analogous improved estimates. 
	First of all, around $p_2$, we recall that $\hat w_k$ is defined in (\ref{hat-wk}) and satisfies (\ref{sec-e}). Then we have, 
	$$
	\int_{B(p_2,\tau)}\rho_kHc_k\xi_k=O(\bar \varepsilon_k).
	$$
	Moreover, around $q_2^k$, we recall that $\bar w_k$ is defined in (\ref{bar-w2k}) and satisfies (\ref{pre-xi}). By using (\ref{pre-xi}) and the expansion of $\bar w_k$ we have
	\begin{equation*}
		\int_{B(q_2^k,\tau)}\rho_k Hc_k\xi_k=O(\bar \varepsilon_k).
	\end{equation*}
	At this point we move back to the evaluation of the oscillation of $\xi_k$ away from blow up points. From the Green's representation formula of $\xi_k$, as in the proof of (\ref{osi-vars}), by using (\ref{new-psi-osc}) for $\psi_{\xi}^k$, (\ref{int-w1}) for $\underline{w}_1^k$, (\ref{sec-e}) for $\hat w_k$ and (\ref{pre-xi}) for $\bar w_k$ we can further improve the estimate of $\xi_k$ as follows,
	\begin{equation}\label{more-xi-o}
		|\xi_k(x)-\xi_k(y)|\leq C(b_0^k\varepsilon_k^2+\varepsilon_k^4+\bar \varepsilon_k),
	\end{equation}
	for $x,y\in M\setminus\{\Omega(p_1,\tau)\cup \Omega(p_2,\tau)\cup \Omega(q,\tau)\}$. 
	Then (\ref{more-xi-o}) further improves the estimate of $\underline{w}_1^k$, so that, by a bootstrap argument which does not involve the leading term proportional to $b_0^k\varepsilon_k^2$ (which in fact comes just from the integration over $\Omega(p_1,\tau)$), after a finite number of iterations, we deduce that,
	\begin{equation*}
		|\underline{w}_1^k(y)|\le C(\delta)(b_0^k\varepsilon_k^2+\bar \varepsilon_k)(1+|y|)^{\delta},
	\end{equation*}
	and the estimate about the oscillation of $\xi_k$ takes the form,
	\begin{equation}\label{best-osci-2}
		|\xi_k(x_1)-\xi_k(x_2)|\le C(b_0^k\varepsilon_k^2+\bar \varepsilon_k),
	\end{equation}
	for all $x_1,x_2\in M\setminus \{ \Omega(p_1,\tau),\cup  \Omega(p_2,\tau)\cup
	\Omega(q,\tau)\}$. 
	
	As an immediately consequence of (\ref{best-osci-2}) for $\psi_{\xi}^k$ we have
	\begin{equation}\label{osc-p-1}
		|\psi_{\xi}^k(\varepsilon_ky)|\leq C(b_0^k\varepsilon_k^3+\varepsilon_k\bar \varepsilon_k)|y|.
	\end{equation}
	
	At this point we prove a more precise estimate about $b_0^k$.

	\begin{lem}\label{b0-va}
		There exists $C>0$ independent of $k$ such that 
		\begin{equation}\label{small-d}
			|b_0^k|\le C e^{-\frac{\lambda_1^k}2+\frac{\lambda_1^k}{1+\alpha}}\quad \mbox{if}\quad \alpha>1.
		\end{equation}
	\end{lem}
	\begin{proof}[Proof of Lemma \ref{b0-va}.] Recall that $w_{2,\xi}^k$ was defined (see \eqref{first:estimate}) to encode the $\Theta_2$ and other separable terms of the order $O(\varepsilon_k^2)$.
	For the estimate of this term we have (\ref{1st-w2k}), but now we include the dependence of $b_0^k$. So (\ref{1st-w2k}) takes the form, 
	\begin{equation*}
		|w_{2,\xi}^k(y)|\le Cb_0^k\varepsilon_k^2(1+|y|)^{-2\alpha}.
	\end{equation*}
	We need at this point another term in the expansion suitable to handle the radial part of order $O(\varepsilon_k^2)b_0^k$, which is why
	we set $z_{0,k}$ to be the radial function satisfying
	\begin{align*}
		\begin{cases}
			\frac{d^2}{dr^2}z_{0,k}+\frac 1r\frac{d}{dr}z_{0,k}+r^{2\alpha} \tilde h_k(0)e^{U_2^k}z_{0,k}=F_{0,k}, \quad 0<r<\tau \varepsilon_k^{-1},\\
			\\
			z_{0,k}(0)=\frac{d}{dr}z_{0,k}(0)=0
		\end{cases}
	\end{align*} where
	\begin{align*}
		F_{0,k}(r)
		=~&- \tilde h_k(0)|y|^{2\alpha}e^{U_2^k}\bigg (w_{0,\xi}^k\big (c_{0,k}+\frac 14\Delta (\log  \tilde h_k)(0)|y|^2\varepsilon_k^2\\
		&+\frac  14\varepsilon_k^2|\nabla \log  \tilde h_k(0)|^2(g_k+r)^2\big )
		+(w_{\xi,1}^k\tilde c_{1,k})_r\bigg )\nonumber
	\end{align*}
	where we use the fact that the projection of $\frac 12(\tilde c_{1,k})^2$ along the constant direction is
	$$\frac 14\varepsilon_k^2|\nabla  \log \tilde h_k(0)|^2(g_k+r)^2.$$ Observe carefully that we already considered this function,
	see \eqref{eq-z0} in subsection \ref{difference} above, but in that case the underlying assumption was $b_0^k\to b_0\neq 0$. Concerning $z_{0,k}$ we have,
	$$|z_{0,k}(r)|\le C(\delta)b_0^k \varepsilon_k^2(1+r)^{\delta},$$
	and then, by defining $\underline{z}_k$ as follows,
	$$\underline{z}_k(y)=\xi_k(\varepsilon_ky)-\psi_{\xi}^k(\varepsilon_k y)-w_{0,\xi}^k(y)-w_{1,\xi}^k(y)-w_{2,\xi}^k(y)-z_{0,k}(y),$$
	we have
	$$\Delta \underline{z}_k+|y|^{2\alpha}h_k(\varepsilon_ky)\varepsilon_k^{2+2\alpha}\tilde c_k(\varepsilon_ky)\underline{z}_k=E_k$$
	with $z_k(0)=0$ and
	$$|E_k|\le Cb_0^k(\varepsilon_k^{2+\epsilon_0}+\bar \varepsilon_k)(1+|y|)^{-2-\frac{\epsilon_0}2},$$
	for some $\epsilon_0>0$.
	Recall that at this moment $\psi_{\xi}^k$ satisfies (\ref{osc-p-1}). By standard potential estimates we have,
	\begin{equation*}
		|z_k(y)|\le C(\delta)(b_0^k\varepsilon_k^{2+\epsilon_0}+\bar \varepsilon_k)(1+|y|)^{\delta}.
	\end{equation*}
	At this point, by using the estimate about $\underline{z}_k$, we deduce that,
	\begin{align*}
		\int_{\Omega(p_1,\tau)}\rho_kHc_k\xi_k
		&=\int_{\Omega_k}\tilde h_k(\varepsilon_ky)|y|^{2\alpha}e^{v_2^k}(\xi_k(\varepsilon_ky)-\psi_{\xi}^k(\varepsilon_ky))dy\\
		&\quad+\int_{\Omega_k}\tilde h_k(\varepsilon_ky)|y|^{2\alpha}e^{U_2^k}\psi_{\xi}^k(\varepsilon_ky)\\
		&=d_{1,k}b_0^k(\Delta \log h_k(0)+O(\varepsilon_k^{\epsilon_0}))\varepsilon_k^2+O(\bar \varepsilon_k),
	\end{align*}
	where the derivation of the above is similar to the derivation of (\ref{p1-leading}) except that we keep track of $b_0^k$ in this estimate. 
	Therefore the integral around $p_1$ reads,
	$$\int_{\Omega(p_1,\tau)}\rho_kH c_k\xi_k=c b_0^k\varepsilon_k^2(\Delta (\log \tilde h_k)(0)+
	O(\varepsilon_k^{\epsilon_0}))+O(\bar \varepsilon_k),\quad c\neq 0.$$
	From $\int_M\rho_kHc_k\xi_k=0$ we readily obtain (\ref{small-d}) by splitting the domain of integration in the three regions near the blow up points
	plus the corresponding complement in $M$. Lemma \ref{b0-va} is established. \end{proof}
	
	\medskip
	
	In view of \eqref{best-osci-2}, by \eqref{small-d} we immediately deduce that Proposition \ref{small-osc-psi} holds. \end{proof}
	
	\medskip
	
	With Proposition \ref{small-osc-psi} we give an estimate of the first term in the approximation of $\xi_k$ around each blowup point. We recall that the first term in the approximation of $\xi_k$ around $p_1$ is proportional to $b_0^k$, around $p_2$ is proportional to $\hat b_0^k$, around $q_2^k$ is proportional to $\bar b_0^k$. For these quantities we have
	\begin{lem}\label{same-b} There exists $C>0$ independent by $k$, such that
		$$|\bar b_0^k|+|\hat b_0^k|\le C\bar \varepsilon_k\varepsilon_k^{-2}=Ce^{-\frac{\lambda_1^k}2+\frac{\lambda_1^k}{1+\alpha}}.$$
	\end{lem}
	
	\begin{proof}[Proof of Lemma \ref{same-b}.] In the proof of Proposition \ref{small-osc-psi} we already established 
	\eqref{small-d} for $b_0^k$. Also the oscillation of $\xi_k$ far away from blow up points is of order $O(\bar \varepsilon_k)$ (Proposition \ref{small-osc-psi}). In the expansion of $\xi_k$ around $q_2^k$ we have 
	$$\xi_k(x)=-\bar b_0^k+O(\bar \varepsilon_k),\quad x\in \partial B(q_2^k,\tau).$$
	On $\partial B(p_1^k,\tau)$ we have
	$$\xi_k=-b_0^k+O(b_0^k\varepsilon_k^2)=-b_0^k+O(\bar \varepsilon_k),$$
	and similarly on $\partial B(p_2,\tau)$, 
	$$\xi_k(x)=\hat b_0^k+O(\bar \varepsilon_k).$$
	Thus $|b_0^k-\hat b_0^k|+|b_0^k-\bar b_0^k|\le C\bar \varepsilon_k.$
	Lemma \ref{same-b} is established. 
\end{proof}

	\medskip
	
	\noindent{\bf Step six: an improved expansion for $\bar \xi_k$.}
	
	\medskip
	
	As a consequence of Proposition \ref{small-osc-psi} %and Lemma \ref{same-b} 
	we now write the equation of $\bar \xi_k$ (see \eqref{xi-k-q}) as follows,
	$$\Delta \bar \xi_k+\bar h_{k,0}(\bar q_2^k+\bar \varepsilon_ky)e^{\bar v_2^k}(1+\sigma_k\bar \xi_k)\bar \xi_k=O(\bar \varepsilon_k^{3/2})(1+|y|)^{-5/2}.$$
	By using $\nabla \log \bar h_{k,0}(\bar q_2^k)=O(\bar \lambda_1^ke^{-\bar\lambda_1^k})$ and the expansion of $\bar v_2^k$ in (\ref{v2-exp-new2})
	we have
	\begin{equation*}
		\Delta \bar \xi_k(y)+\bar h_{k,0}(\bar q_2^k)e^{U_2^k}(1+\sigma_k\bar \xi_k)\bar \xi_k=O(\bar \varepsilon_k^{3/2})(1+|y|)^{-5/2}.
	\end{equation*}
	At this point it is enough to use $\sigma_k=O(\bar \varepsilon_k)$ so that by standard potential estimates, we have,
	\begin{equation}\label{q-rough}
		|\bar \xi_k(y)-\bar w_{0,\xi}^k|\leq C(\delta)\bar \varepsilon_k (1+|y|)^{\delta}.
	\end{equation}
	There is no need to subtract $\bar \psi_{\xi}^k(\bar \varepsilon_ky)$ here because of Proposition \ref{small-osc-psi}
	$$|\bar \psi_{\xi}^k(\bar \varepsilon_ky)|\le C\bar \varepsilon_k^2|y|.$$
	In view of (\ref{q-rough}) we write the equation of $\bar \xi_k$
	as follows,
	$$\Delta \bar \xi_k+\bar h_{k,0}(\bar q_2^k)e^{U_2^k}\bar \xi_k=-\sigma_k\bar h_{k,0}(\bar q_2^k)e^{U_2^k}(\bar w_{0,\xi}^k)^2+O(\bar \varepsilon_k^{3/2})(1+|y|)^{-5/2},$$
	which immediately leads to 
	\begin{equation*}
		(\Delta +\bar h_{k,0}(\bar q_2^k)e^{U_2^k})(\bar \xi_k-\bar w_{0,\xi}^k)
		=-\sigma_k\bar h_{k,0}(\bar q_2^k)e^{U_2^k}(\bar w_{0,\xi}^k)^2+O(\bar \varepsilon_k^{3/2})(1+|y|)^{-5/2},
	\end{equation*}
	in $\Omega_k$. 
	We now construct the next term in the expansion as follows,
	\begin{equation*}
		\Delta \bar w_{1,\xi}^k+\bar h_{k,0}(\bar q_2^k)e^{U_2^k}\bar w_{1,\xi}^k=-\sigma_k\bar h_{k,0}(\bar q_2^k)e^{U_2^k}(\bar w_{0,\xi}^k)^2_{\theta,2},
	\end{equation*}
	where $(\bar w_{0,\xi}^k)^2_{\theta,2}$ means the projection of $(\bar w_{0,\xi}^k)^2$ onto $e^{2i\theta}$. Note that, because of the vanishing rate of $\bar b_0^k$, we can safely disregard the projection of $(\bar w_{0,\xi}^k)^2$ along  $e^{i\theta}$, just because it is proportional to $\bar b_0^k$. 
	We use a separable function to handle this term. The radial part of $\bar w_{1,\xi}^k$ is represented by $\bar f_1^k$, which satisfies
	for $0<r<\tau\bar \varepsilon_k^{-1}$
	\begin{equation*}
		\frac{d^2}{dr^2}\bar f_1^k(r)+\frac 1r\frac{d}{dr}\bar f_1^k(r)+(\bar h_{k,0}(\bar q_2^k)e^{U_2^k}-\frac{4}{r^2})\bar f_1^k(r)=\bar E_k
	\end{equation*}
	where 
	$$\bar E_k=\sigma_k\bar h_{k,0}(\bar q_2^k)\frac{r^2}{(1+\frac{\bar h_{k,0}(\bar q_2^k)}8r^2)^4}.$$
	Let $\bar f_{1a},\bar f_{1b}$ be a pair of fundamental solutions which behaves as follows: 
	$|\bar f_{1a}(r)|\sim r^2$ near $0$ and $\infty$,  
	$     |\bar f_{1b}(r)|\sim r^{-2}$  
	near $0$ and $\infty$, 
	we construct $\bar f_1^k$ as follows,
	$$\bar f_1^k=\bar f_{1a}(r)\int_r^{L_k}\frac{W_1}{W}ds+\bar f_{1b}(r)\int_0^r\frac{W_2}{W}ds,\quad L_k=\tau \bar \varepsilon_k^{-1},$$
	where one can easily verify that
	$$W(r)\sim -\frac{1}{r},\quad |W_1(r)|\le C\sigma_k(1+r)^{-8},\quad |W_2(r)|\le C\sigma_k r^4(1+r)^{-8}.$$
	Then it is readily seen that a first estimate for $\bar f_1^k(r)$ reads, 
	$$|\bar f_1^k(r)|\le C\bar \varepsilon_kr(1+r)^{-5},\quad 0<r<\tau \bar \varepsilon_k^{-1}.$$
	Here we remark that even though $\frac{d}{dr}\bar f_1^k(0)\neq 0$, it does not interfere with $\nabla \bar w_{1,\xi}^k(0)$ because $\bar f_1^k$ is multiplied by $e^{2i\theta}$ in $\bar w_{1,\xi}^k$. Thus 
	\begin{equation}\label{decay-w1}
		|\bar w_{1,\xi}^k(y)|=O(\bar \varepsilon_k)(1+|y|)^{-4},
	\end{equation}
	where $O(\bar\varepsilon_k)$ comes from the estimate of $\sigma_k$. Next we use a radial term to eliminate the radial part of
	$\sigma_k\bar h_{k,0}(\bar q_2^k)e^{U_2^k}(\bar w_{0,\xi})^2_r$ for $r\in(0,\tau\bar\varepsilon_k^{-1})$:
	\begin{equation*}
		%\label{c0xi}
		\begin{cases}
			\frac{d^2}{dr^2}c_{0,\xi}^k+\frac 1r\frac{d}{dr}c_{0,\xi}^k+\bar h_{k,0}(\bar q_2^k)e^{U_2^k}c_{0,\xi}^k=\sigma_k\bar h_{k,0}(\bar q_2^k)e^{U_2^k}(\bar w_{0,\xi})^2_r,\\
			\\
			c_{0,\xi}^k(0)=\frac{d}{dr}c_{0,\xi}^k(0)=0.
		\end{cases}
	\end{equation*}
	Then it is not difficult to check that $c_{0,\xi}^k$ satisfies
	\begin{equation}\label{estc0xi}
		|c_{0,\xi}^k(r)|\le C\bar \varepsilon_k\log (1+r),
		\quad |\frac{d}{dr}c_{0,\xi}^k(r)|\le C \bar \varepsilon_k(1+r)^{-1}.
	\end{equation}
	Let 
	$$z_{0,\xi}^k=\bar \xi_k-\bar w_{0,\xi}^k-\bar w_{1,\xi}^k-c_{0,\xi}^k,$$ 
	we have 
	\begin{equation}\label{e-z0xi}
		(\Delta +\bar h_{k,0}(\bar q_2^k)e^{U_2^k})(z_{0,\xi}^k- \psi^k_{\xi}(\bar\varepsilon_ky))=O(\bar \varepsilon_k^{1+\epsilon_0})(1+|y|)^{-5/2}
	\end{equation}
	in $\Omega_k$ for some $\epsilon_0>0$. Note that the oscillation of $z_{0,\xi}^k$ is $O(\bar \varepsilon_k^2)$. If we use an harmonic function to eliminate it, well this harmonic function is so small that it can be safely neglected just by considering it as a part of $\bar\psi_{\xi}^k(\bar \varepsilon_ky)$. As a consequence in particular 
	we may assume without loss of generality that $z_{0,\xi}^k-\bar\psi_{\xi}^k(\bar \varepsilon_ky)$ has no oscillation on $\partial \Omega_k$. 
	Thus we deduce from \eqref{e-z0xi} that,
	\begin{equation}\label{bar-xi-e}
		|z_{0,\xi}^k-\bar \psi_{\xi}^k(\bar \varepsilon_ky)|\leq C(\delta)\bar \varepsilon_k^{1+\epsilon_0}(1+|y|)^{\delta},
		\quad |y|<\tau\bar \varepsilon_k^{-1},
	\end{equation}
	for some $\epsilon_0>0$. Then we apply standard elliptic estimates based on (\ref{e-z0xi}) and (\ref{bar-xi-e}) to obtain,
	\begin{equation}\label{small-z0xi}
		|\nabla (z_{0,\xi}^k-\bar \psi_{\xi}^k(\bar\varepsilon_k\cdot))|\le C\bar \varepsilon_k^{2+\epsilon_0}\quad \mbox{on}\quad \partial B_{\frac{\tau}2\bar\varepsilon_k^{-1}}
	\end{equation}
	for some $\varepsilon_0>0$. From (\ref{small-z0xi}) we see that for $r\sim \bar\varepsilon_k^{-1}$, the leading term in $\nabla \xi_k$ is $\bar w_{0,\xi}^k+c_{0,\xi}^k$: 
	\begin{align}\label{xi-k-large}
		\begin{cases}
			\partial_1\bar \xi_k=\frac{b_1^k(1+a_kr^2)-2a_ky_1(b_1^ky_1+b_2^ky_2)}{(1+a_kr^2)^2}+\frac{d}{dr}c_{0,\xi}^k(r)\frac{y_1}r+O(\bar \varepsilon_k^{2+\delta_0}),\\
			\\
			\partial_2\bar \xi_k=\frac{b_2^k(1+a_kr^2)-2a_ky_2(b_2^ky_2+b_1^ky_1)}{(1+a_kr^2)^2}+\frac{d}{dr}c_{0,\xi}^k(r)\frac{y_2}r+O(\bar \varepsilon_k^{2+\delta_0}) 
		\end{cases}
	\end{align}
	where $a_k=\frac{\bar h_{k,0}(\bar q_2^k)}8$ and $\delta_0$ is a small positive constant. 
	
	\medskip
	
	\noindent{\bf Step seven: Evaluation of a simplified Pohozaev identity}
	
	\medskip
	
	In a local coordinates system around $q$, $0=x(q)$, we recall the set up in local coordinates around $q_2^k$: $u_i^k={\rm w}_i^k-f_k$ where $f_k$ is defined in (\ref{def-fk}) and $u_k$ is a solution of (\ref{around-q}).  The term $e^{\phi}$ in (\ref{def-fk}) is just
	the local conformal factor of the metric (see \eqref{def-phi}). Let us recall that $\rm w_i^k$ satisfies
	$$\Delta_g {\rm w}_i^k+\rho_k(He^{\rm w_i^k}-1)=0,~\ \int_M He^{\rm w_i^k}=1,\quad i=1,2,$$
	where
	$$H=he^{-4\pi\alpha G(x,p_1)-4\pi \beta G(x,p_2)}.$$
	
	Since ${\rm w}_i^k=u_i^k+f_k$, we have $\xi_k=\sigma_k^{-1}({\rm w}_1^k-{\rm w}_2^k)$.
	Recall that, for the time being, we are working just with three blowup points, which is why we can choose $f_k$ as follows,
	\begin{equation}
		\label{exp-f-3}
		\begin{aligned}
			f_k(x)=~&\frac{\rho_k}{3+\alpha+\beta}\left(R(x,\bar q_{2}^k)-R(\bar q_{2}^k,\bar q_{2}^k)\right)\\
			&+\frac{\rho_k}{3+\alpha+\beta}\left( (1+\alpha)(G(x,p_1)-G(\bar q_{2}^k,p_1))\right)\\
			&+\frac{\rho_k}{3+\alpha+\beta}\left((1+\beta)(G(x,p_2)-G(\bar q_{2}^k,p_2))\right).
		\end{aligned}
	\end{equation}
	%and then let us set
	%$$
	%h_{k,*}(x)=\rho_k h(x)e^{\phi-4\pi\alpha G(x,p_1)-4\pi\beta G(x,p_2)}.
	%$$
	
	\begin{rem} {\it
			In the general case the corresponding definitions around a regular blow up point $q_j$ would be,
			\begin{align*}
				f_{k,j}(x)=~&\frac{\rho_k}{m+\sum\limits_{j=1}^\tau\alpha_j}
				\sum\limits_{\ell\neq j, \ell=\tau+1}^m (G(x,\bar q_{2,\ell}^k)-G(\bar q_{2,j}^k,\bar q_{2,\ell}^k))\\
				&+\frac{\rho_k}{m+\sum\limits_{j=1}^\tau\alpha_j}\sum\limits_{i=1}^\tau(1+\alpha_i)(G(x,p_i)-G(\bar q_{2,j}^k,p_i))\\
				&+\frac{\rho_k}{m+\sum\limits_{j=1}^\tau\alpha_j}\left(R(x,\bar q_{2,j}^k)-R(\bar q_{2,j}^k,\bar q_{2,j}^k)\right)
		\end{align*}}
		%$$
		%h_{k,*}(x)=\rho_k h(x)e^{\scriptscriptstyle \phi-\sum\limits_{i=1}^\tau 4\pi\alpha_i G(x,p_i)},
		%$$
	\end{rem}
	
	%The proof makes a crucial use of the following Pohozaev-type identity recently obtained in \cite{lin-yan-uniq}. Compared to \cite{bart-4}
	%we have some additional terms due to the fact that Green functions are not harmonic far away from the singular sources. 
	
	Because of the refined estimates we have obtained, we can use an easier version of the Pohozaev identity (as opposed to the argument in \cite{lin-yan-uniq,bart-4} etc). We claim that, with the choice of $f_k$ as in \eqref{exp-f-3}, the oscillation of 
	$u_2^k$ on $\partial B(q_2^k,\tau/2)$ is negligible. 
	Indeed, from the Green representation of $\rm w_i^k$ we have 
	\begin{align*}
		{\rm w_i^k}(x)=~&\overline{ \rm w_i^k}+\int_{\Omega(p_1,\tau)}G(x,\eta)\rho_k He^{\rm w_i^k}+\int_{\Omega(p_2,\tau)}G(x,\eta)\rho_k He^{\rm w_i^k}\\
		&+\int_{\Omega(q,\tau)}G(x,\eta)\rho_k He^{\rm w_i^k}+O(e^{-\lambda_1^k}).   
	\end{align*} 
	By using our refined estimates to evaluate each integral we have 
	$${\rm w_i^k}(x)=\overline{ \rm w_i^k}+8\pi\left((1+\alpha)G(x,p_1)+ (1+\beta)G(x,p_2)+  G(x,\bar q_2^k)\right)+O(\varepsilon_k^2). $$
	Thus, as far as $x$ is far away from blow up points, we have,
	\begin{align*}
		u_i^k(x)=~&\overline{\rm w_i^k} +8\pi (1+\alpha)G(x,p_1)+8\pi(1+\beta)G(x,p_2)\\
		&+8\pi G(x,\bar q_2^k)- f_k+O(\varepsilon_k^2).
	\end{align*}
	By the expression of $f_k$ in (\ref{exp-f-3}), we see that if we use $\phi_2^k$ to denote the local harmonic function defined by the oscillation of $u_2^k$ on $\partial B_{\tau/2}(\bar q_2^k)$, we have
	$$|\phi_2^k(\bar \varepsilon_ky)|\le C\bar \varepsilon_k\varepsilon_k^2|y|.$$ 
	Here we used $\rho_*-\rho_k=O(\varepsilon_k^2)$, which follows from (\ref{small-rho-k}).
	
	Then for $j=1,2$, we use the standard Pohozaev identity:
	\begin{equation}
		\label{new-pi-1}
		\begin{aligned}
			&\int_{\partial \Omega_k}(-\partial_{\nu}\bar v_1^k\partial_j\bar v_1^k+\frac 12|\nabla \bar v_1^k|^2\gamma_j)\\
			&=\int_{\partial \Omega_k}e^{v_1^k}\bar h_{k,0}(\bar q_2^k+\bar \varepsilon_ky)\gamma_jdS-\bar \varepsilon_k\int_{\Omega_k}\partial_j\bar h_{k,0}(\bar q_2^k+\bar \varepsilon_ky)e^{\bar v_1^k}.  
		\end{aligned}
	\end{equation}
	\begin{equation}
		\label{new-pi-2}
		\begin{aligned}
			&\int_{\partial \Omega_k}(-\partial_{\nu}\bar v_2^k\partial_j\bar v_2^k+\frac 12|\nabla \bar v_2^k|^2\gamma_j)\\
			&=\int_{\partial \Omega_k}e^{\bar v_2^k}\bar h_{k,0}(\bar q_2^k+\bar \varepsilon_ky)\gamma_jdS-\bar\varepsilon_k\int_{\Omega_k}\partial_j\bar h_{k,0}(\bar q_2^k+\bar \varepsilon_ky)e^{\bar v_2^k}. 
		\end{aligned}
	\end{equation}
	Since by definition,
	$$\bar v_1^k(y)=\bar v_2^k(y)+\sigma_k\bar \xi_k(y), \quad y\in \Omega_k, $$
	then from (\ref{new-pi-1}) and (\ref{new-pi-2}) we have
	\begin{equation}
		\label{impor-2}
		\begin{aligned}
			&\frac{1}{\bar \varepsilon_k}\int_{\partial \Omega_k}(-\partial_{\nu}\bar v_2^k\partial_j\bar \xi_k-\partial_{\nu}\bar \xi_k\partial_j\bar v_2^k+(\nabla \bar v_2^k\cdot \nabla \bar \xi_k)\gamma_j)dS\\
			&=-\int_{\Omega_k}\partial_j\bar h_{k,0}(\bar q_2^k+\bar\varepsilon_ky)e^{\bar v_2^k}\bar \xi_k+O(\bar \varepsilon_k^{1+\delta_0}),\quad j=1,2. 
		\end{aligned}
	\end{equation}
	Concerning $\nabla \bar v_2^k$ we use 
	\begin{equation*}
		\nabla \bar v_2^k(y)=-(4+O(\lambda_1^ke^{-\lambda_1^k}))\left(\frac{y_1}{r^2},\frac{y_2}{r^2}\right),\quad |y|\sim \bar \varepsilon_k^{-1}. 
	\end{equation*}
	Then (\ref{impor-2}) is reduced to 
	\begin{equation}\label{pi-final}
		\frac{1}{\bar \varepsilon_k}\int_{\partial \Omega_k}\frac{4}{r}\partial_j\bar \xi_k=-\int_{\Omega_k}\partial_j\bar h_{k,0}(\bar q_2^k+\bar\varepsilon_ky)e^{\bar v_2^k}\bar \xi_k+O(\bar \varepsilon_k^{1+\epsilon_0}). 
	\end{equation}
	To evaluate the left hand side of (\ref{pi-final})  we use
	(\ref{xi-k-large}) to deduce that,
	$$\frac{1}{\bar \varepsilon_k}\int_{\partial \Omega_k}\frac{4}{r}\partial_j\bar \xi_k=O(\bar \varepsilon_k^{1+\epsilon_0}).$$
	To evaluate the right hand side of (\ref{pi-final}), we use,
	$$\partial_j\bar h_{k,0}(\bar q_2^k+\bar \varepsilon_ky)=\partial_j\bar h_{k,0}(\bar q_2^k)+\sum_{s=1}^2\partial_{js}\bar h_{k,0}(\bar q_2^k)\bar \varepsilon_ky_s+O(\bar \varepsilon_k^2)|y|^2,$$
	$$e^{\bar v_2^k}=\frac{1}{(1+a_k|y|^2)^2}(1+O((\log \bar \varepsilon_k)^2\bar \varepsilon_k^2))$$
	and (\ref{bar-xi-e}) for $\bar \xi_k$,
	(\ref{estc0xi}) for $c_{0,\xi}^k$,
	(\ref{decay-w1}) for $\bar w_{1,\xi}^k$. Then by a straightforward evaluation we have, 
	$$\bar \varepsilon_k\left(\begin{matrix}
		\partial_{11}\bar h_{k,0}(\bar q_2^k) & \partial_{12}\bar h_{k,0}(\bar q_2^k) \\
		\partial_{12}\bar  h_{k,0}(\bar q_2^k) & \partial_{22}\bar h_{k,0}(\bar q_2^k)
	\end{matrix}
	\right)\left(\begin{matrix}
		b_1^k\\
		b_2^k
	\end{matrix}
	\right)=O(\bar \varepsilon_k^{1+\epsilon_0}).$$
	Since the underlying assumption is $(b_1^k,b_2^k)\to (b_1,b_2)\neq (0,0)$ 
	and since the non-degeneracy condition is ${\rm det}(D^2\bar h_k(0))\neq 0$, then we readily obtain a contradiction. Theorem \ref{main-theorem-2} is established under the assumptions that $\alpha>1$ and there are three blowup points consisting of one positive singular source, one negative singular source and one regular point. 
	The case where $\tau=m$ is easier and we omit the details to avoid repetitions.
	The more general case is worked out by minor changes, essentially just about notations, which proves Theorem \ref{main-theorem-2} as far as $\alpha>1$.
	
	\medskip
	
	On the other side, whenever $0<\alpha<1$, since
	$$\varepsilon_k^2=e^{-\frac{\lambda_1^k}{1+\alpha}}=o(e^{-\lambda_1^k/2}), $$
	then the proof of $b_1=b_2=0$ around each regular blowup point follows as in \cite{bart-4}.
	Theorem \ref{main-theorem-2} is established under the assumption that the set of blowup points
	consists of one positive singular source, one negative singular source and one regular point. 
	The case where $\tau=m$ is easier and we omit the details to avoid repetitions. At last, as above, 
	the more general case is worked out by minor changes, essentially just about notations, which proves Theorem \ref{main-theorem-2}. $\hfill\Box$

	\section{Proof of Other Theorems}
	
	\begin{proof}[Proof of Theorem \ref{mainly-case-2}.]
	We first consider the case $\alpha_M=0$, so that 
	the set of blowup points consists either of both negative singular sources and regular points, or just of regular points. As far as all blowup points are regular points, we refer the reader to \cite{bart-4} for the proof. Here we mention the two situations pursued in \cite{bart-4}:
	\begin{enumerate}\item
		$L(\mathbf{p})\neq 0$, ${\rm det}(D^2f^*)(p_{\tau+1},...,p_m))\neq 0$
		\item $L(\mathbf{p})=0$, ${\rm det}(D^2f^*)(p_{\tau+1},...,p_m))\neq 0$, $D(\mathbf{p})\neq 0$.
	\end{enumerate}
	
	A minor difference with respect to the argument in \cite{bart-4} is about the proof of $b_0=0$.
	To simplify the exposition we just use one negative singular pole $p$ with a negative index $-1<\beta<0$ and a regular blowup point $q$. Then we have,
	$$
	\rho_k=\int_M\rho_kHe^{u_i^k}=\int_{B(p,\tau)}\rho_kHe^{u_i^k}+\int_{B(q,\tau)}\rho_kHe^{u_i^k}+\int_{M\setminus (B(p,\tau)\cup B(q,\tau))}\rho_kHe^{u_i^k},
	$$
	and our local expansions imply, as far as either $L(\mathbf{q})\neq 0$ or $L(\mathbf{q})=0$ but $D(\mathbf{q})\neq 0$, that 
	$|\lambda_1^k-\lambda_2^k|=O(e^{-(1-\epsilon_0)\lambda_1^k})$ for some 
	$\epsilon_0>0$. At this point, if $b_0\neq 0$, the key observation is that the integral around $p$ satisfies,
	$$A_{p}^k:=\int_{B(p,\tau)}\rho_kc_k(x)\xi_k=b_0 c(\tau)e^{-\lambda_1^k}+o(e^{-\lambda_1^k})$$
	where $c(\tau)\to 0$ as $\tau\to 0$, which comes from the local integration around a negative singular source. This means in either one of the cases above ($L(\mathbf{q})\neq 0$ or $L(\mathbf{q})=0,D(\mathbf{q})\neq 0$) that $A_p^k$  is not the leading term. As a consequence the whole argument in \cite{bart-4} can be applied exactly as it stands, 
	because $A_p^k$ does not play any role.  On then other side, as far as there are no regular blowup points, then all the blowup points are negative singular sources and a careful evaluation based on our local estimates shows that the leading term is the one proportional to $D(\mathbf{q})\neq 0$, which comes, see \cite{bart-4} for details, as a result of a global integral. In other words, local bubbles contributions are minor compared to that proportional to $D(\mathbf{q})$, which is why the proof of $b_0=0$ follows again by the same argument in \cite{bart-4}. Theorem \ref{mainly-case-2} is established. \end{proof}
	
	\medskip

	Finally, we prove the theorems about Dirichlet problems.
	
	\begin{proof}[Proof of Theorem \ref{main-theorem-4} and Theorem \ref{main-theorem-1}.]
		For the blowup solutions to (\ref{equ-flat}), the corresponding estimates as in section \ref{preliminary} have been also obtained in \cite{chen-lin,zhang2} for $\alpha_j\in\mathbb{R}^+\setminus\mathbb{N}$ and in \cite{chen-lin-sharp,zhang1,gluck} for $\alpha_j=0$. Those preliminary estimates have almost the same form except for $\phi_j=0$ and $K\equiv 0$, where $\phi_j$ are the conformal factor at $p_j$ and $K$ is the Gaussian curvature of $M$.
		
		Then, under the assumption of regularity about $\partial\Omega$ and $q_j\in \Omega$ $(1\leq j\leq N)$, it has been shown in \cite{ma-wei} via the moving planes method that the blowup points of (\ref{equ-flat}) are far away from $\partial\Omega$. On the other hand, the crucial estimates obtained in section \ref{difference} and \ref{pf-uni-2} are in fact just local estimates. Consequently, a careful inspection of the proof shows that the new facts pushed forward here are not at all affected by the presence of the boundary of $\Omega$ throughout the argument. As a consequence Theorem \ref{main-theorem-4} and Theorem \ref{main-theorem-1} can be proved by the same argument.
	\end{proof}

    \section{Appendix}
    In this appendix we prove a useful lemma. The notations in this section overlap with those in other sections, so the readers should treat these notations as being independent of those in other sections. 
    \begin{lem}\label{osc-1}
Let $u_k$ satisfy 
\[\Delta u_k+h_ke^{u_k}=0,\quad \mbox{in}\quad B(0,\tau)\]
where $h_k$ is $C^5$, bounded above and below and there is a uniform bound on $\int_{B_{\tau}}h_ke^{u_k}$. Suppose $0$ is the only blowup point of $u_k$ and $u_k$ is constant on $\partial B_{\tau}$. Then we assume that $q_k=O(\epsilon_k^2)$ is the local maximum of $u_k$, and 
$\epsilon_k=e^{-\frac{u_k(q_k)}2}$. We also require $\nabla h_k(q_k)=O(\epsilon_k^2\log \frac{1}{\epsilon_k})$. Let $B(q_k,\tilde \tau_k)$ be inside $B(0,\tau_k)$ and these two balls are tangent. Then $\tau-\tilde \tau_k=|q_k|$. Then we claim that 
\begin{align*}
\frac{1}{2\pi}\int_0^{2\pi}u_k(q_k+\tilde \tau_ke^{i\theta})\cos \theta d\theta=-\frac{2q_1^k}{\tau}+O(\epsilon_k^3),\\
\frac{1}{2\pi}\int_0^{2\pi}u_k(q_k+\tilde \tau_ke^{i\theta}) \sin \theta d\theta=-\frac{2q_2^k}{\tau}+O(\epsilon_k^3),
\end{align*}
where $q_k=(q_1^k,q_2^k)$.
\end{lem}
\begin{rem}\label{omega-1-oscillation} If we use  $\phi_E^k$ to remove the oscillation of $u_k$ on $\partial B(q_k,\tilde \tau_k)$ with $\phi_E^k(q_k)=0$, we have
\[\phi_{E}^k(q_k+\epsilon_ky)=-\frac{2q_1^k}{\tau}\epsilon_ky_1-\frac{2q_2^k}{\tau}\epsilon_ky_2+O(\epsilon_k^3)|y|+\phi_{E,2}(\epsilon_ky)   \]
where $\phi_{E,2}^k$ is orthogonal to $e^{i\theta}$. By adding the estimates on the projections of $e^{ik\theta}$ ($k\neq 1$) we have
\[|\phi_{E,2}^k(\epsilon_ky)|\le C\epsilon_k^4|y|^2.\]
\end{rem}

    \begin{proof}[Proof of Lemma \ref{osc-1}.] 
    First we remark that the oscillation of $u_k$ on $\partial B(q_k,\tilde \tau_k)$ is $O(\epsilon_k^2)$, because for any two points $x_1,x_2$ on $\partial B(q_k,\tilde \tau_k)$, there exists $x_1^k,x_2^k$ on $\partial B(0,\tau)$ such that $|x_1-x_1^k|+|x_2-x_2^k|\le C\epsilon_K^2$. Since $|\nabla u_k|\sim 1$ near $\partial B(0,\tau)$ and $u_k(x_1^k)=u_k(x_2^k)$, we see that 
    \[|u_k(x_1)-u_k(x_2)|=|u_k(x_1)-u_k(x_1^k)-(u_k(x_2)-u_k(x_2^k)|\le C\epsilon_K^2.\]

    Let 
    \[v_k(y)=u_k(\epsilon_ky)+2\log \epsilon_k,\quad y\in \Omega_k:=B(0,\tau\epsilon_k^{-1}),\]
    and $\Gamma_k$ be the image of $\partial B(q_k,\tilde \tau_k)$ after scaling, we first observe that the oscillation of $v_k$ on $\partial \Omega_{1,k}$ is $O(\epsilon_k^2)$. In order to identify the leading term in the projection to $e^{i\theta}$,  we use $G_k$ to denote the Green's function on $\Omega_k$:
  \begin{equation}\label{green-ap}
G_k(y,\eta)=-\frac 1{2\pi}\log |y-\eta |+\frac 1{2\pi}\log \frac{|y|}{L_k}+\frac 1{2\pi}\log |y^*-\eta |
\end{equation}
where $L_k=\tau\epsilon_k^{-1}$, $y^*=\frac{L_k^2y}{|y|^2}$. Now the image of $q_k$ becomes $q_k/\epsilon_k$.
For $y\in \Gamma_k$ the Green representation formula gives
\begin{equation}\label{tem-8}
v_k(y)=\int_{\Omega_k}G_k(y,\eta)\tilde h_ke^{v_k}d\eta-\int_{\partial \Omega_k}\partial_{\nu}G_k v_kdS,
\end{equation} 

Our goal is to prove 
\begin{equation}\label{goal-ap}
\frac{1}{2\pi}\int_0^{2\pi}v_k(\frac{q_k}{\epsilon_k}+re^{i\theta})\cos \theta d\theta=2\frac{q_1^k}{\tau}+O(\epsilon_k^3),\quad r=\tilde \tau\epsilon_k^{-1},
\end{equation}
and the corresponding equation with $\sin \theta$. 

 If we use (\ref{tem-8}) as the expression for $v_k$ to evaluate (\ref{goal-ap}), we observe that the last term in (\ref{tem-8}) is minor, because it is a constant. 
 Next we deal with the term $\frac 1{2\pi}\log \frac{|y|}{L_k}$, which corresponds to the following integral:
 \begin{equation}\label{tem-12}\frac{1}{2\pi}\int_0^{2\pi}\frac 1{4\pi}\log \frac{|y|^2}{L_k^2}\cos \theta d\theta \int_{\Omega_k}\tilde h_ke^{v_k}
 \end{equation}
 Here $y=q_k/\epsilon_k+\tilde \tau\epsilon_k^{-1}e^{i\theta}$.
 Direct computation gives
 \[\frac 1{4\pi}\log \frac{|y|^2}{L_k^2}=\frac 1{4\pi}\log \frac{|q_k|^2+\tilde \tau^2+2q_1^k\tilde \tau\cos\theta+2q_2^k\tilde \tau\sin\theta}{\tau^2}.\]
 Using this expression in the evaluation of (\ref{tem-12}) we have
 \begin{equation}\label{tem-13}
 \begin{aligned}
\frac{1}{2\pi}\int_0^{2\pi}\frac 1{4\pi}\log \frac{|y|^2}{L_k^2}\cos \theta d\theta \int_{\Omega_k}\tilde h_ke^{v_k}
 =~&\frac{1}{4\pi}\frac{q_1^k\tilde \tau}{\tau^2}\int_{\Omega_k}\tilde h_ke^{v_k} \\
 =~&2\frac{q_1^k}{\tau}+O(\epsilon_k^3).  
 \end{aligned}
 \end{equation}
 Next we evaluate 
 \begin{equation}\label{tem-14}\frac 1{2\pi}\int_0^{2\pi}\bigg (\int_{\Omega_k}\frac 1{2\pi}\log \frac{|y^*-\eta |}{|y-\eta |}\tilde h_ke^{v_k}d\eta \bigg )\cos \theta d\theta.
 \end{equation}
 Now we evaluate the integral by cutting $\Omega_k$ in several subsets. First, for each fixed $y\in \Gamma_k$ and $\eta\in B(y,\epsilon_k^{-1}/2)$, we use \[\frac 1{2\pi}\log \frac{|y^*-\eta |}{|y-\eta |}\le C\epsilon_k/|y^{**}-\eta |\] where we used $|y-y^*|\le C\epsilon_k$ and $y^{**}$ is between $y$ and $y^*$. Over this region $e^{v_k}\sim \epsilon_k^4$. So the integral of this part is $O(\epsilon_k^4)$. Next we consider for fixed $y\in \Gamma_k$, the integration over $\Omega_k\setminus (B(y,\epsilon_k^{-1}/2)\cup B(0,\epsilon_k^{-\frac 12})$. When $\eta$ is in this region, $|y^{**}-\eta |\sim \epsilon_k^{-1}$, then one sees that the integral over this region is $O(\epsilon_k^3)$. Finally we consider the integral over $B(0,\epsilon_k^{-\frac 12})$ when $y\in \Gamma_k$ is fixed. 
For this case we use

\begin{equation}
\label{tem-15}
\begin{aligned}
\frac 1{4\pi}\log \frac{|y^*-\eta |^2}{|y-\eta |^2}
=&\frac{1}{4\pi}\log (1+\frac{|y^*|^2-|y|^2+2(y^*-y)\cdot \eta}{|y-\eta |^2}) \\
=&\frac 1{4\pi}\log (1+\frac{|y^*|^2-|y|^2}{|y|^2}+O(\epsilon_k^3)|\eta |)\\
=&\frac 1{4\pi}\log \big(1+(\frac{|L_k|^4}{|y|^4}-1)+O(\epsilon_k^3)|\eta |\big ) 
\end{aligned}
\end{equation}
where in the second step we used 
\begin{equation}\label{tem-16}
\frac{1}{|y-\eta |^2}-\frac{1}{|y|^2}=O(\bar\varepsilon_k^3)|\eta |
\end{equation}
for $|\eta |<C\bar\varepsilon_k^{-\frac 12}$. Using (\ref{tem-16}) and 
\[|y|^2=\epsilon_k^{-1}\bigg (\tilde \tau^2+2q_1^k\tau\cos\theta+2q_2^k\tau \sin \theta +O(\epsilon_k^4)\bigg ),\]
we have
\[\frac 1{2\pi}\int_0^{2\pi}\bigg (\int_{B(0,\epsilon_k^{-\frac 12})}\frac 1{2\pi}\log\frac{|y^*-\eta |}{|y-\eta |}\tilde h_ke^{v_k}d\eta)\cos\theta d\theta=-4\frac{q_1^k}{\tau}+O(\epsilon_k^3).\]
The proof about the last estimate in the claim is similar and we skip it here to avoid repetitions.\\
Lemma \ref{osc-1} is established. 
\end{proof}
	\medskip
	
	On behalf of all authors, the corresponding author states that there is no conflict of interest.
	
	\medskip
	
	Data availability statements: This article is a complete mathematical proof that requires no data.


\begin{thebibliography}{99}
		
		\bibitem{ambjorn} J. Ambjorn, P. Olesen, Anti-screening of large magnetic fields by vector bosons. \emph{Phys. Lett. B}, \textbf{214} (1988), no. 4, 565--569.
		
		\bibitem{Aub} T. Aubin, "Nonlinear analysis on Manifolds Monge-Amp\'ere equations." Grundlehren der Mathematischen Wissenschaften {\bf 252},
		Springer-Verlag, New York, 1982.
		
		\bibitem{B1} D. Bartolucci, On the Best Pinching Constant of Conformal Metrics
		on $\mathbb{S}^2$ with One and Two Conical Singularities.
		\emph{The Journal of Geometric Analysis} 23  (2013), 855-877.
		
		\bibitem{BCLT} D. Bartolucci, C.C. Chen, C.S. Lin, G. Tarantello, Profile of blow-up solutions to mean field equations with singular data.
		\emph{Comm. P.D.E.} \textbf{29} (2004),  no. 7-8, 1241-1265.
		
		\bibitem{BdM1}
		D. Bartolucci, F. De Marchis, {On the Ambjorn-Olesen electroweak condensates},
		{\em  Jour. Math. Phys.} {\bf 53} 073704 (2012).
		
		\bibitem{BdM2} D. Bartolucci, F. De Marchis,
		{ Supercritical Mean Field Equations on convex domains and the Onsager's
			statistical description of two-dimensional turbulence}, {\em Arch. Rat. Mech. An.}, {\bf 217} (2015), 525-570.
		
		\bibitem{BdMM} D. Bartolucci, F. De Marchis, A. Malchiodi, {Supercritical conformal metrics on
			surfaces with conical singularities}, {\em Int. Math. Res. Not.}, (2011) {\bf (24)}, 5625-5643.
		
		\bibitem{BGJM} D. Bartolucci, C.Gui, A. Jevnikar, A. Moradifam,
		{A singular sphere covering inequality: uniqueness and symmetry of solutions to singular Liouville-type equations,}
		{\em Math. Ann.}  {\bf 374} (2019), 1883-1922.
		
		\bibitem{bart-5} D. Bartolucci, A. Jevnikar, Y. Lee, W. Yang, Non degeneracy, Mean Field Equations and the
		Onsager theory of 2D turbulence, \emph{Arch. Ration. Mech. An.}  \textbf{230} (2018), 397-426.
		
		\bibitem{bart-4} D. Bartolucci, A. Jevnikar, Y. Lee, W. Yang, Uniqueness of bubbling solutions of mean field equations. \emph{J. Math. Pures Appl.},  \textbf{123} (2019), 78-126.
		
		\bibitem{bart-4-2} D. Bartolucci, A. Jevnikar, Y. Lee, W. Yang, Local Uniqueness of Blowup solutions of mean field equations with singular data.
		\emph{Jour. Diff. Eqs.}, \textbf{269} (2020), 2057-2090.
		
		\bibitem{bart-4-3} D. Bartolucci, A. Jevnikar, Y. Lee, W. Yang, Local uniqueness of m-bubbling sequences for the Gel’fand equation. \emph{Comm. P.D.E.}, \textbf{44} (2019), 447-466.
		
		\bibitem{bj-lin}
		D. Bartolucci, A. Jevnikar, C.S. Lin, Non-degeneracy and uniqueness of solutions to singular mean field equations on
		bounded domains. {\em Jour. Diff. Eqs.} \textbf{266} (2019) 716–741, 
		
		\bibitem{bart-lin} D. Bartolucci, C.S. Lin, Uniqueness Results for Mean Field Equations with Singular Data.  \emph{Comm. P. D. E.}, \textbf{34}(7) (2009), 676-702.
		
		\bibitem{BLin2} D. Bartolucci, C.S. Lin, { Sharp existence results for mean field equations with singular data}. \emph{Jour. Diff. Eq.} \textbf{252} (2012), 4115-4137.
		
		\bibitem{BLin3} D. Bartolucci, C.S. Lin, { Existence and uniqueness for
			Mean Field Equations on multiply connected domains at the critical parameter}.
		{\em Math. Ann.}, {\bf 359} (2014), 1-44; DOI 10.1007/s00208-013-0990-6.
		
		\bibitem{BLT} D. Bartolucci, C.S. Lin, G. Tarantello, { Uniqueness and symmetry results for
			solutions of a mean field equation on ${\mathbb{S}}^{2}$ via a new bubbling phenomenon}.
		{\em Comm. Pure Appl. Math.} {\bf 64}(12) (2011), 1677-1730.
		
		\bibitem{BMal} D. Bartolucci, A. Malchiodi, {An improved geometric
			inequality via vanishing moments, with applications to singular
			Liouville equations}, {\em Comm. Math. Phys.} {\bf 322} (2013), 415-452.
		
		
		\bibitem{BM3} D. Bartolucci, E. Montefusco, { Blow up analysis,
			existence and qualitative properties of solutions for the two
			dimensional Emden-Fowler equation with singular potential}, {\em M$^{2}$.A.S. } 
		{\bf 30}(18) (2007), 2309-2327.
		
		\bibitem{BT} D. Bartolucci, G. Tarantello, Liouville type
		equations with singular data and their applications to periodic
		multivortices for the electroweak theory. {\em Comm. Math. Phys.} 229
		(2002), 3-47.
		
		\bibitem{BT-2}  D. Bartolucci; G. Tarantello, Asymptotic blow-up analysis for singular Liouville type equations with applications. \emph{J. Differential Equations}, \textbf{262} (2017), 3887-3931.
		
		\bibitem{luca-b} L. Battaglia, M. Grossi, A. Pistoia, Non-uniqueness of blowing-up solutions to the Gelfand problem. {\em Calc. Var. \& P.D.E.} \textbf{58} (2019), art. 163
		
		\bibitem{BM} H. Brezis, F. Merle, Uniform estimates and blow-up behavior
		for solutions of $-\Delta u=v(x)e^u$ in two dimensions. {\em Comm. P. D. E.}, 16 (1991)
		1223-1253.
		
		%\bibitem{CaYang} L. Caffarelli, Y. Yang, Vortex condensation in the Chern-Simons Higgs model: an existence theorem.  Comm. Math. Phys.  168  (1995),  no. 2, 321--336.
		
		%\bibitem{caglioti-1} E. Caglioti, P.L. Lions, C. Marchioro , M. Pulvirenti, A special class of stationary flows for two-dimensional Euler equations: A statistical mechanics description, \emph{Comm. Math. Phys.}, \textbf{143} (1992), 501--525.
		
		\bibitem{caglioti-2} E. Caglioti, P.L. Lions, C. Marchioro , M. Pulvirenti, A special class of stationary flows for two-dimensional Euler equations: A statistical mechanics description, part II. \emph{Comm. Math. Phys.}, \textbf{174} (1995), 229--260.
		
		\bibitem{chai} C.C. Chai, C.S. Lin, C.L.Wang, Mean field equations, hyperelliptic curves, and modular forms: I, \emph{Camb. J. Math.}. \textbf{3}(1-2) (2015), 127-274.
		
		\bibitem{chan-fu-lin} H. Chan, C.C. Fu, C.S. Lin, Non-topological multi-vortex solutions to the self-dual Chern-Simons-Higgs equation, \emph{Comm. Math. Phys.}, \textbf{231} (2002), no. 2, 189-221.
		
		\bibitem{chang-chen-lin} A. Chang, C.C. Chen, C.S. Lin, Extremal functions for a mean field equation in two dimension. New Stud. Adv. Math., 2
		International Press, Somerville, MA, 2003, 61–93.
		
		\bibitem{CL1} W.X. Chen, C.M. Li, { Classification of solutions of some nonlinear elliptic equations,}
		{\em Duke Math. J.}  {\bf 63}(3) (1991), 615-622.
		
		\bibitem{CL2} W.X. Chen, C.M. Li, { Qualitative properties of solutions of
			some nonlinear elliptic equations in $R^2$}, {\em Duke Math. J.} 71(2) (1993), 427-439.
		
		\bibitem{chen-lin-sharp} C.C. Chen, C.S. Lin, Sharp estimates for solutions of multi-bubbles in compact Riemann surface. \emph{Comm. Pure Appl. Math.}, \textbf{55} (2002), 728-771.
		
		\bibitem{chen-lin-deg} C.C. Chen, C.S. Lin, Topological degree for a mean field equation on Riemann surfaces. \emph{Comm. Pure Appl. Math.}, \textbf{56} (2003), 1667-1727.
		
		\bibitem{chen-lin-wang} C.C. Chen, C.S. Lin, G.Wang, Concentration phenomena of two-vortex solutions in a Chern-Simons model. \emph{Ann. Sc. Norm. Super. Pisa Cl. Sci.}, (5) \textbf{3} (2004), 2, 367397.
		
		\bibitem{chen-lin} C.C. Chen, C.S. Lin, Mean field equation of Liouville type with singularity data: Sharper estimates, \emph{Discrete and Continuous Dynamic Systems-A}, \textbf{28} (2010), 1237-1272
		
		\bibitem{chen-lin-deg-2} C.C. Chen, C.S. Lin, Mean field equation of Liouville type with singular data: topological degree. \emph{Comm. Pure Appl. Math.}, \textbf{68} (2015), 6, 887-947.
		
		\bibitem{Zchen-lin-1} 
		Z. Chen, C.S. Lin, Critical points of the classical Eisenstein
		series of weight two. \emph{J. Differential Geom.} 
		\textbf{113} (2019), 189-226. 
		
		\bibitem{Zchen-lin} 
		Z. Chen, C.S. Lin, Sharp nonexistence results for curvature
		equations with four singular sources on rectangular tori. \emph{Amer. J. Math.} 
		\textbf{142} (2020), 1269-1300. 
		
		%\bibitem{chen-kuo-lin} Z.J. Chen, T.J. Kuo, C.S. Lin, Hamiltonian system for the elliptic form of Painleve VI equation, \emph{J. Math. Pure App.}, \textbf{106}(3) (2016), 546-581.
		
		%\bibitem{coddington} E.A. Coddington, N. Levinson, Theory of ordinary differential equations. McGraw-Hill Book Company, Inc., New York-Toronto-London, 1995.
		
		\bibitem{dAWZ} T. D'Aprile, J. Wei, L. Zhang,  {On the construction of non-simple blow-up solutions for the singular Liouville equation with a potential}, {\em Calc. Var. \& P.D.E.} to appear.
		
		\bibitem{DKM} M. del Pino, M. Kowalczyk, M. Musso,  {Singular limits in
			Liouville-type equations}, {\em Calc. Var. \& P.D.E.} {\bf 24} (2005), 47-81.
		
		\bibitem{DJ} Z. Djadli, {Existence result for the mean field problem on Riemann surfaces of all genuses},{\em Comm. Contemp. Math.} \textbf{10} (2008), 205-220.
		
		\bibitem{EGP} P. Esposito, M. Grossi \& A. Pistoia, {On the existence of blowing-up solutions for a mean field equation}, {\em Ann. Inst. H. Poincare Anal. Nonlinear }, {\bf 22} (2005), 227-257.
		
		%\bibitem{Gladiali} F. Gladiali, M. Grossi, Some results for the Gelfand's problem. \emph{Comm. Partial Differential Equations}, 29 (2004), no. 9-10, 1335–1364.
		
		\bibitem{gu-zhang-1} Y. Gu, L. Zhang, Degree counting theorems for singular Liouville systems. {\em Ann.
			Sc. Norm. Super. Pisa Cl. Sci.} (5) Vol. XXI (2020), 1103-1135.
		
		\bibitem{gu-zhang-2} Y. Gu, L. Zhang, Structure of bubbling solutions of Liouville systems with negative singular
		sources, preprint 2022. https://arxiv.org/abs/2112.10031.
		
		\bibitem{fang-lai} H. Fang, M. Lai, On curvature pinching of conic 2-spheres, \emph{Calc. Var.  P.D.E.}, \textbf{55} (2016), 118.
		
		\bibitem{gluck} M. Gluck, Asymptotic behavior of blow up solutions to a class of prescribing Gauss curvature equations. \emph{Nonlinear Anal.}, \textbf{75} (2012), 5787-5796.
		
		\bibitem{GM1} C. Gui, A. Moradifam, { The Sphere Covering Inequality and Its Applications}, {\em Invent. Math.} \textbf{214, (2018) 1169-1204}.
		
		\bibitem{huang-zhang-cvpde} H.Y. Huang, L. Zhang,  On Liouville systems at critical parameters, Part 2: Multiple bubbles. {\em Calc. Var. Partial Differential Equations}, 61 (2022), no. 1, Paper No. 3, 26 pp.
		
		\bibitem{KW} J. L. Kazdan, F. W. Warner,
		{Curvature functions for compact 2-manifolds}, {\em Ann. Math.} {\bf 99}  (1974), 14-74.
		
		\bibitem{kuo-lin} T.J. Kuo, C.S. Lin, Estimates of the mean field equations with integer singular sources: non-simple blow up, \emph{Jour. Diff. Geom.}, \textbf{103} (2016), 377-424.
		
		\bibitem{lee-lin-jfa} Y. Lee, C.S. Lin, Uniqueness of bubbling solutions with collapsing singularities. \emph{J. Funct. Anal.}, 277 (2019), no. 2, 522–557.
		
		\bibitem{li-cmp} Y.Y. Li, Harnack type inequality: the method of moving planes, \emph{Comm. Math. Phys.}, \textbf{200} (1999), 421--444.
		
		\bibitem{ls} Y.Y. Li, I.Shafrir, {Blow-up analysis for Solutions of $-\Delta u = V(x)e^{u}$
			in dimension two}, {\em Ind. Univ. Math. J.}  {\bf 43} (1994), 1255-1270.
		
		
		\bibitem{Lin1} C.S. Lin, {Uniqueness of solutions to the mean field equation for the
			spherical Onsager Vortex}, \emph{Arch. Rat. Mech. An.} {\bf 153} (2000), 153-176.
		
		
		\bibitem{lin22}  C.S. Lin, Spherical metrics with one singularity and odd integer angle on flat tori I, preprint 2022.
		
		\bibitem{lin-yan-uniq}  C.S. Lin, S.S. Yan, On the mean field type bubbling solutions for Chern-Simons-Higgs equation. \emph{Adv. Math.}, \textbf{338} (2018), 1141--1188.
		
		\bibitem{lin-yan-cs} C.S. Lin, S.S. Yan,
		Existence of bubbling solutions for Chern-Simons model on a torus. {\em Arch. Ration. Mech. An.} \textbf{207} (2013), 353–392.
		
		\bibitem{lin-Lwang} C.S. Lin, 
		C.L. Wang. Elliptic functions, Green functions
		and the mean field equations on tori.  \emph{Ann. Math.}, \textbf{172} (2010), 911-954.
		
		\bibitem{lin-wang-cga} C.S. Lin, 
		C.L. Wang. Geometric quantities arising from bubbling analysis of mean field equations.  \emph{Comm An. Geom.}, \textbf{28} (2020), 1289-1313.
		
		\bibitem{lin-zhang-jfa} C.S. Lin, L. Zhang, On Liouville systems at critical parameters, Part 1: One bubble
		{\em J. Funct. Anal.} 264 (2013), no. 11, 2584–2636.
		
		\bibitem{ma-wei} L. Ma, J.C. Wei, Convergence for a Liouville equation. \emph{Comment. Math. Helv.}, 76 (2001), no. 3, 506–514.
		
		\bibitem{Mal1} A. Malchiodi, {Topological methods for an elliptic equation with exponential nonlinearities}, \emph{Discr. Cont. Dyn. Syst.}
		{\bf 21} (2008), 277-294.
		
		\bibitem{Mal2} A. Malchiodi, {Morse theory and a scalar field equation on compact
			surfaces}, {\em Adv. Diff. Eq.} {\bf 13} (2008), 1109-1129.
		
		\bibitem{mal-ruiz} A. Malchiodi, D. Ruiz,
		New improved Moser-Trudinger inequalities and singular Liouville equations on compact surfaces.
		\emph{Geom. Funct. Anal.}, 21 (2011), no. 5, 1196–1217.
		
		\bibitem{nolasco-taran} M. Nolasco, G. Tarantello, On a sharp Sobolev-type Inequality on two-dimensional compact manifold, \emph{Arch. Rat. Mech. An.}, \textbf{145} (1998), 161-195.
		
		\bibitem{pot} A. Poliakovsky, G. Tarantello, {On a planar Liouville-type problem in the study of
			selfgravitating strings}, \emph{J. Differential Equations} {\textbf{252}} (2012), 3668-3693.
		
		\bibitem{Parjapat-Tarantello} J. Prajapat, G. Tarantello, On a class of elliptic problems in $\mathbb{R}^2$: Symmetry and uniqueness results, \emph{Proc. Roy. Soc. Edinburgh Sect. A}, \textbf{131} (2001), 967-985.
		
		\bibitem{SO} R.A. Smith, T.M. O'Neil, 
		nonaxisymmetric thermal equilibria of cylindrically bounded guiding-center plasmas or discrete vortex system, \emph{Phys. Fluids B}, \textbf{2} (1990), 2961–2975.
		
		\bibitem{spruck-yang} J. Spruck, Y. Yang, On Multivortices in the Electroweak Theory I:Existence of Periodic Solutions, \emph{Comm. Math. Phys.}, \textbf{144} (1992), 1-16.
		
		\bibitem{suzuki} T. Suzuki, Global analysis for a two-dimensional elliptic eiqenvalue problem with the exponential nonlinearly, \emph{Ann. Inst. H. Poincare Anal. Nonlinear }, \textbf{9}(4) (1992), 367-398.
		
		\bibitem{T3} G. Tarantello,
		{Analytical aspects of Liouville type equations with singular sources}, Handbook Diff. Eqs., North Holland,
		Amsterdam, Stationary partial differential equations, {\bf I} (2004), 491-592.
		
		\bibitem{taran-1} G. Tarantello, Multiple condensate solutions for the Chern-Simons-Higgs theory, \emph{J. Math. Phys.}, \textbf{37} (1996), 3769-3796.
		
		\bibitem{taran-2} G. Tarantello, Self-Dual Gauge Field Vortices: An Analytical Approach, PNLDE 72, Birkhauser Boston, Inc., Boston, MA, 2007.
		
		\bibitem{taran-3} G. Tarantello, Blow-up analysis for a cosmic strings equation, \emph{Jour. Funct. An.}, \textbf{272} (1) (2017) 255-338.
		
		\bibitem{taran-4} G. Tarantello, Asymptotics for minimizers of a Donaldson functional and mean curvature 1-immersions of surfaces into hyperbolic 3-manifolds, \emph{Adv. Math.}, \textbf{425} (2023) 109090.
		
		\bibitem{troy} M. Troyanov, Prescribing curvature on compact surfaces with conical singularities, \emph{Trans. Amer. Math. Soc.}, \textbf{324} (1991), 793-821.
		
		\bibitem{TuY} A. Tur, V. Yanovsky, Point vortices with a rational necklace: New exact stationary solutions
		of the two-dimensional Euler equation, \emph{Phys. Fluids }, \textbf{16} (2004), 2877–2885.
		
		\bibitem{wei-zhang-pacific} J.C. Wei, L. Zhang, Nondegeneracy of the Gauss curvature equation with negative conic singularity. \emph{Pacific J. Math.}, 297 (2018), no. 2, 455–475.
		
		\bibitem{wei-zhang-19} J.C. Wei, L. Zhang, Estimates for Liouville equation with quantized singularities, {\em Adv. Math.}, \textbf{380} (2021), 107606. 
		
		\bibitem{wei-zhang-22} J.C. Wei, L. Zhang, Vanishing estimates for Liouville equation with quantized singularities, {\em Proc. Lond. Math. Soc.}, \textbf{3} (2022), 1-26. 
		
		\bibitem{wei-zhang-jems} J.C.Wei, L. Zhang,
		Laplacian Vanishing Theorem for Quantized Singular Liouville Equation. To appear on Journal of European Mathematical Society.  
		
		\bibitem{wolan} G. Wolansky, On steady distributions of self-attracting clusters under friction and fluctuations, \emph{Arch. Rational Mech. An.}, \textbf{119} (1992), 355--391.
		
		\bibitem{wu-zhang-ccm} L.N. Wu, L. Zhang, Uniqueness of bubbling solutions of mean field equations with non-quantized singularities. Commun. Contemp. Math. 23 (2021), no. 4, Paper No. 2050038, 39 pp.
		
		\bibitem{y-yang} Y. Yang, Solitons in Field Theory and Nonlinear Analysis, Springer Monographs in Mathematics, Springer, New York, 2001.
		
		\bibitem{zhang1} L. Zhang, Blowup solutions of some nonlinear elliptic equations involving
		exponential nonlinearities, \emph{Comm. Math. Phys}, \textbf{268} (2006), 105-133.
		
		\bibitem{zhang2} L. Zhang, Asymptotic behavior of blowup solutions for elliptic equations with exponential nonlinearity and singular data, \emph{Commun. Contemp. Math}, \textbf{11} (2009), 395-411.
		
	\end{thebibliography}
\end{document}